\documentclass{amsart}

\usepackage{color}
\usepackage[dvips]{graphicx}
\usepackage{amscd}
\usepackage{amsmath}
\usepackage{ascmac}
\usepackage{cases}
\usepackage{amssymb}
\usepackage{amsfonts}
\usepackage{amsthm}
\usepackage[all]{xy}
\usepackage{enumerate}
\usepackage{wrapfig}

\theoremstyle{definition}
\newtheorem{Thm}{{\bf Theorem}}[section]

\newtheorem{Lem}[Thm]{{\bf Lemma}}
\newtheorem{Prop}[Thm]{{\bf Proposition}}
\newtheorem{Cor}[Thm]{{\bf Corollary}}

\newtheorem{Def}[Thm]{{\bf Definition}}
\newtheorem{Rem}[Thm]{{\bf Remark}}

\newcounter{Exami}

\newcommand{\balpha}{{\boldsymbol \alpha}}

\newcommand{\bnu}{{\boldsymbol \nu}}

\newcommand{\bt}{{\bold t}}

\newcommand{\bH}{{\mathbf H}}

\newcommand{\cA}{{\mathcal A}}
\newcommand{\cG}{{\mathcal G}}
\newcommand{\cO}{{\mathcal O}}

\newcommand{\cH}{{\mathcal H}}
\newcommand{\cM}{{\mathcal M}}

\newcommand{\cD}{{\mathcal D}}

\newcommand{\cC}{{\mathcal C}}

\newcommand{\cF}{{\mathcal F}}
\newcommand{\cE}{{\mathcal E}}

\newcommand{\cL}{{\mathcal L}}

\newcommand\rank{\mathop{\rm rank}\nolimits}

\newcommand{\Tr}{\mathop{\rm Tr}\nolimits}

\newcommand\Spec{\mathop{\rm Spec}\nolimits}

\newcommand{\res}{\mathop{\sf res}\nolimits}

\newcommand\lra{\longrightarrow}
\newcommand\ra{\rightarrow}

\title[Hamiltonian structures of isomonodromic deformations]{Hamiltonian structures 
of isomonodromic deformations on moduli spaces of parabolic connections}
\author{Arata Komyo}
\address{Center for Mathematical and Data Sciences, Kobe Univ.,
1-1 Rokkodai-cho, Nada-ku, Kobe, 657-8501, Japan}
\email{akomyo@math.kobe-u.ac.jp}

\begin{document}

\maketitle

\begin{abstract}
In this paper, we treat moduli spaces of parabolic connections.
We take an affine open covering of the moduli spaces,
and we construct a Hamiltonian structure of an algebraic vector field determined 
by the isomonodromic deformation for each affine open subset.
\end{abstract}

\section{Introduction}

Let $(C, \bt )$ ($\bt=(t_1, \ldots,t_n)$) be an $n$-pointed smooth projective curve of genus $g$ 
over $\mathbb{C}$, where $t_1, \ldots,t_n$ are distinct points.
We take a positive integer $r$ and
an element $\bnu=(\nu^{(i)}_j )^{1\le i \le n}_{0\le j \le r-1} \in \mathbb{C}^{nr}$ 
such that $\sum_{i,j} \nu^{(i)}_j = -d \in \mathbb{Z}$.
We say $(E,\nabla, \{ l^{(i)}_* \}_{1\le i \le n})$ is a $(\bt, \bnu)$-parabolic connection of rank $r$ if 
\begin{itemize}
\item[(1)] $E$ is a rank $r$ algebraic vector bundle of degree $d$ on $C$, 
\item[(2)] $\nabla\colon E \ra E\otimes \Omega^1_{C}(t_1+\cdots+t_n)$ is a connection, and
\item[(3)] for each $t_i$, $l^{(i)}_* $ is a filtration $E|_{t_i} = l_0^{(i)} \supset l_1^{(i)}
 \supset \cdots \supset l_r^{(i)}=0$ such that $\dim (l_j^{(i)}/l_{j+1}^{(i)})=1$
and $(\res_{t_i}(\nabla)-\nu^{(i)}_j \mathrm{id}_{E|_{t_i}}) (l^{(i)}_j) \subset l_{j+1}^{(i)}$ 
for $j=0,\ldots,r-1$.
\end{itemize}
Inaba--Iwasaki--Saito \cite{IIS} (for the general case see Inaba \cite{Inaba}) 
introduces the $\balpha$-stability for $(\bt, \bnu)$-parabolic connections,
and constructs the moduli scheme of $\balpha$-stable $(\bt, \bnu)$-parabolic connections 
of rank $r$, denoted by $M_C^{\balpha}(\bt, \bnu)$.
Moreover, let $T$ be a smooth algebraic scheme which is an \'etale covering of the moduli stack 
of $n$-pointed smooth projective curves of genus $g$ over $\mathbb{C}$ and
take a universal family $(\cC, \tilde{t}_1,\ldots, \tilde{t}_n)$ over $T$.
Let $N^{(n)}_r(d)$ be the set of $\bnu=(\nu^{(i)}_j ) \in \mathbb{C}^{nr}$ 
such that $\sum_{i,j} \nu^{(i)}_j = -d \in \mathbb{Z}$.
Then we can construct a relative fine moduli scheme $M_{\cC/T}^{\balpha}(\tilde{\bt},r,d) 
\ra  T \times N^{(n)}_r(d)$ of 
$\balpha$-stable parabolic connections of rank $r$ and of degree $d$, which is smooth 
and quasi-projective \cite[Theorem 2.1]{Inaba}. 
The moduli space $M_C^{\balpha}(\bt, \bnu)$ is a fiber of $M_{\cC/T}^{\balpha}(\tilde{\bt},r,d) 
\ra  T \times N^{(n)}_r(d)$
and is equipped with a natural symplectic structure.

The moduli space $M_{\cC/T}^{\balpha}(\tilde{\bt},r,d)$ gives a geometric description 
of the differential equation determined by the isomonodromic deformation.
We fix $\bnu$.
We can regard $M_{\cC/T}^{\balpha}(\tilde{\bt},r,d)_{\bnu}\ra T$ as a phase space 
of the differential equation determined by the isomonodromic deformation, 
and $T$ as a space of time variables.
A fiber of $M_{\cC/T}^{\balpha}(\tilde{\bt},r,d)_{\bnu} \ra  T$ becomes a space of initial conditions. 
In fact, for the case of $C=\mathbb{P}^1$, $r=2$ and $n=4$, these fibers coincide 
with the spaces of initial conditions for the Painlev\'e VI equation constructed 
by Okamoto \cite{Okamoto} 
(see \cite{IIS2}). 
Inaba--Iwasaki--Saito \cite{IIS} (for rank $2$ and $\mathbb{P}^1$ cases) 
and Inaba \cite{Inaba} (for general cases)
show that the Riemann-Hilbert correspondence induces 
a proper surjective bimeromorphic morphism between
the moduli space of $\balpha$-stable parabolic connections 
and the moduli space of certain equivalence classes of representations 
of the fundamental group $\pi_1(C\setminus\{ t_1,\ldots,t_n \}, *)$.
By this property of the Riemann-Hilbert correspondence, 
they show that the differential equation determined by the isomonodromic deformation 
satisfies the geometric Painlev\'e property (see \cite{IIS} and \cite{Inaba}).
Note that geometric descriptions of the isomonodromic deformation
are also given by Hitchin \cite{Hit} and Boalch \cite{Boal1}, \cite{Boal2} et al.\ 
from symplectic points of view.

One of the important properties of the Painlev\'e equations is that they can be written 
as (\textit{non-autonomous}) \textit{Hamiltonian systems}
(\cite{Okamoto2}, \cite{Okamoto3}, \cite{KO}).
The purpose of this paper is to give Hamiltonian descriptions of the vector fields determined 
by the isomonodromic deformations
 on the moduli space $M_{\cC/T}^{\balpha}(\tilde{\bt},r,d)_{\bnu}$, which is a phase space.
Hamiltonian descriptions of the vector fields determined by the isomonodromic deformations 
on moduli spaces of certain connections 
 were essentially considered by Krichever \cite{Krich} and Hurtubise \cite{Hurt}. 
 (Wong \cite{Wong} generalizes those results in the case of principal $G$-bundles).
We apply their ideas to the moduli space $M_{\cC/T}^{\balpha}(\tilde{\bt},r,d)_{\bnu}$ 
of parabolic connections.
Accordingly, we give an affine open covering $\{ M_i\}_i$ of 
$M_{\cC/T}^{\balpha}(\tilde{\bt},r,d)_{\bnu}$ 
and a Hamiltonian structure on each $M_i$.
Namely, we construct a 2-form $\omega$ on $M_{\cC/T}^{\balpha}(\tilde{\bt},r,d)_{\bnu}$ 
such that the kernel $\mathrm{Ker}(\omega)$ induces the vector fields determined 
by the isomonodromic deformations and $\omega$ is the symplectic form fiberwise 
(Proposition \ref{Kernel is IMD}).
The 2-form $\omega$ is considered in \cite{Iwa1} and \cite{Iwa2}.
We define Hamiltonian functions $H_i$ for $i=1,2,\ldots, \dim T$ on each $M_i$ 
(Definition \ref{Definition of Hamil}).
If we may take good coordinates on $M_{\cC/T}^{\balpha}(\tilde{\bt},r,d)_{\bnu}$,
we obtain a Hamiltonian description of the vector field on $M_i$ induced 
by the isomonodromic deformation (Corollary \ref{main theorem}).

In the argument of the Hamiltonian description of the isomonodromic deformations,
it is important to 
describe the vector fields determined by the isomonodromic deformations
on $M_{\cC/T}^{\balpha}(\tilde{\bt},r,d)_{\bnu}$ 
in terms of the hypercohomology of some complex.
Indeed, we give such a description of the isomonodromic deformations
on the moduli space $M_{\cC/T}^{\balpha}(\tilde{\bt},r,d)_{\bnu}$.
Biswas--Heu--Hurtubise \cite{BHH}
have described the vector fields determined by the isomonodromic deformations
in terms of the hypercohomology of some complex.
Essentially, our description is same as their description.
However, their description of the isomonodromic deformations are only for 
certain generic connections.
On the other hand, our moduli space $M_{\cC/T}^{\balpha}(\tilde{\bt},r,d)_{\bnu}$ 
is the \textit{full} phase space of 
the isomonodromic deformations, since 
the isomonodromic deformations on this moduli space have geometric Painlev\'e property.
So we give a description of the isomonodromic deformations
for the full phase space.

The organization of this paper is as follows.
In Section \ref{Pre}, we recall the description of the tangent space of $M_C^{\balpha}(\bt, \bnu)$ 
and of the natural symplectic structure on $M_C^{\balpha}(\bt, \bnu)$
in terms of the hypercohomology of a certain complex.
In Section \ref{SS AA}, we recall the \textit{Atiyah algebra}. 
In Section \ref{SS d t s of m by c c}, 
we discuss a description of the tangent space 
of $M_{\cC/T}^{\balpha}(\tilde{\bt},r,d)_{\bnu}$ in terms of the hypercohomology 
of a certain complex.
We use the Atiyah algebra in a definition of this complex.
In Section \ref{SS IMD by c c}, 
we describe the vector field determined by the isomonodromic deformation 
in terms of the hypercohomology.
In Section \ref{S Hamil}, we give Hamiltonian descriptions of the vector fields determined 
by the isomonodromic deformations. 
In Section \ref{SS const nabla0}, we take an affine open covering $\{ M_i\}_i$ of
$M_{\cC/T}^{\balpha}(\tilde{\bt},r,d)_{\bnu}$ and
construct a relative initial connection $\nabla_0$ on each $M_i$.
Note that $\nabla_0$ is defined Zariski locally.
In this paper, we will use only the Zariski topology.
In Section \ref{SS alge split nabla 0}, we define vector fields on each $M_i$ associated 
to the relative initial connection $\nabla_0$.
We use these vector fields instead of vector fields associated to time variables.
In Section \ref{SS Hamil}, 
we give the 2-form $\omega$ on $M_{\cC/T}^{\balpha}(\tilde{\bt},r,d)_{\bnu}$
and define Hamiltonian functions on each $M_i$.
The Hamiltonian functions depend on the choice of the relative initial connection $\nabla_0$.
Finally, we obtain a Hamiltonian description of the vector field on each $M_i$
induced by the isomonodromic deformation.

\section{Preliminaries}\label{Pre}

\subsection{Moduli space of stable parabolic connections}
Let $C$ be a smooth projective curve of genus $g$.
We put 
\begin{equation*}
T_n := \{ (t_1, \ldots,t_n  ) \in  C \times \cdots \times C \mid t_i \neq t_j \text{ for $i\neq j$} \} 
\end{equation*}
for a positive integer $n$. 
For integers $d,r$ with $r>0$, we put 
\begin{equation*}
N^{(n)}_r(d) := \left\{  (\nu^{(i)}_j)^{1\le i \le n}_{0\le j \le r-1}  \in \mathbb{C}^{nr} 
\middle|  d+ \sum_{i,j} \nu^{(i)}_j =0  \right\}.
\end{equation*}
Take members $\bt= (t_1,\ldots,t_n) \in T_n$ and $\bnu=( \nu^{(i)}_j)_{1\le i \le n,0\le j \le r-1} 
\in N^{(n)}_r(d)$.

\begin{Def}
We say $(E,\nabla, \{ l_*^{(i)} \}_{1\le i\le n})$ is a \textit{$(\bt , \bnu)$-parabolic connection 
of rank $r$ and of degree $d$ over $C$} if 
\begin{itemize}
\item[(1)] $E$ is a rank $r$ algebraic vector bundle on $C$,
\item[(2)] $\nabla \colon E \ra E \otimes \Omega^1_C(t_1 + \cdots + t_n)$ is a connection, 
that is, $\nabla$ is a homomorphism of sheaves 
satisfying $\nabla(fa)=a \otimes df + f \nabla(a)$ for $f \in \cO_C$ and $a \in E$, and 
\item[(3)] for each $t_i$, $l_*^{(i)}$ is a filtration $E|_{t_i} = l_0^{(i)} \supset l_1^{(i)} 
\supset \cdots \supset l_r^{(i)}=0$ such that $\dim (l_j^{(i)}/l_{j+1}^{(i)})=1$
and $(\res_{t_i}(\nabla)-\nu^{(i)}_j \mathrm{id}_{E|_{t_i}}) (l^{(i)}_j) \subset l_{j+1}^{(i)}$ 
for $j=0,\ldots,r-1$.
\end{itemize}
\end{Def}

The filtration $l_*^{(i)}$ ($1\le i \le n $) is said
to be a \textit{parabolic structure} on the vector bundle $E$.

\begin{Rem}
We have
\begin{equation*}
\deg E = \deg (\det(E)) = -\sum^n_{i=1} \res_{t_i} (\nabla_{\det (E)}) 
= - \sum^n_{i=1} \sum^{r-1}_{j=0} \nu_j^{(i)} =d.
\end{equation*}
\end{Rem}

Take rational numbers
\begin{equation*}
0 < \alpha^{(i)}_1 < \alpha^{(i)}_2 < \cdots < \alpha^{(i)}_r <1
\end{equation*}
for $i=1,\ldots , n$ satisfying $\alpha^{(i)}_j \neq \alpha^{(i')}_{j'}$ for $(i,j) \neq (i', j')$.
We choose a sufficiently generic $\balpha=(\alpha_j^{(i)})$.

\begin{Def}
A parabolic connection $(E,\nabla, \{ l^{(i)}_* \}_{1\le i \le n})$ is \textit{$\balpha$-stable} 
(resp. \textit{$\balpha$-semistable}) if for any proper nonzero subbundle $F\subset E$
satisfying $\nabla(F) \subset F \otimes \Omega^1_C(t_1 + \cdots +t_n)$, the inequality
\begin{equation*}
\frac{\deg F+\sum^n_{i=1}\sum^r_{j=1} \alpha^{(i)}_j \dim ( (F|_{t_i} \cap l_{j-1}^{(i)})/(F|_{t_i}
\cap l_j^{(i)}) ) }{\rank F} 
\underset{(\text{resp. $\le$})}{<} \frac{\deg E+\sum^n_{i=1}\sum^r_{j=1}\alpha^{(i)}_j
 \dim ( l_{j-1}^{(i)}/ l_j^{(i)} )}{ \rank E}
\end{equation*}
holds.
\end{Def}

Let $T$ be a smooth algebraic scheme which is an \'etale covering 
of the moduli stack $\cM_{g,n}$ of $n$-pointed smooth projective curves of genus $g$ 
over $\mathbb{C}$ 
and take a universal family $(\cC,\tilde{t}_1, \ldots,\tilde{t}_n)$ over $T$.

\begin{Def}
We denote the pull-back of $\cC$ and $\tilde{\bt}$ by the morphism $T \times N^{(n)}_r(d) \ra T$ 
by the same characters $\cC$
and $\tilde{\bt}=\{ \tilde{t}_1 ,\ldots,\tilde{t}_n \}$.
Then $D(\tilde{\bt}):=\tilde{t}_1 +\cdots+\tilde{t}_n$ becomes an effective Cartier divisor 
on $\cC$ flat over $T \times N^{(n)}_r(d)$.
We also denote by $\tilde{\bnu}$ the pull-back of the universal family on $N^{(n)}_r(d)$ 
by the morphism $T \times N^{(n)}_r(d) \ra N^{(n)}_r(d)$.
We define a functor $\cM_{\cC/T}^{\balpha}(\tilde{\bt},r,d)$ from the category 
of locally noetherian schemes over $T \times N^{(n)}_r(d)$ to the category of sets by 
\begin{equation*}
\cM_{\cC/T}^{\balpha}(\tilde{\bt},r,d)(S):=\left\{ (E,\nabla, \{l^{(i)}_j \}) \right\}/\sim
\end{equation*}
for a locally noetherian scheme $S$ over $T \times N^{(n)}_r(d)$, where 
\begin{itemize}
\item[(1)] $E$ is a rank $r$ algebraic vector bundle on $\cC_S$,
\item[(2)] $\nabla \colon E \ra E \otimes \Omega^1_{\cC_S/S}(D(\tilde{\bt})_S)$ 
is a relative connection, 
\item[(3)] for each $(\tilde{t}_i)_S$, $l_*^{(i)}$ is a filtration by subbundles
 $E|_{(\tilde{t}_i)_S} = l_0^{(i)} \supset l_1^{(i)} \supset \cdots \supset l_r^{(i)}=0$ 
such that $(\res_{(\tilde{t}_i)_S}(\nabla)-(\tilde{\nu}^{(i)}_j)_S \mathrm{id}_{E|_{t_i}}) (l^{(i)}_j) 
\subset l_{j+1}^{(i)}$ for $j=0,\ldots,r-1$, and
\item[(4)] for any point $s \in S$, $\dim (l_j^{(i)}/l_{j+1}^{(i)})\otimes k(s)=1$ for any $i,j$ 
and $(E,\nabla,\{ l^{(i)}_j \}) \otimes k(s)$ is $\balpha$-stable.
\end{itemize} 
Here $(E, \nabla,\{ l^{(i)}_j \})\sim (E', \nabla',\{ l'^{(i)}_j \})$ if there exist a line bundle $\cL$ on $S$ 
and an isomorphism $\sigma \colon E \xrightarrow{\sim} E'\otimes \cL$
such that $\sigma|_{(\tilde{t}_i)_S}(l^{(i)}_j) = l'^{(i)}_j \otimes \cL$ for any $i,j$ and the diagram 
\begin{equation*}
\xymatrix{
E \ar[r]^-{\nabla} \ar[d]^-{\sigma} & E \otimes \Omega^1_{\cC/T}(D(\tilde{\bt})) 
\ar[d]^-{\sigma\otimes \mathrm{id}} \\
E' \otimes \cL \ar[r]^-{\nabla'\otimes \cL} 
& E' \otimes  \Omega^1_{\cC/T}(D(\tilde{\bt}))\otimes \cL
}
\end{equation*} 
commutes.
\end{Def}

\begin{Thm}[{\cite[Theorem 2.1]{Inaba}}]
\textit{
For the moduli functor $\cM_{\cC/T}^{\balpha}(\tilde{\bt},r,d)$,
there exists a fine moduli scheme }
\begin{equation*}
\varpi\colon M_{\cC/T}^{\balpha}(\tilde{\bt},r,d) \lra T \times N^{(n)}_r(d) 
\end{equation*}
\textit{of $\balpha$-stable parabolic connections of rank $r$ and degree $d$, 
which is smooth and quasi-projective.
The fiber $M_{C_x}^{\balpha}(\tilde{\bt}_x, \bnu)$ over $(x,\bnu) \in T \times N^{(n)}_r(d)$ 
is the moduli space of $\balpha$-stable
$(\tilde{\bt}_x, \bnu)$-parabolic connections whose dimension is }
\begin{equation*}
2r^2(g-1) +nr(r-1)+2
\end{equation*}
\textit{
if it is non-empty.
}
\end{Thm}

\subsection{Description of the relative tangent sheaf 
of the phase space over the space of time variables}\label{Inf defor fixed curve}

We recall the description of the relative tangent sheaf 
$\Theta_{M_{\cC/T}^{\balpha}(\tilde{\bt} ,r,d)/T\times N^{(n)}_r(d)}$ 
in terms of the hypercohomology of a certain complex
(\cite[the proof of Theorem 2.1]{Inaba}).
Let $(\tilde{E},\tilde{\nabla} ,\{ \tilde{l}^{(i)}_j\})$ be a universal family 
on $\cC \times_TM_{\cC/T}^{\balpha}(\tilde{\bt} ,r,d)$.
First, we define a complex $\cF^{\bullet}$ by 
\begin{equation*}
\begin{aligned}
&\cF^0 := \left\{  s \in \cE nd (\tilde{E})  
\middle| s |_{\tilde{t}_i\times M_{\cC/T}^{\balpha}(\tilde{\bt} ,r,d)} (\tilde{l}^{(i)}_j) 
\subset \tilde{l}^{(i)}_j  \text{ for any $i,j$} \right\} \\
&\cF^1 :=  \left\{  s \in  \cE nd (\tilde{E})\otimes \Omega^1_{\cC/T}(D(\tilde{\bt}))  
\middle|  \res_{\tilde{t}_i\times M_{\cC/T}^{\balpha}(\tilde{\bt} ,r,d)} (s) (\tilde{l}^{(i)}_j) 
\subset \tilde{l}^{(i)}_{j+1} \text{ for any $i,j$}  \right\} \\
&\nabla_{\cF^{\bullet}} \colon \cF^0 \lra \cF^1; 
\quad \nabla_{\cF^{\bullet}} (s) = \tilde{\nabla}\circ s - s \circ \tilde{\nabla}.
\end{aligned}
\end{equation*}
Second, we take an affine open set $M\subset M_{\cC/T}^{\balpha}(\tilde{\bt} ,r,d)$ 
and an affine open covering $\cC_M = \bigcup_{\alpha} U_{\alpha}$ such that 
$\tilde{E}|_{U_{\alpha}} \cong \cO^{\oplus r}_{U_{\alpha}}$ for any $\alpha$,
$\sharp\{ i \mid \tilde{t}_i |_{\cC_M} \cap U_{\alpha} \neq \emptyset \} \le 1$ for any $\alpha$ and 
$\sharp\{ \alpha \mid \tilde{t}_i |_{\cC_M} \cap U_{\alpha} \neq \emptyset \} \le 1$ for any $i$.
Take a relative tangent vector field 
$v \in \Theta_{M_{\cC/T}^{\balpha}(\tilde{\bt} ,r,d)/T\times N^{(n)}_r(d)} (M)$.
The field $v$ corresponds to a member 
$(E_{\epsilon},\nabla_{\epsilon}, \{ (l_{\epsilon})^{(i)}_j \}) 
\in M_{\cC/T}^{\balpha}(\tilde{\bt} ,r,d)(\Spec \cO_M[\epsilon])$
such that $(E_{\epsilon},\nabla_{\epsilon}, 
\{ (l_{\epsilon})^{(i)}_j \}) \otimes \cO_{M}[\epsilon]/(\epsilon) 
\cong (\tilde{E},\tilde{\nabla},\{ \tilde{l}^{(i)}_j \})|_{\cC_M}$,
where $\cO_M[\epsilon]= \cO_M[t]/(t^2)$.
There is an isomorphism
\begin{equation*}\label{equation isom verphi 1}
\varphi_{\alpha} \colon  E_{\epsilon}|_{U_{\alpha}\times \Spec \cO_M[\epsilon] }
\xrightarrow{\sim} \cO^{\oplus r}_{U_{\alpha}\times \Spec \cO_M[\epsilon]} 
\xrightarrow{\sim} \tilde{E}|_{U_{\alpha}} \otimes \cO_M[\epsilon]
\end{equation*}
such that $\varphi_{\alpha}\otimes \cO_M[\epsilon]/(\epsilon ) 
\colon E_{\epsilon}\otimes \cO_M [\epsilon]/(\epsilon)|_{U_{\alpha}} 
\xrightarrow{\sim} \tilde{E}|_{U_{\alpha}}\otimes \cO_M [\epsilon]/(\epsilon)=\tilde{E}|_{U_{\alpha}}$ 
is the given isomorphism and 
that $\varphi_{\alpha}|_{t_i\otimes \cO_M[\epsilon]} ((l_{\epsilon})_j^{(i)}) 
= \tilde{l}_j^{(i)}|_{U_{\alpha} \times \Spec \cO_M[\epsilon]}$ 
if $\tilde{t}_i |_{\cC_M}\cap U_{\alpha} \neq \emptyset$.
We put 
\begin{equation*}
\begin{aligned}
u_{\alpha\beta} :=&\ \varphi_{\alpha} \circ \varphi_{\beta}^{-1} 
- \mathrm{id}_{\tilde{E}|_{U_{\alpha\beta}\times \Spec \cO_M[\epsilon]}},\\
v_{\alpha} :=&\ (\varphi_{\alpha}\otimes \mathrm{id}) \circ 
\nabla_{\epsilon} |_{U_{\alpha}\times \Spec \cO_M[\epsilon]} \circ \varphi^{-1}_{\alpha} 
- \tilde{\nabla}|_{U_{\alpha}\times \Spec \cO_M[\epsilon]}.
\end{aligned}
\end{equation*}
Then $\{ u_{\alpha\beta} \} \in C^1 ((\epsilon) \otimes \cF^0_M)$, 
$\{ v_{\alpha} \} \in C^0((\epsilon) \otimes \cF^1_M)$ and 
\begin{equation*}
d\{ u_{\alpha \beta} \} = \{ u_{\beta \gamma}-u_{\alpha \gamma} +u_{\alpha\beta} \} = 0;
\quad \nabla_{\cF^{\bullet}} \{ u_{\alpha\beta} \} = \{ v_{\beta}-v_{\alpha} \} =d\{ v_{\alpha} \}.
\end{equation*}
So $[(\{ u_{\alpha\beta} \},\{ v_{\alpha} \} )]$ determines an element $\sigma_M(v)$ 
of $\bH^1(\cF_M^{\bullet})$.
We can check that $v \mapsto \sigma_M (v)$ determines an isomorphism
\begin{equation*}
\sigma_M \colon \Theta_{M_{\cC/T}^{\balpha}(\tilde{\bt} ,r,d)/T\times N^{(n)}_r(d)} (M) 
\xrightarrow{\, \sim \, } \bH^1(\cF_M^{\bullet});
\quad v \longmapsto \sigma_M(v).
\end{equation*}
The isomorphism $\sigma_M$ induces a canonical isomorphism 
\begin{equation*}
\sigma \colon \Theta_{M_{\cC/T}^{\balpha}(\tilde{\bt} ,r,d)/T\times N^{(n)}_r(d)} 
\xrightarrow{\, \sim \, } \bold{R}^1(\pi_{M_{\cC/T}^{\balpha}(\tilde{\bt} ,r,d)})_*(\cF^{\bullet}),
\end{equation*}
where $\pi_{M_{\cC/T}^{\balpha}(\tilde{\bt} ,r,d)} \colon 
\cC_{M_{\cC/T}^{\balpha}(\tilde{\bt} ,r,d)} \rightarrow M_{\cC/T}^{\balpha}(\tilde{\bt} ,r,d)$ 
is the natural morphism.

\subsection{Symplectic structure}\label{S2 symp form}

For each affine open subset $U \subset M_{\cC/T}^{\balpha}(\tilde{\bt} ,r,d)$, we define a pairing
\begin{equation}\label{2020.11.7.15.48}
\begin{aligned}
\bH^1(\cC \times_T U, \cF_U^{\bullet}) \otimes \bH^1(\cC \times_T U, \cF_U^{\bullet}) 
&\lra \bH^2(\cC \times_T U, \Omega_{\cC\times_T U/U}^{\bullet}) \cong H^0(\cO_U)\\
[(\{ u_{\alpha\beta} \}, \{ v_{\alpha} \})]\otimes[(\{ u_{\alpha\beta}' \}, \{ v_{\alpha}' \} )] &\longmapsto 
[ (\{  \Tr( u_{\alpha\beta} \circ u_{\beta\gamma}') \}, -\{ \Tr (u_{\alpha\beta} \circ v_{\beta}') 
- \Tr (v_{\alpha} \circ u'_{\alpha\beta}) \} )],
\end{aligned}
\end{equation}
considered in \v{C}ech cohomology with respect to an affine open covering $\{ U_{\alpha} \}$ 
of $\cC \times_TU$, $\{ u_{\alpha\beta} \} \in C^1(\cF^0)$, $\{ v_{\alpha} \} \in C^0(\cF^1)$
and so on.
This pairing determines a pairing 
\begin{equation*}
\omega \colon \bold{R}^1 (\pi_{M_{\cC/T}^{\balpha}(\tilde{\bt} ,r,d)})_*(\cF^{\bullet}) 
\otimes \bold{R}^1 (\pi_{M_{\cC/T}^{\balpha}(\tilde{\bt} ,r,d)})_*(\cF^{\bullet}) \lra 
\cO_{M_{\cC/T}^{\balpha}(\tilde{\bt} ,r,d)}.
\end{equation*}
This pairing is a nondegenerate relative 2-form.
This fact follows from the Serre duality and the five lemma:
\begin{equation*}
\xymatrix{
H^0(\cF^0_x) \ar[r] \ar[d]^-{\sim} & H^0(\cF^1_x) \ar[r] \ar[d]^-{\sim} 
& \bH^1(\cF_x^{\bullet}) \ar[r] \ar[d]^-{\sim} & H^1(\cF^0_x) \ar[r] \ar[d]^-{\sim} 
& H^1(\cF^1_x) \ar[d]^-{\sim} \\
H^1(\cF^1_x)^{\vee} \ar[r] & H^1(\cF^0_x)^{\vee} \ar[r] 
& \bH^1(\cF_x^{\bullet})^{\vee} \ar[r] & H^0(\cF^1_x)^{\vee} \ar[r] 
& H^0(\cF^0_x)^{\vee}
}
\end{equation*}
for any point $x \in M_{\cC/T}^{\balpha}(\tilde{\bt} ,r,d)$.
Moreover we have $d\omega=0$ (see \cite[Proposition 7.3]{Inaba}).

\section{Isomonodromic deformation}\label{S IM}

In this section, we consider a description of the vector field determined 
by the isomonodromic deformation in terms of a \v{C}ech cohomology.
In Section \ref{SS AA}, we recall the Atiyah algebra.
By the Atiyah algebra, we obtain descriptions of first-order deformations of pairs of 
pointed curves 
and vector bundles.
In Section \ref{SS d t s of m by c c}, we show that the relative tangent sheaf 
$\Theta_{M_{\cC/T}^{\balpha}(\tilde{\bt} ,r,d)/ N^{(n)}_r(d)}$ 
is isomorphic to the hypercohomology of a certain complex using the Atiyah algebra.
In Section \ref{SS IMD by c c}, 
we consider the integrable deformations of parabolic connections 
when the $n$-pointed curves vary.
The integrable deformations of parabolic connections mean the isomonodromic deformations 
of the corresponding relative parabolic connections.
We describe the vector field determined by the isomonodromic deformation 
in terms of a \v{C}ech cohomology.

\subsection{Atiyah algebra}\label{SS AA}

We recall the \textit{Atiyah algebra}.
(For details, for example see \cite{BS}).
Let $C$ be a smooth projective curve and $\Theta_C$ be the tangent sheaf.
Let $E$ be a vector bundle of rank $r$ on $C$.
We define a sheaf of differential operators on $E$ as follows.
(For detail, see \cite[Section 1]{BB}).
We take an affine open subset $U \subset C$.
For the $\mathcal{O}(U)$-bimodule $\mathrm{End}_{\mathbb{C}}(E(U),E(U))$,
we define $\mathrm{End}_{\mathbb{C}}(E(U),E(U))_i^{\vee}$ (for $i \ge -1$)
by induction
as follows:
$$
\begin{aligned}
\mathrm{End}_{\mathbb{C}}(E(U),E(U))_{-1}^{\vee} &:=0 \\
\mathrm{End}_{\mathbb{C}}(E(U),E(U))_{i}^{\vee} &:=
\left\{ s \in \mathrm{End}_{\mathbb{C}}(E(U),E(U))
\ \middle|\ 
\begin{array}{l}
r \cdot s - s \cdot r \in \mathrm{End}_{\mathbb{C}}(E(U),E(U))_{i-1}^{\vee}\\
\text{for any $r \in \mathcal{O}(U)$} 
\end{array}
\right\}.
\end{aligned}
$$
\begin{Def}[{see \cite[1.1.6]{BB}}]
We say a $\mathbb{C}$-linear morphism $s \colon E \rightarrow E$
is a \textit{differential operator of degree $i$} if 
for any affine subset $U \subset C$ 
the corresponding morphism $s_U \colon E(U) \rightarrow E(U)$
lies in $\mathrm{End}_{\mathbb{C}}(E(U),E(U))_{i}^{\vee}$
and is not contained in 
$\mathrm{End}_{\mathbb{C}}(E(U),E(U))_{i-1}^{\vee}$.
\end{Def}

Put $\cD_E=\mathcal{D}\it{iff}(E,E)=\bigcup_i \cD_i$, 
where $\cD_i$ is the sheaf of differential operators of degree $\le i$ on $C$.
We have $\cD_i/\cD_{i-1}= \mathcal{E}nd E \otimes S^i(\Theta_C)$, 
where $S^i(\Theta_C)$ is the $i$-th symmetric product of $\Theta_C$.
We denote by $\mathrm{symb}_1$ 
the composition 
\[
\cD_1 
\longrightarrow \cD_1/\cD_0
\xrightarrow{\ =\ }
\mathcal{E}nd E \otimes \Theta_C.
\]
For the differential operator $v$ of degree $1$,
the image $\mathrm{symb}_1(v)$ is called the {\it symbol} of $v$.
\begin{Def}
We define the \textit{Atiyah algebra} of $E$ as
\begin{equation*}
\cA_E = \{ \partial \in \cD_1 \mid  \mathrm{symb}_1(\partial) 
\in  \mathrm{id}_E\otimes \Theta_C \subset \cE nd(E) \otimes \Theta_C \}.
\end{equation*}
\end{Def}
We have inclusions $\cD_0= \mathcal{E}nd E \subset \cA_E \subset \cD_1$ 
and a short exact sequence
\begin{equation*}\label{ES Atiyah}
0\lra \cE nd (E) \lra \cA_E \xrightarrow{\mathrm{symb}_1} \Theta_C \lra 0.
\end{equation*}
Fix a positive integer $n$.
Let $D=t_1+\cdots+t_n$ be an effective divisor of $C$, 
where $t_1 ,\ldots ,t_n$ are distinct points of $C$.
We put $\cA_E(D):=\mathrm{symb}_1^{-1} (\Theta_C(-D))$.
Then we have the following exact sequence
\begin{equation*}
0\lra \cE nd (E) \lra \cA_E(D) \xrightarrow{\mathrm{symb}_1} \Theta_C(-D) \lra 0.
\end{equation*}
By this exact sequence, we have the following exact sequence
\begin{equation*}
0\lra H^1(C ,\cE nd (E) ) \lra H^1(C, \cA_E (D) ) 
\xrightarrow{\mathrm{symb}_1} H^1(C, \Theta_C (-D) ) \lra 0
\end{equation*}
when $2g-2+n>0$.
Then $H^1(C, \cA_E (D) )$ means the set of infinitesimal deformations of the pair 
$((C,t_1,\ldots,t_n),E)$ of the $n$-pointed curve $(C,t_1,\ldots,t_n)$
and the vector bundle $E$ on $C$.

For a connection $\nabla \colon E \rightarrow E \otimes \Omega^1_{C}(D)$,
we define a splitting 
\begin{equation}\label{spritting nabla}
\iota (\nabla) \colon \Theta_C(-D) \lra \cA_E(D)
\end{equation}
by
\begin{equation*}
(\iota (\nabla) (v)) (u) = \langle v, \nabla (u) \rangle \in E
\end{equation*}
for $v \in \Theta_C(-D)$ and $u \in E$.
The splitting $\iota(\nabla)$ is locally described as follows.
Let $U$ be an affine open subset of $C$ 
where we have a trivialization $E|_{U} \cong \cO_U^{\oplus r}$. 
We denote by $Af^{-1} df$ a connection matrix of $\nabla$ on $U$, 
where $f$ is a local defining equation of $t_i$ and $A \in M_r (\cO_U)$. 
For an element $g\frac{\partial}{\partial f} \in \Theta_C(-D)(U)$,
we have an equality $\iota (\nabla)(g\frac{\partial}{\partial f}) = 
g \left( \frac{\partial}{\partial f}  +   A f^{-1} \right) \in \cA_E(D)(U)$.
Conversely, 
a splitting of $\mathrm{symb}_1 \colon \cA_E(D) \rightarrow \Theta_C(-D)$
gives a connection as follows.
A splitting of $\mathrm{symb}_1$
gives a splitting of $\cA_E(D) \otimes \Omega_C(D) \rightarrow  \mathcal{O}_C$.
Let $(f\frac{\partial}{\partial f} +A) \otimes \frac{df}{f} $ be 
a local description of the image of $1 \in \mathcal{O}_C$ by this splitting.
Remark that the image of $(f\frac{\partial}{\partial f} +A) \otimes \frac{df}{f} $ 
under the morphism $\cA_E(D) \otimes \Omega_C(D) \rightarrow  \mathcal{O}_C$ is $1$.
The image of $1$ gives a morphism
\begin{equation*}
\begin{aligned}
E &\longrightarrow E \otimes \Omega^1_C(D)\\
a &\longmapsto 
\left(\frac{\partial}{\partial f}(a) +\frac{A}{f}(a) \right) \otimes df.
\end{aligned}
\end{equation*}
This is a connection $E \rightarrow E \otimes \Omega_C(D)$.

\subsection{Description of the total tangent sheaf 
of the phase space}\label{SS d t s of m by c c}

We discuss a description of the relative tangent sheaf 
$\Theta_{M_{\cC/T}^{\balpha}(\tilde{\bt} ,r,d)/ N^{(n)}_r(d)}$ 
in terms of the hypercohomology of a certain complex.
Let $(\tilde{E},\tilde{\nabla} ,\{ \tilde{l}^{(i)}_j\})$ be a universal family 
on $\cC \times_TM_{\cC/T}^{\balpha}(\tilde{\bt},r,d)$.
We put
\begin{equation*}
\begin{aligned}
&\cG^0 := \left\{  s \in \cA_{\tilde{E}}(D(\tilde{\bt})) \ 
\middle|\ s |_{\tilde{t}_i\times_T M_{\cC/T}^{\balpha}(\tilde{\bt} ,r,d)} (\tilde{l}^{(i)}_j) 
\subset \tilde{l}^{(i)}_j  \text{ for any $i,j$} \right\} \\
&\cG^1 :=  \left\{  s \in  \cE nd (\tilde{E})\otimes \Omega^1_{\cC/T}(D(\tilde{\bt}))\  
\middle| \  \res_{\tilde{t}_i\times_T M_{\cC/T}^{\balpha}(\tilde{\bt} ,r,d)} (s) (\tilde{l}^{(i)}_j) 
\subset \tilde{l}^{(i)}_{j+1} \text{ for any $i,j$}  \right\},
\end{aligned}
\end{equation*}
where $\cA_{\tilde{E}}(D(\tilde{\bt}))$ is the relative Atiyah algebra which is the extension
\begin{equation*}
0\lra \cE nd(\tilde{E}) \lra \cA_{\tilde{E}}(D(\tilde{\bt})) 
\lra \Theta_{\cC \times_T M_{\cC/T}^{\balpha}(\tilde{\bt},r,d)/ M_{\cC/T}^{\balpha}(\tilde{\bt},r,d)} 
(-D(\tilde{\bt})) \lra 0.
\end{equation*}
Here note that   
$s |_{\tilde{t}_i\times_T M_{\cC/T}^{\balpha}(\tilde{\bt} ,r,d)} (\tilde{l}^{(i)}_j)$
is well-defined for $s \in \cA_{\tilde{E}}(D(\tilde{\bt}))$.
Indeed, we take an affine open subset $U$ of $\cC \times_T M_{\cC/T}^{\balpha}(\tilde{\bt},r,d)$
such that $(\tilde{t}_i\times_T M_{\cC/T}^{\balpha}(\tilde{\bt} ,r,d))\cap U\neq \emptyset$,
and we take a trivialization $\bar{\phi} \colon  \tilde{E}|_{U} \rightarrow \mathcal{O}_U^{\oplus r}$.
If we take a local description $s=s_0 (\partial/\partial f) +s_1$ 
(where
$s_0 (\partial/\partial f) \in  \Theta_{U/ M_{\cC/T}^{\balpha}(\tilde{\bt},r,d)} 
(-D(\tilde{\bt})|_{U})$
and $s_1 \in \mathcal{E}nd (\mathcal{O}_U^{\oplus r})$),
then we define 
$$s |_{(\tilde{t}_i\times_T M_{\cC/T}^{\balpha}(\tilde{\bt} ,r,d))\cap U} (\tilde{l}^{(i)}_j):=
\bar{\phi}^{-1} ( (s_1 \circ \bar{\phi})
 |_{(\tilde{t}_i\times_T M_{\cC/T}^{\balpha}(\tilde{\bt} ,r,d))\cap U} (\tilde{l}^{(i)}_j)).$$
This definition is independent of the choice of a trivialization of $\tilde{E}|_{U}$,
since 
$$
g^{-1} \circ \left( s_0 \frac{\partial}{\partial f}  +s_1 \right)\circ g
=s_0 \frac{\partial}{\partial f}+
s_0 g^{-1} \frac{\partial g}{\partial f}  +g^{-1} s_1 g
$$
for $g \in \mathrm{Aut}(\mathcal{O}_U^{\oplus r})$
and $s_0 g^{-1} \frac{\partial g}{\partial f}$ vanishes 
at $(\tilde{t}_i\times_T M_{\cC/T}^{\balpha}(\tilde{\bt} ,r,d))\cap U$.
We define a complex $\cG^{\bullet}$ by
\begin{equation*}
\begin{aligned}
\nabla_{\cG^{\bullet}} \colon \cG^0 &\lra \cG^1; \\
 s &\longmapsto  
 \tilde{\nabla}\circ \left(s- \iota(\tilde{\nabla}) \circ \mathrm{symb}_1 (s) \right) 
- \left( \left(s- \iota(\tilde{\nabla}) \circ \mathrm{symb}_1 (s) \right) \otimes 1\right)
 \circ \tilde{\nabla}.
 \end{aligned}
\end{equation*}
In fact, we use the description \eqref{2020.5.22.23.09} below
for calculation of the complex $\cG^{\bullet}$.

Let $U$ be an affine open subset of $\cC \times_T M_{\cC/T}^{\balpha}(\tilde{\bt},r,d)$.
we take a trivialization of $\tilde{E}|_{U}$.
Let $d$ be the relative exterior derivative on $U \rightarrow M_{\cC/T}^{\balpha}(\tilde{\bt},r,d)$. 
We take $s \in \cG^0(U)$, which has a local description $s_0 (\partial/\partial f) +s_1$.
Here $s_1 \in \mathcal{E}nd (\mathcal{O}_{U}^{\oplus r})$ and 
$\mathrm{symb}_1(s) = s_0 (\partial/\partial f)$.

\begin{Def}
For the section $s=s_0 (\partial/\partial f) +s_1\in \cG^0(U)$, 
we define a homomorphism of sheaves 
$s\colon \mathcal{O}_{U}^{\oplus r}\otimes \Omega^1_{\cC/T}(D(\tilde{\bt}))
\rightarrow \mathcal{O}_{U}^{\oplus r} \otimes \Omega^1_{\cC/T}(D(\tilde{\bt}))$
by 
\begin{equation}\label{2019.5.22.12.11}
\left(s_0 \frac{\partial}{\partial f} +s_1 \right)( a )
 := d \left\langle 1\otimes s_0 \frac{\partial}{\partial f} ,a \right\rangle 
 +( s_1 \otimes 1) a 
\end{equation}
for any $a \in  \mathcal{O}_U^{\oplus r}\otimes \Omega^1_{\cC/T}(D(\tilde{\bt}))$.
Here, $d \langle 1\otimes s_0 \partial/\partial f , \rangle $ means 
the composition
\begin{equation*}
\mathcal{O}_U^{\oplus r}\otimes \Omega^1_{\cC/T}(D(\tilde{\bt}))
\xrightarrow{\langle 1\otimes s_0 \partial/\partial f , \rangle} \mathcal{O}_U^{\oplus r}
\xrightarrow{\ d\ } \mathcal{O}_U^{\oplus r}\otimes \Omega^1_{\cC/T}
\longrightarrow \mathcal{O}_U^{\oplus r}\otimes \Omega^1_{\cC/T}(D(\tilde{\bt})).
\end{equation*}
\end{Def}

Since 
$\cC \times_T M_{\cC/T}^{\balpha}(\tilde{\bt},r,d)\ra M_{\cC/T}^{\balpha}(\tilde{\bt},r,d)$ 
is a family of smooth projective curves,
the definition \eqref{2019.5.22.12.11} is independent of the choice of the 
parameter $f$. 
Indeed, 
let $a_0 df$ be a element of $\mathcal{O}_U^{\oplus r}\otimes \Omega^1_{\cC/T}(D(\tilde{\bt}))$.
Then $\left(s_0 \frac{\partial}{\partial f} +s_1 \right)( a_0df )=
d \left(  s_0  a_0  \right)
 +( s_1 \otimes 1) a_0 df $.
We have the following equalities:
\begin{equation*}
\begin{aligned}
\left(s_0 \frac{\partial g}{\partial f}\frac{\partial}{\partial g} +s_1 \right)
( a_0 df  )
&= d \left\langle 1\otimes s_0 \frac{\partial g}{\partial f}\frac{\partial}{\partial g}  ,
a_0 \frac{\partial f}{\partial g} dg \right\rangle 
 +( s_1 \otimes 1) a_0 df \\
 &=d \left(  s_0  a_0  \right)
 +( s_1 \otimes 1) a_0 df .
 \end{aligned}
\end{equation*}
These equality mean that 
the definition \eqref{2019.5.22.12.11} is independent of the choice of the 
parameter $f$.
The definition \eqref{2019.5.22.12.11} comes from 
the local description of the Lie derivative.

\begin{Lem}\label{2020.5.24.14.29}
\textit{
We take a trivialization of $\tilde{E}$ on $U$.
Let $Af^{-1} df$ be the connection matrix of $\tilde{\nabla}$ on $U$ with respect to 
the trivialization. Then we have the following equality:
\begin{equation*}
\begin{aligned}
&\left( d + A \frac{df}{f} \right) \left( \left(s_1 - s_0 \frac{A}{f}  \right)(a)  \right)-
\left(s_1 - s_0 \frac{A}{f}  \right)\left( \left( d + A \frac{df}{f} \right)(a)  \right)\\
&=\left( d + A \frac{df}{f} \right) \left( \left(s_0 \frac{\partial}{\partial f} +s_1 \right)(a)\right)
-\left(s_0 \frac{\partial}{\partial f} +s_1 \right)\left( \left( d + A \frac{df}{f} \right) (a) \right)
\end{aligned}
\end{equation*}
for any $a \in  \mathcal{O}_U^{\oplus r}\otimes \Omega^1_{\cC/T}(D(\tilde{\bt}))$.
That is, we have 
\begin{equation}\label{2020.5.22.23.09}
\nabla_{\cG^{\bullet}} (s)(a) = \tilde{\nabla} \circ s (a) - s \circ \tilde{\nabla} (a) 
\end{equation}
for any $a \in  \tilde{E}|_U\otimes \Omega^1_{\cC/T}(D(\tilde{\bt}))$.
}
\end{Lem}

\begin{proof}
From the right hand side, we compute as follows:
\begin{equation*}
\begin{aligned}
&\left( d + A \frac{df}{f} \right) \left( \left(s_0 \frac{\partial}{\partial f} +s_1 \right)(a)\right)
-\left(s_0 \frac{\partial}{\partial f} +s_1 \right)\left( \left( d + A \frac{df}{f} \right) (a) \right)\\
&=
\left( d + A \frac{df}{f} \right) \left( s_0 \frac{\partial a}{\partial f} \right)
+\left( d + A \frac{df}{f} \right) \left(s_1(a) \right)
- s_0 \frac{\partial}{\partial f}  \left( d(a) + A a \frac{df}{f} \right) 
-s_1  \left( d(a) + A a \frac{df}{f} \right)\\
&=
\left( d + A \frac{df}{f} \right) \left( s_0 \frac{\partial a}{\partial f} \right)
+\left( d + A \frac{df}{f} \right) \left(s_1(a) \right)\\
&\qquad - d\left\langle s_0 \frac{\partial}{\partial f}  , d(a) \right\rangle 
-  d \left\langle s_0 \frac{\partial}{\partial f} , Aa \frac{df}{f} \right\rangle  
-s_1  \left( d(a) + A a \frac{df}{f} \right)\\
&=
 s_0 \frac{A}{f}  \frac{\partial a}{\partial f} df 
+\left( d + A \frac{df}{f} \right) \left(s_1(a) \right) 
-   d \left(  s_0  \frac{A}{f} (a) \right) 
-s_1  \left( d(a) + A a \frac{df}{f} \right).
\end{aligned}
\end{equation*}
On the other hand, we have the following equalities:
\begin{equation*}
\begin{aligned}
&\left( d + A \frac{df}{f} \right) \left( \left(s_1 - s_0 \frac{A}{f}  \right)(a)  \right)-
\left(s_1 - s_0 \frac{A}{f}  \right)\left( \left( d + A \frac{df}{f} \right)(a)  \right)\\
&=    \left( d + A \frac{df}{f} \right) (s_1(a)) - 
\left( d + A \frac{df}{f} \right) \left( s_0  \frac{A}{f}  (a) \right)  -
\left(s_1 - s_0 \frac{A}{f}  \right)\left( \left( d(a) + A a \frac{df}{f} \right)  \right)\\
&=    \left( d + A \frac{df}{f} \right) (s_1(a)) - 
 d  \left( s_0  \frac{A}{f}  (a) \right)
-  A \frac{df}{f}  \left( s_0  \frac{A}{f}  (a) \right)  \\
&\qquad- s_1  \left( d(a) + A a \frac{df}{f} \right)  
+  s_0 \frac{A}{f}   \left( d(a) + A a \frac{df}{f} \right)  \\
&=
 s_0 \frac{A}{f}  \frac{\partial a}{\partial f} df 
+\left( d + A \frac{df}{f} \right) \left(s_1(a) \right) 
-   d \left(  s_0  \frac{A}{f} (a) \right) 
-s_1  \left( d(a) + A a \frac{df}{f} \right),
\end{aligned}
\end{equation*}
since $s_0$ is a scalar.
By the calculation above, we obtain the equalities in the assertion of this lemma.
\end{proof}

\begin{Prop}\label{Prop d t s of m by c c}
\textit{
There exists a canonical isomorphism 
\begin{equation*}
\varsigma \colon \Theta_{M_{\cC/T}^{\balpha}(\tilde{\bt} ,r,d)/N^{(n)}_r(d)} 
\xrightarrow{\, \sim \, } \bold{R}^1(\pi_{M_{\cC/T}^{\balpha}(\tilde{\bt} ,r,d)})_*(\cG^{\bullet}),
\end{equation*}
where $\pi_{M_{\cC/T}^{\balpha}(\tilde{\bt} ,r,d)} \colon 
\cC_{M_{\cC/T}^{\balpha}(\tilde{\bt} ,r,d)} \rightarrow M_{\cC/T}^{\balpha}(\tilde{\bt} ,r,d)$ 
is the natural morphism.
}
\end{Prop}

\begin{proof}
We take an affine open set $M\subset M_{\cC/T}^{\balpha}(\tilde{\bt} ,r,d)$.
We also denote by $(\tilde{E},\tilde{\nabla} ,\{ \tilde{l}^{(i)}_j\})$ the family 
on $\cC_M=\cC \times_{T \times N^{(n)}_r(d)} M$ induced by the universal family.
Let $D(\tilde{\bt})_{M}$ be the pull-back of $D(\tilde{\bt})$ by the morphism $\cC_{M} \ra \cC$.
We take an affine open covering $\cC_M = \bigcup_{\alpha} U_{\alpha}$ such that 
we have $\bar{\phi}_{\alpha} \colon  \tilde{E}|_{U_{\alpha}} 
\xrightarrow{\sim} \cO^{\oplus r}_{U_{\alpha}}$ for any $\alpha$,
$\sharp\{ i \mid \tilde{t}_i |_{\cC_M} \cap U_{\alpha} \neq \emptyset \} \le 1$ for any $\alpha$ and 
$\sharp\{ \alpha \mid \tilde{t}_i |_{\cC_M} \cap U_{\alpha} \neq \emptyset \} \le 1$ for any $i$.

Take a relative tangent vector field 
$v \in \Theta_{M_{\cC/T}^{\balpha}(\tilde{\bt} ,r,d)/N^{(n)}_r(d)} (M)$.
Put $M[\epsilon]=\Spec \cO_M[\epsilon]$, 
where $\cO_M[\epsilon]= \cO_M[t]/(t^2)$.
The field $v$ corresponds to a morphism 
$M[\epsilon] \rightarrow M$ over $N^{(n)}_r(d)$. 
Let  
$(E_{\epsilon},\nabla_{\epsilon}, \{ (l_{\epsilon})^{(i)}_j \}) $
be the flat family of parabolic connections on
$\cC \times_{T \times N^{(n)}_r(d)}  M[\epsilon]$ over $M[\epsilon]$
induced by this morphism $M[\epsilon] \rightarrow M$
and the flat family $(\tilde{E},\tilde{\nabla} ,\{ \tilde{l}^{(i)}_j\})$ 
on $\cC \times_{T \times N^{(n)}_r(d)}M$ over $M$.
Set 
\begin{itemize}
\item $\mathcal{C}_{\epsilon}:= \mathcal{C} \times_{T \times N^{(n)}_r(d)} M[\epsilon]$,
\item $(\tilde{t}_i)_{\epsilon} := \tilde{t}_i \times_{T \times N^{(n)}_r(d)} M[\epsilon]$ ($i=1,\ldots,n$), and
\item $D(\tilde{\bt})_{\epsilon}:= (\tilde{t}_1)_{\epsilon}+\cdots+(\tilde{t}_n)_{\epsilon}$.
\end{itemize}
The tuple $(\mathcal{C}_{\epsilon},(\tilde{t}_i)_{\epsilon},\ldots,(\tilde{t}_i)_{\epsilon})$
is the family of $n$-pointed curves over $M[\epsilon]$
induced by the composition $M[\epsilon] \rightarrow M\rightarrow 
T \times N^{(n)}_r(d)$ over $N^{(n)}_r(d)$
and the family $(\cC,\tilde{t}_1,\ldots,\tilde{t}_n)$ over $T \times N^{(n)}_r(d)$.
Remark that, 
since the morphism $T \rightarrow \mathcal{M}_{g,n}$ is \'etale,
giving a morphism $M[\epsilon] \rightarrow T \times  N^{(n)}_r(d)$ over $N^{(n)}_r(d)$ 
is equivalent to giving a flat family 
$(\mathcal{C}_{\epsilon},(\tilde{t}_1)_{\epsilon},\ldots,(\tilde{t}_n)_{\epsilon})$
of $n$-pointed curves satisfying 
$(\mathcal{C}_{\epsilon},(\tilde{t}_1)_{\epsilon},\ldots,(\tilde{t}_n)_{\epsilon})
\otimes \mathcal{O}_M[\epsilon]/(\epsilon) = 
(\cC,\tilde{t}_1,\ldots,\tilde{t}_n)$. 
Let $\cC_{\epsilon} = \bigcup_{\alpha} U^{\epsilon}_{\alpha}$ 
be the open covering corresponding to the affine open covering of $\cC_M$.
There is an $M[\epsilon]$-morphism $\sigma_{\alpha} \colon U^{\epsilon}_{\alpha} 
\ra U_{\alpha}\times \Spec \mathbb{C}[\epsilon]$ 
which is a lift of $\mathrm{id}_{U_{\alpha}}$ 
preserving the divisor $D(\tilde{\bt})_{\epsilon} \cap U^{\epsilon}_{\alpha}$ 
and $(D(\tilde{\bt})\cap U_{\alpha})\times \Spec \mathbb{C}[\epsilon]$.
If we put $d_{\alpha\beta} := (\sigma_{\alpha}^*)^{-1} \circ \sigma_{\beta}^* - 1 
\colon \cO_{U_{\alpha\beta}} \ra \epsilon \otimes \cO_{U_{\alpha\beta}}$,
then $d_{\alpha\beta}$ becomes a derivation and $[\{ d_{\alpha\beta} \}] 
\in H^1(\{ U_{\alpha} \}, \Theta_{\cC_M/M} (-D(\tilde{\bt})_M ))$ 
gives the Kodaira--Spencer class induced by $M[\epsilon] \xrightarrow{v} M \ra T \ra \cM_{g,n}$.
If we take a frame $\phi_{\alpha} \colon E_{\epsilon}|_{U_{\alpha}^{\epsilon}} 
\rightarrow \cO_{U_{\alpha}^{\epsilon}}^{\oplus r}$ such that
$\bar{\phi}_{\alpha}= \phi_{\alpha}\ (\text{mod } \epsilon)$,
then there is a composition of isomorphisms
\begin{equation}\label{isom varphi}
\varphi_{\alpha} \colon \tilde{E}|_{U_{\alpha}} \otimes_{\mathbb{C}} \mathbb{C}[\epsilon] 
\xrightarrow{\bar{\phi}_{\alpha}\otimes 1}  
\cO^{\oplus r}_{U_{\alpha}} \otimes_{\mathbb{C}} \mathbb{C}[\epsilon] 
= 
\cO^{\oplus r}_{U_{\alpha}} \otimes_{\cO_{U_{\alpha}}} (\cO_{U_{\alpha}} 
\otimes_{\mathbb{C}} \mathbb{C}[\epsilon] )
\xrightarrow{1 \otimes \sigma_{\alpha}^*} 
\cO^{\oplus r}_{U_{\alpha}} \otimes_{\cO_{U_{\alpha}}} \cO_{U^{\epsilon}_{\alpha}} 
\xrightarrow{\phi_{\alpha}^{-1}} 
E_{\epsilon}|_{U^{\epsilon}_{\alpha}}.
\end{equation} 
Since $\sharp\{ \alpha \mid \tilde{t}_i |_{\cC_M} \cap U_{\alpha} \neq \emptyset \} \le 1$,
if $\tilde{t}_i |_{\cC_M}\cap U_{\alpha} \neq \emptyset$,
then $\tilde{t}_i|_{\cC_M} \subset U_{\alpha}$.
By a change of the frame of $E_{\epsilon}|_{U_{\alpha}^{\epsilon}}$ 
for each $\alpha$ such that $\tilde{t}_i |_{\cC_M}\cap U_{\alpha} \neq \emptyset$,
we can assume that $\varphi_{\alpha} |_{\tilde{t}_i\times \Spec \mathbb{C}[\epsilon]} (\tilde{l}_j^{(i)}
\otimes_{\mathbb{C}} \mathbb{C}[\epsilon]) 
= (l_{\epsilon})_j^{(i)} $ if $\tilde{t}_i |_{\cC_M}\cap U_{\alpha} \neq \emptyset$.
We put 
\begin{equation*}
u_{\alpha\beta} :=  \varphi^{-1}_{\alpha} \circ \varphi_{\beta}  - \mathrm{id} \colon
\tilde{E}|_{U_{\alpha\beta}}  \lra \epsilon \otimes \tilde{E}|_{U_{\alpha\beta}} .
\end{equation*}
Put $\bar{g}_{\alpha\beta}:= (\bar{\phi}_{\alpha}\otimes 1)(\bar{\phi}_{\beta}\otimes 1)^{-1} $ 
and
$g_{\alpha\beta}:= \phi_{\alpha} \phi_{\beta}^{-1}$.
We define $g'_{\alpha\beta}$ and $b_{\alpha\beta}$ by
\begin{equation*}
\bar{g}_{\alpha\beta} (\mathrm{id} + \epsilon g'_{\alpha\beta}) = 
(1 \otimes \sigma_{\alpha}^*)^{-1} g_{\alpha\beta} 
(1 \otimes \sigma_{\alpha}^*)
\quad
\text{and}
\quad 
\bar{g}_{\alpha\beta} (\mathrm{id} + \epsilon b_{\alpha\beta}) = 
(1 \otimes \sigma_{\alpha}^*)^{-1} g_{\alpha\beta} 
(1 \otimes \sigma_{\beta}^*).
\end{equation*}
Here $g'_{\alpha\beta} 
\in \mathrm{End}( \mathcal{O}^{\oplus r}_{U_{\alpha\beta}}, 
\mathcal{O}^{\oplus r}_{U_{\alpha\beta}} )$ and 
$b_{\alpha\beta} \colon \cO^{\oplus r}_{U_{\alpha\beta}} \ra \cO^{\oplus r}_{U_{\alpha\beta}}$ 
is a differential operator of degree $\leq 1$ 
satisfying 
\begin{equation}\label{2020.5.22.21.58}
\epsilon
b_{\alpha\beta}(fa) = d_{\alpha\beta}(f)a 
+\epsilon f b_{\alpha\beta}(a)
\end{equation}
for $f 
\in \cO_{U_{\alpha\beta}}$ and $a \in \cO_{U_{\alpha\beta}}^{\oplus r}$.
Indeed we have
\begin{equation*}
\begin{aligned}
(\mathrm{id}+ \epsilon b_{\alpha\beta})(fa)
&=\Big(\bar{g}_{\alpha\beta}^{-1}  
(1\otimes \sigma_{\alpha}^*)^{-1}  g_{\alpha\beta}
  (1\otimes \sigma_{\beta}^*) \Big) (fa)  \\
&=\big(  (\sigma_{\alpha}^*)^{-1} \sigma_{\beta}^*(f) \big) \cdot  \Big(
\bar{g}_{\alpha\beta}^{-1} 
(1\otimes \sigma_{\alpha}^*)^{-1} g_{\alpha\beta}
 (1\otimes \sigma_{\beta}^*) \Big) (a)\\
&=\big(  f +d_{\alpha\beta}(f) \big) \cdot  (\mathrm{id}+ \epsilon b_{\alpha\beta}) (a)\\
&= fa +\big( d_{\alpha\beta}(f) a+ \epsilon f b_{\alpha\beta}(a)\big).
\end{aligned}
\end{equation*}
Since 
\begin{equation*}
\begin{aligned}
\bar{g}_{\alpha\beta} (\mathrm{id} + \epsilon b_{\alpha\beta}) &= 
(1 \otimes \sigma_{\alpha}^*)^{-1} g_{\alpha\beta} 
(1 \otimes \sigma_{\beta}^*)\\
&=(1 \otimes \sigma_{\alpha}^*)^{-1} g_{\alpha\beta} 
(1 \otimes \sigma_{\alpha}^*)
(1 \otimes \sigma_{\alpha}^*)^{-1}
(1 \otimes \sigma_{\beta}^*) \\
&=\bar{g}_{\alpha\beta} (\mathrm{id} + \epsilon g'_{\alpha\beta}) 
(1\otimes (1 +d_{\alpha\beta}))  \\
&= \bar{g}_{\alpha\beta} (\mathrm{id} + d_{\alpha\beta} + \epsilon g'_{\alpha\beta}),
\end{aligned}
\end{equation*}
we have $\epsilon b_{\alpha\beta} = d_{\alpha\beta} + \epsilon g'_{\alpha\beta} $ .
Since
\begin{equation*}
\begin{aligned}
u_{\alpha\beta} &=  \varphi^{-1}_{\alpha} \circ \varphi_{\beta}  - \mathrm{id}  \\
&= (\phi_{\alpha}^{-1} \circ (1\otimes \sigma_{\alpha}^*) \circ ( \bar{\phi}_{\alpha} \otimes 1))^{-1} \circ 
(\phi_{\beta}^{-1} \circ (1\otimes \sigma_{\beta}^*) \circ ( \bar{\phi}_{\beta} \otimes 1))
  - \mathrm{id} \\
  &= ( \bar{\phi}_{\alpha} \otimes 1)^{-1}  \circ 
(\bar{g}_{\alpha\beta} (\mathrm{id} + \epsilon b_{\alpha\beta}))
\circ ( \bar{\phi}_{\beta} \otimes 1)
  - \mathrm{id} \\
&  =( \bar{\phi}_{\beta} \otimes 1)^{-1} \circ
(\epsilon b_{\alpha\beta})  \circ ( \bar{\phi}_{\beta} \otimes 1)
\end{aligned}
\end{equation*}
and $\mathrm{symb}_1(u_{\alpha\beta})= 1\otimes d_{\alpha\beta} \in \mathrm{id} \otimes 
\Theta_{\cC_M/M} (-D(\tilde{\bt})_M )$,
we have $u_{\alpha\beta} \in \cA_{\tilde{E}}(D(\tilde{\bt}))|_{U_{\alpha\beta}}$.
Moreover we have
\begin{equation*}
v_{\alpha}:=  (\varphi_{\alpha}^{-1}\otimes \mathrm{id}) 
\circ \nabla_{\epsilon} |_{U^{\epsilon}_{\alpha}} \circ \varphi_{\alpha} 
- \tilde{\nabla}|_{U_{\alpha}\times \Spec \mathbb{C}[\epsilon]} .
\end{equation*}
We can check $v_{\alpha} 
\in \epsilon \otimes \cE nd (\tilde{E}) \otimes \Omega^1_{\cC/T} (D(\tilde{\bt}))$
as follows.
We have $v_{\alpha}(a)=0\text{ (mod $\epsilon$)}$ for $a \in \tilde{E}|_{U_{\alpha}}$, 
that is,
$v_{\alpha}(a) \in 
\epsilon \otimes (\tilde{E}\otimes \Omega^1_{\cC/T}(D(\tilde{\bt})))|_{U_{\alpha}}$.
Since
$d( \sigma_{\alpha}^*(f))=\sigma_{\alpha}^*(df)$ 
for $f \in \cO_{U_{\alpha}}$,
we have
\begin{equation*}
\begin{aligned}
v_{\alpha}(fa) &= ((\varphi_{\alpha}^{-1}\otimes \mathrm{id}) 
\circ \nabla_{\epsilon} |_{U^{\epsilon}_{\alpha}} \circ \varphi_{\alpha}) (fa)
-\tilde{\nabla}|_{U_{\alpha}\times \Spec \mathbb{C}[\epsilon]} (fa) \\
&= ((\varphi_{\alpha}^{-1}\otimes \mathrm{id}) 
\circ \nabla_{\epsilon} |_{U^{\epsilon}_{\alpha}} ) ( \sigma_{\alpha}^*(f)  \varphi_{\alpha} (a) )
-  a \otimes df  - f  \tilde{\nabla}|_{U_{\alpha}\times \Spec \mathbb{C}[\epsilon]} (a) \\
&= (\varphi_{\alpha}^{-1}\otimes \mathrm{id}) 
\circ \left(  \varphi_{\alpha}(a) \otimes d( \sigma_{\alpha}^*(f))+
\sigma_{\alpha}^*(f) \nabla_{\epsilon} |_{U^{\epsilon}_{\alpha}}   \varphi_{\alpha} (a) \right)
-  a \otimes df  - f  \tilde{\nabla}|_{U_{\alpha}\times \Spec \mathbb{C}[\epsilon]} (a) \\
&= f (\varphi_{\alpha}^{-1}\otimes \mathrm{id}) 
\circ  \nabla_{\epsilon} |_{U^{\epsilon}_{\alpha}} \circ   \varphi_{\alpha} (a) 
  - f  \tilde{\nabla}|_{U_{\alpha}\times \Spec \mathbb{C}[\epsilon]} (a)  = f v_{\alpha} (a).
\end{aligned}
\end{equation*}
Then we have $v_{\alpha} 
\in \epsilon \otimes \cE nd (\tilde{E}) \otimes \Omega^1_{\cC/T} (D(\tilde{\bt}))$.
Since $\varphi_{\alpha} |_{\tilde{t}_i\times \Spec \mathbb{C}[\epsilon]} (\tilde{l}_j^{(i)}
\otimes_{\mathbb{C}} \mathbb{C}[\epsilon]) 
= (l_{\epsilon})_j^{(i)} $ if $\tilde{t}_i |_{\cC_M}\cap U_{\alpha} \neq \emptyset$,
we have $v_{\alpha} 
\in \epsilon \otimes \mathcal{G}^1_M$.
We can check the equalities
\begin{equation}\label{2019.5.17.14.03}
u_{\beta\gamma} - u_{\alpha\gamma} + u_{\alpha\beta} =0 \quad  
\text{and} \quad \tilde{\nabla}\circ u_{\alpha\beta} - u_{\alpha\beta} \circ \tilde{\nabla} 
= v_{\beta} -v_{\alpha}.
\end{equation}
In fact, the first equality follows from 
the equality $g_{\alpha\beta} g_{\beta\gamma}= g_{\alpha \gamma}$.
Next, by the equality
\begin{equation}\label{2020.10.22.13.36}
(\varphi_{\beta}^{-1} \circ \varphi_{\alpha} ) \circ \varphi_{\alpha}^{-1}  \circ 
\nabla_{\epsilon} |_{U^{\epsilon}_{\alpha\beta}} \circ \varphi_{\alpha}  
\circ (\varphi_{\alpha}^{-1} \circ \varphi_{\beta} )
 =   
 \varphi_{\beta}^{-1} \circ
\nabla_{\epsilon} |_{U^{\epsilon}_{\alpha\beta}}
\circ \varphi_{\beta} 
\end{equation}
and the definition of $u_{\alpha\beta}$ and
$v_{\alpha}$, we have the equality
\begin{equation}\label{2020.10.22.13.37}
\begin{aligned}
(\mathrm{id} - u_{\alpha\beta} ) \circ 
(\tilde{\nabla}|_{U_{\alpha}\times \Spec \mathbb{C}[\epsilon]} + v_{\alpha})  
\circ (\mathrm{id} + u_{\alpha\beta} )
 &=   
\tilde{\nabla}|_{U_{\alpha\beta}\times \Spec \mathbb{C}[\epsilon]} + v_{\beta} \\ 
\tilde{\nabla}|_{U_{\alpha\beta}\times \Spec \mathbb{C}[\epsilon]} 
+ \left( v_{\alpha}
+ \tilde{\nabla}|_{U_{\alpha\beta}\times \Spec \mathbb{C}[\epsilon]}  \circ u_{\alpha\beta}
-u_{\alpha\beta} \circ \tilde{\nabla}|_{U_{\alpha\beta}\times \Spec \mathbb{C}[\epsilon]} \right)
 &=   
\tilde{\nabla}|_{U_{\alpha\beta}\times \Spec \mathbb{C}[\epsilon]} + v_{\beta}
\end{aligned}
\end{equation}
Then we have the second equality.
So $[(\{ u_{\alpha\beta} \},\{ v_{\alpha} \})]$ determines 
an element $\varsigma_M (v)$ of $\bH^1(\cG_M^{\bullet})$.

Conversely, by an element $[(\{ u_{\alpha\beta} \},\{ v_{\alpha} \})] \in \bH^1(\cG_M^{\bullet})$ 
we have a tangent vector field of $M$ as follows.
The class $[\{ \mathrm{symb}_1 (u_{\alpha\beta}) \}] \in H^1 (\Theta_{\cC_M/M} (-D(\tilde{\bt}) )_M)$ 
determines a first-order deformation $(\cC_{\epsilon}, D(\tilde{\bt})_{\epsilon})$ 
of $(\cC_M, D(\tilde{\bt})_M)$.
We can take an $M[\epsilon]$-morphism $\sigma_{\alpha}\colon U_{\alpha}^{\epsilon} 
\xrightarrow{\sim} U_{\alpha} \times \Spec \mathbb{C}[\epsilon]$
satisfying $\mathrm{symb}_1(u_{\alpha\beta})=
(\sigma_{\alpha}^*)^{-1} \circ \sigma_{\beta}^* - 1 $.
We put $\tilde{E}_{\alpha}:= 
(1\otimes \sigma^*_{\alpha}) (\tilde{E}|_{U_{\alpha}} \otimes_{\mathbb{C}} \mathbb{C}[\epsilon])$.
Let $\nabla_{\alpha}$ be the homomorphism of sheaves defined 
by the following composition
\begin{equation*}
\nabla_{\alpha} \colon \tilde{E}_{\alpha} 
\xrightarrow{(1\otimes \sigma^*_{\alpha})^{-1}} \tilde{E}|_{U_\alpha} 
\otimes_{\mathbb{C}} \mathbb{C}[\epsilon]
\xrightarrow{\tilde{\nabla}+ v_{\alpha}} (\tilde{E}|_{U_\alpha}  
\otimes \Omega^1_{\cC/T}(D(\tilde{\bt}))  )\otimes_{\mathbb{C}} \mathbb{C}[\epsilon] 
\xrightarrow{1\otimes \sigma_{\alpha}^*} \tilde{E}_{\alpha}
\otimes 
\Omega^1_{\cC_{\epsilon} / \Spec \cO_M[\epsilon]} \left(D(\tilde{\bt})_{\epsilon} \right).
\end{equation*}
Since for $f \in \mathcal{O}_{U_{\alpha}^{\epsilon} }$ 
and $a \in \tilde{E}_{\alpha} $, 
\begin{equation*}
\begin{aligned}
\nabla_{\alpha} (fa)&= (1\otimes \sigma_{\alpha}^*) \circ \Big( 
\tilde{\nabla}( (\sigma^*_{\alpha})^{-1}(f)  (1\otimes \sigma^*_{\alpha})^{-1}(a) )
+ v_{\alpha} ( (\sigma^*_{\alpha})^{-1}(f)  (1\otimes \sigma^*_{\alpha})^{-1}(a) )\Big)\\
&=(1\otimes \sigma_{\alpha}^*) \circ \Big( (1\otimes \sigma^*_{\alpha})^{-1}(a) \otimes
d((\sigma^*_{\alpha})^{-1}(f) )
+  (\sigma^*_{\alpha})^{-1}(f) \tilde{\nabla}
(   (1\otimes \sigma^*_{\alpha})^{-1}(a)) \\
&\qquad +  (\sigma^*_{\alpha})^{-1}(f)  v_{\alpha}
(   (1\otimes \sigma^*_{\alpha})^{-1}(a)) 
 \Big) \\
&=(1\otimes \sigma_{\alpha}^*) \circ  ((1\otimes \sigma^*_{\alpha})^{-1}(a) \otimes
((\sigma^*_{\alpha})^{-1}(df) ))  
+f (1\otimes \sigma_{\alpha}^*) \circ  (\tilde{\nabla}+ v_{\alpha})
(   (1\otimes \sigma^*_{\alpha})^{-1}(a))  \\
&=a\otimes d f +  f \nabla_{\alpha} (a),
\end{aligned}
\end{equation*}
the homomorphism $\nabla_{\alpha}$ becomes a connection.
Let $g_{\alpha\beta}$ be the composition
\begin{equation}\label{2020.5.18.10.27}
g_{\alpha\beta} \colon E_{\beta}|_{U^{\epsilon}_{\alpha\beta}} 
\xrightarrow{(1\otimes \sigma^*_{\beta})^{-1} |_{U^{\epsilon}_{\alpha\beta}}} \tilde{E}|_{U_{\alpha\beta}} 
\otimes_{\mathbb{C}} \mathbb{C}[\epsilon]
\xrightarrow{\mathrm{id}+u_{\alpha\beta}} \tilde{E}|_{U_{\alpha\beta}} 
\otimes_{\mathbb{C}} \mathbb{C}[\epsilon]
\xrightarrow{1\otimes \sigma^*_{\alpha} |_{U_{\alpha\beta} \times \Spec \mathbb{C}[\epsilon]}}
E_{\alpha}|_{U^{\epsilon}_{\alpha\beta}} .
\end{equation}
By the composition $g_{\alpha\beta}$ and the condition (\ref{2019.5.17.14.03}), we can
glue $(\tilde{E}_{\alpha}, \tilde{\nabla}_{\alpha})$.
For the subbundle $\tilde{l}_{j}^{(i)} \subset \tilde{E}_{(\tilde{t}_i)_S}$,
we put $(l_{\epsilon})_j^{(i)} :=
(1\otimes \sigma^*_{\alpha}) |_{(\tilde{t}_i)_S\times \Spec \mathbb{C}[\epsilon]} (\tilde{l}_j^{(i)}
\otimes_{\mathbb{C}} \mathbb{C}[\epsilon])$,
which gives a parabolic structure.
Then we obtain a flat family of parabolic connections on $(\cC_{\epsilon}, D(\tilde{\bt})_{\epsilon})$
over $\Spec \cO_M[\epsilon]$ up to isomorphism. 
This flat family gives a tangent vector field on $M$.
So $v \mapsto \varsigma_M (v)$ determines an isomorphism
\begin{equation*}
\varsigma_M \colon \Theta_{M_{\cC/T}^{\balpha}(\tilde{\bt} ,r,d)/N^{(n)}_r(d)} (M) 
\xrightarrow{\, \sim \, } \bH^1(\cG_M^{\bullet});
\quad v \longmapsto \varsigma_M(v).
\end{equation*}
The isomorphism $\varsigma_M$ induces a canonical isomorphism 
\begin{equation*}
\varsigma \colon \Theta_{M_{\cC/T}^{\balpha}(\tilde{\bt} ,r,d)/N^{(n)}_r(d)} 
\xrightarrow{\, \sim \, } \bold{R}^1(\pi_{M_{\cC/T}^{\balpha}(\tilde{\bt} ,r,d)})_*(\cG^{\bullet}).
\end{equation*}
\end{proof}

\subsection{Isomonodromic deformation}\label{SS IMD by c c}

Let $p_1 \colon  T\times N^{(n)}_r(d) \ra T$ be the projection.
There exists an algebraic splitting
\begin{equation}\label{eq1 of ss IMD}
D \colon  (p_1 \circ \varpi)^* (\Theta_T) \lra
 \Theta_{M_{\cC/T}^{\balpha}(\tilde{\bt}, r,d)/N^{(n)}_r(d)}
\end{equation}
of the tangent map 
$\Theta_{M_{\cC/T}^{\balpha}(\tilde{\bt}, r,d)/N^{(n)}_r(d)} \ra 
(p_1 \circ \varpi)^* (\Theta_T)$.
Here an image of (\ref{eq1 of ss IMD}) means an algebraic vector field determined 
by the isomonodromic deformation.
(See \cite[Proposition 8.1]{Inaba}).
We will define the algebraic splitting \eqref{eq1 of ss IMD}
rigorously below (see \eqref{2020.10.11.17.50}).
We describe this algebraic splitting in terms of the description 
of $\Theta_{M_{\cC/T}^{\balpha}(\tilde{\bt}, r,d)/N^{(n)}_r(d)}$ 
in Proposition \ref{Prop d t s of m by c c}.

Take any affine open set $U \subset T$ and a vector field $v \in H^0(U,\Theta_T)$.
Then $v$ corresponds to a morphism $\iota^v \colon \mathrm{Spec}\, \cO_U[\epsilon] \ra T$ 
with $\epsilon^2=0$ 
such that the composite $U \hookrightarrow \mathrm{Spec}\, \cO_U[\epsilon] \ra T$ 
is just the inclusion $U \hookrightarrow T$.
We denote the restriction of the universal family to 
$\cC \times_T (p_1 \circ \varpi)^{-1}(U)$ simply 
by $(\tilde{E},\tilde{\nabla},\{\tilde{l}^{(i)}_j \})$.
Consider the fiber product 
$\cC \times_T \mathrm{Spec}\, \cO_{(p_1 \circ \varpi)^{-1}(U)} [\epsilon]$ 
with respect to the canonical projection $\cC \ra T$ and 
the composite $\Spec\cO_{(p_1 \circ \varpi)^{-1}(U)}[\epsilon] 
\ra \Spec\cO_{U}[\epsilon] \xrightarrow{\iota^{v}}  T$.
We denote the pull-back of $D(\tilde{\bt})$ 
by the morphism $\cC \times_T \mathrm{Spec}\, 
\cO_{(p_1 \circ \varpi)^{-1}(U)} [\epsilon] \ra \cC$ simply 
by $D(\tilde{\bt})_{\cO_{(p_1 \circ \varpi)^{-1}(U)} [\epsilon]}$.
In this section, let $d$ be the relative exterior derivative 
on $\cC \times_T \mathrm{Spec}\, \cO_{(p_1 \circ \varpi)^{-1}(U)} [\epsilon]
\rightarrow (p_1 \circ \varpi)^{-1}(U)$.
We set 
$$ 
\begin{aligned}
&\hat{\Omega}^1 :=
\Omega^1_{\cC \times_T \mathrm{Spec}\, 
\cO_{(p_1 \circ \varpi)^{-1}(U)} [\epsilon]/(p_1 \circ \varpi)^{-1}(U)}\quad \text{and} \\
&\hat{\Omega}^2 :=
\Omega^2_{\cC \times_T \mathrm{Spec}\, 
\cO_{(p_1 \circ \varpi)^{-1}(U)} [\epsilon]/(p_1 \circ \varpi)^{-1}(U)}.
\end{aligned}
$$
Here $\hat{\Omega}^1$ is the sheaf of relative differential with respect to
the composition
$$
\cC \times_T \mathrm{Spec}\, \cO_{(p_1 \circ \varpi)^{-1}(U)} [\epsilon] \longrightarrow 
\mathrm{Spec}\, \cO_{(p_1 \circ \varpi)^{-1}(U)} [\epsilon] 
\longrightarrow (p_1 \circ \varpi)^{-1}(U)
$$
of the trivial projection.

\begin{Def}\label{Def horizontal lift}
We call $(\cE, \nabla^{\cE}, \{ (l_{\cE})^{(i)}_j \})$ \textit{a horizontal lift} 
of $(\tilde{E},\tilde{\nabla},\{\tilde{l}^{(i)}_j \})$ if 
\begin{itemize}
\item[(1)] $\cE$ is a vector bundle on 
$\cC \times_T \mathrm{Spec}\, \cO_{(p_1 \circ \varpi)^{-1}(U)} [\epsilon]$,
\item[(2)] $\cE|_{\tilde{t}_i \times \cO_{(p_1 \circ \varpi)^{-1}(U)} [\epsilon]} 
= (l_{\cE})^{(i)}_0 \supset \cdots \supset (l_{\cE})^{(i)}_r =0$ is a filtration 
by subbundles for $i=1,\ldots,n$,
\item[(3)] $\nabla^{\cE} \colon \cE \ra \cE \otimes 
\hat{\Omega}^1
\left(\log\left(D(\tilde{\bt})_{\cO_{(p_1 \circ \varpi)^{-1}(U)} [\epsilon]}\right) \right)$ 
is a connection satisfying
\begin{itemize}
\item[(a)] $\nabla^{\cE}(F^{(i)}_j(\cE) ) \subset F^{(i)}_j (\cE) \otimes 
\hat{\Omega}^1
\left(\log\left(D(\tilde{\bt})_{\cO_{(p_1 \circ \varpi)^{-1}(U)} [\epsilon]}\right) \right)$, 
where $F_j^{(i)} (\cE)$ is given by 
$F^{(i)}_j(\cE) := \mathrm{Ker}\left( \cE 
\ra \cE|_{\tilde{t}_i \times_T 
\cO_{(p_1 \circ \varpi)^{-1}(U)} [\epsilon]/
(p_1 \circ \varpi)^{-1}(U)}/(l^{\cE})^{(i)}_{j+1} \right)$,
\item[(b)] the curvature $\nabla^{\cE} \circ \nabla^{\cE} \colon \cE \ra \cE \otimes
\hat{\Omega}^2
\left(\log\left(D(\tilde{\bt})_{\cO_{(p_1 \circ \varpi)^{-1}(U)} [\epsilon]}\right) \right)$ is zero,
\item[(c)] $(\res_{\tilde{t}_i \times_T \Spec \cO_{(p_1 \circ \varpi)^{-1}(U)}[\epsilon]}
(\tilde{\nabla}^{\cE}) 
-  \tilde{\nu}^{(i)}_j) ((l^{\cE})^{(i)}_j) \subset (l^{\cE})^{(i)}_{j+1}$ for any $i, j$, where 
$\tilde{\nabla}^{\cE}$ is the relative connection over 
$\Spec \cO_{(p_1 \circ \varpi)^{-1}(U)}[\epsilon]$ 
induced by $\nabla^{\cE}$ and 
\item[(d)] $(\cE, \tilde{\nabla}^{\cE}, \{ (l_{\cE})^{(i)}_j \}) 
\otimes \cO_{(p_1 \circ \varpi)^{-1}(U)}[\epsilon]/(\epsilon) 
\cong (\tilde{E},\tilde{\nabla},\{\tilde{l}^{(i)}_j \})$.
\end{itemize}
\end{itemize}
Here, we define the sheaf $\hat{\Omega}^1
\left(\log\left(D(\tilde{\bt})_{\cO_{(p_1 \circ \varpi)^{-1}(U)} [\epsilon]}\right) \right)$ 
as the coherent subsheaf of $\hat{\Omega}^1
 \left(D(\tilde{\bt})_{\cO_{(p_1 \circ \varpi)^{-1}(U)} [\epsilon]}\right)$ 
locally generated by
$\tilde{g}^{-1} d\tilde{g}$ and $d \epsilon$ for a local defining equation $\tilde{g}$ 
of $D(\tilde{\bt})_{\cO_{(p_1 \circ \varpi)^{-1}(U)} [\epsilon]}$
and the sheaf $\hat{\Omega}^2
\left(\log\left(D(\tilde{\bt})_{\cO_{(p_1 \circ \varpi)^{-1}(U)} [\epsilon]}\right) \right)$
as the coherent subsheaf of $\hat{\Omega}^2 
\left(D(\tilde{\bt})_{\cO_{(p_1 \circ \varpi)^{-1}(U)} [\epsilon]}\right)$ 
locally generated by
$\tilde{g}^{-1} d\tilde{g} \wedge d \epsilon$.
\end{Def}

Let $\cC \times_T (p_1 \circ \varpi)^{-1}(U)=\bigcup_{\alpha}U_{\alpha}$ 
be an affine open covering such that 
we have $\bar{\phi}_{\alpha}\colon \tilde{E}|_{U_{\alpha}} 
\xrightarrow{\sim} \cO_{U_{\alpha}}^{\oplus r}$ for any $\alpha$,
$\sharp\{ i \mid \tilde{t}_i |_{\cC \times_T (p_1 \circ \varpi)^{-1}(U)} \cap U_{\alpha}
 \neq \emptyset \} \le 1$ 
for any $\alpha$ and 
$\sharp\{ \alpha \mid \tilde{t}_i |_{\cC \times_T (p_1 \circ \varpi)^{-1}(U)} \cap U_{\alpha} 
\neq \emptyset \} \le 1$ 
for any $i$.
Assume that the parabolic connection $(\tilde{E}, \tilde{\nabla}, \{ \tilde{l}_{j}^{(i)} \})$ 
is locally given in the affine subset $U_{\alpha}$ by a connection matrix
$A_{\alpha} f_{\alpha}^{-1} df_{\alpha}$, where 
\begin{itemize}
\item $f_{\alpha}$ is a local defining equation of 
$(\tilde{t}_i \times (p_1 \circ \varpi)^{-1}(U))\cap U_{\alpha}$, 
\item $A_{\alpha} \in M_r(\cO_{U_{\alpha}})$,
\item $A_{\alpha}((\tilde{t}_i \times_T (p_1 \circ \varpi)^{-1}(U))\cap U_{\alpha})$ 
is an upper triangular matrix, and 
\item the parabolic structure $\{ \tilde{l}^{(i)}_j \}_{U_{\alpha}}$ is given 
by $(\tilde{l}_j^{(i)})_{U_{\alpha}}=(*,*,\ldots,*,0,\ldots,0)$.
\end{itemize}
Put $U_{\alpha}^{\epsilon}:= 
U_{\alpha} \times_T\Spec \cO_{(p_1 \circ \varpi)^{-1}(U)}[\epsilon]$.
We take an $M[\epsilon]$-morphism $\sigma_{\alpha} \colon U^{\epsilon}_{\alpha} 
\ra U_{\alpha}\times \Spec \mathbb{C}[\epsilon]$
which is a lift of $\mathrm{id}_{U_{\alpha}}$ 
preserving the divisor $D(\tilde{\bt})_{\cO_{(p_1 \circ \varpi)^{-1}(U)} [\epsilon]} 
\cap U^{\epsilon}_{\alpha}$ 
and $(D(\tilde{\bt})\cap U_{\alpha})\times \Spec \mathbb{C}[\epsilon]$.
Put $f_{\alpha}^{\epsilon}:= \sigma_{\alpha}^* f_{\alpha}$, $\tilde{E}_{\alpha}:= 
(1\otimes \sigma^*_{\alpha}) (\tilde{E}|_{U_{\alpha}} \otimes_{\mathbb{C}} \mathbb{C}[\epsilon])$,
and $A^{\epsilon}_{\alpha} := (1\otimes \sigma^*_{\alpha}) (A_{\alpha}) 
\in M_r (\cO_{U_{\alpha}^{\epsilon}})$.
If we denote the composite
\begin{equation*}
\cO_{U_{\alpha}^{\epsilon}} \xrightarrow{\, d \,} 
\Omega^1_{U_{\alpha}^{\epsilon} / (p_1 \circ \varpi)^{-1}(U)}
= \cO_{U_{\alpha}^{\epsilon}}df_{\alpha}^{\epsilon} \oplus \cO_{U_{\alpha}^{\epsilon}} d \epsilon 
\lra \cO_{U_{\alpha}^{\epsilon}} d \epsilon
\end{equation*}
by $d_{\epsilon}$, then we have $d_{\epsilon}(A^{\epsilon}_{\alpha})=0$ (which means integrable).
Then $\tilde{E}_{\alpha}$ 
and the connection matrix $A^{\epsilon}_{\alpha}(f^{\epsilon}_{\alpha})^{-1}d f^{\epsilon}_{\alpha}$ 
give a \textit{local} horizontal lift of $(\tilde{E}, \tilde{\nabla}, \{ \tilde{l}_{j}^{(i)} \})|_{U_{\alpha}}$.
By \cite[Proposition 8.1]{Inaba},
the obstruction for the patching the local horizontal lifts vanishes and 
the global horizontal lift of $(\tilde{E},\tilde{\nabla},\{\tilde{l}^{(i)}_j \})$ 
to $\cC \times_T (p_1 \circ \varpi)^{-1}(U)$ 
is unique for a vector field $v \in H^0(U,\Theta_T)$.
Let $M$ be an affine open subset of $M_{\cC/T}^{\balpha}(\tilde{\bt}, r,d)$.
If we have $v \in H^0(M,  (p_1 \circ \varpi)^* \Theta_T)$, 
then we have the relative connection 
$(E^{\text{hor}}, \tilde{\nabla}^{\text{hor}}, \{ (\tilde{l}^{\text{hor}})^{(i)}_j \})$
on $\cC\times_{T \times N^{(n)}_r(d)} M[\epsilon]$ over $M[\epsilon]$ 
induced by the global horizontal lift $(\cE, \nabla^{\cE}, \{ (l_{\cE})^{(i)}_j \})$ 
on $\cC\times_{T \times N^{(n)}_r(d)} M[\epsilon]$
with respect to $v$.
We denote by $D(v) \in  \Theta_{M_{\cC/T}^{\balpha}(\tilde{\bt}, r,d)/N^{(n)}_r(d)}(M)$ 
the vector filed on $M$ corresponding to the 
relative connection $(E^{\text{hor}}, \tilde{\nabla}^{\text{hor}}, \{ (\tilde{l}^{\text{hor}})^{(i)}_j \})$.
We define a morphism $D$ as 
\begin{equation}\label{2020.10.11.17.50}
\begin{aligned}
D \colon  (p_1 \circ \varpi)^* (\Theta_T) &\lra
 \Theta_{M_{\cC/T}^{\balpha}(\tilde{\bt}, r,d)/N^{(n)}_r(d)} \\
v &\longmapsto D(v),
\end{aligned}
\end{equation}
which is the algebraic splitting \eqref{eq1 of ss IMD}.

\begin{Rem}
The connection $\nabla^{\cE}$ on $\cE$ on $\cC\times_T\Spec 
\cO_{(p_1 \circ \varpi)^{-1}(U)}[\epsilon]$ 
satisfies the \textit{integrability condition} (Definition \ref{Def horizontal lift} (3) (b)).
The integrability means that the relative connection associated to $\nabla^{\cE}$ is 
an \textit{isomonodromic family}. 
(See for example \cite[0.16.6]{Sab}).
\end{Rem}

Let
\begin{equation*}
\mu  \colon (p_1 \circ \varpi)^* \Theta_T  \lra  
\bold{R}^1 (\pi_{M_{\cC/T}^{\balpha}(\tilde{\bt} ,r,d)})_*(\Theta_{\cC \times_T M_{\cC/T}^{\balpha}
(\tilde{\bt},r,d)/ M_{\cC/T}^{\balpha}(\tilde{\bt},r,d)} (-D(\tilde{\bt}))) 
\end{equation*}
be the Kodaira--Spencer map, 
where $\pi_{M_{\cC/T}^{\balpha}(\tilde{\bt} ,r,d)} \colon \cC_{M_{\cC/T}^{\balpha}(\tilde{\bt} ,r,d)} 
\rightarrow M_{\cC/T}^{\balpha}(\tilde{\bt} ,r,d)$ 
is the natural morphism.
We obtain the desired algebraic splitting $D$ as follows.
(The statement of the following proposition is essentially pointed out 
in \cite[Section 4]{BHH}).
\begin{Prop}\label{Prop explicit description IMD}
\textit{Let $M$ be an affine open subset of $M^{\balpha}_{\cC/T}(\tilde{\bt},r,d)$.}
\textit{The following morphism}
\begin{equation}\label{map which means IMD}
\begin{aligned}
(p_1 \circ \varpi)^* \Theta_T (M) &\lra \bH^1(\cG^{\bullet}_M) 
\cong \Theta_{M_{\cC/T}^{\balpha}(\tilde{\bt} ,r,d)/N^{(n)}_r(d)}(M) \\
v&\longmapsto [\{\iota(\tilde{\nabla}) \left( \mu_M(v) \right) \}, \{ 0 \}]
\end{aligned}
\end{equation}
\textit{coincides with the algebraic splitting \eqref{2020.10.11.17.50}.}
\textit{Here}
\begin{equation*}
\iota(\tilde{\nabla}) \colon  H^1(\Theta_{\cC\times_T M/M}(-D(\tilde{\bt})_M)) \lra H^1(\cG^{0}_M) 
\end{equation*}
\textit{is induced by the splitting associated to the universal family $\tilde{\nabla}$
defined in Section} \ref{SS AA}.
\end{Prop}
\begin{proof}
Take an affine open set $U \subset T$ and a vector field $v \in H^0(U,\Theta_T)$.
We denote the restriction of the universal family 
to $\cC \times_T (p_1 \circ \varpi)^{-1}(U)$ simply by $(\tilde{E},\tilde{\nabla},\{\tilde{l}^{(i)}_j \})$.
We take a horizontal lift $(\cE, \nabla^{\cE}, \{ (l_{\cE})^{(i)}_j \})$ 
of $(\tilde{E},\tilde{\nabla},\{\tilde{l}^{(i)}_j \})$ corresponding to $v$.

We take an affine open set $M\subset (p_1 \circ \varpi)^{-1}(U)$ and 
put $M[\epsilon] = \mathrm{Spec}\, \cO_{M} [\epsilon]$.
We denote the restriction of the horizontal lift $(\cE, \nabla^{\cE}, \{ (l_{\cE})^{(i)}_j \})$ 
by $(\cE_M, \nabla_M^{\cE}, \{ ((l_{\cE})_M)^{(i)}_j \})$.
As in the proof of Proposition \ref{Prop d t s of m by c c},
take an affine open covering $\{ U_{\alpha} \}$ of $\cC \times_T M$.
Let $\{ U_{\alpha}^{\epsilon} \}$ be the affine open covering 
of $\cC \times_T M[\epsilon]$ corresponding to $\{ U_{\alpha} \}$ and 
let $\sigma_{\alpha} \colon U^{\epsilon}_{\alpha} \ra U_{\alpha}\times \Spec \mathbb{C}[\epsilon]$ 
be an isomorphism as in the proof of Proposition \ref{Prop d t s of m by c c}.
Put $d_{\alpha\beta} := (\sigma_{\alpha}^*)^{-1}\circ \sigma_{\beta}^*-1 
\colon \cO_{U_{\alpha\beta}} \ra \epsilon \otimes \cO_{U_{\alpha\beta}}$.
Then $[\{ d_{\alpha\beta} \}]$
is the Kodaira--Spencer class corresponding to $\cC \times_T M[\epsilon] \ra M[\epsilon]$.
If we take a frame $\phi_{\alpha}\colon \cE_M|_{U_{\alpha}^{\epsilon}} 
\rightarrow \cO_{U_{\alpha}^{\epsilon}}^{\oplus r}$ and 
put $\bar{\phi}_{\alpha}:= \phi_{\alpha}\ (\text{mod } \epsilon)$, there is a commutative diagram 
\begin{equation*}
\xymatrix@C=36pt{
\tilde{E}|_{U_{\alpha\beta}} \otimes_{\mathbb{C}} \mathbb{C}[\epsilon] 
\ar[r]^-{\bar{\phi}_{\beta}\otimes 1} 
\ar[d]^-{\mathrm{id}+ u_{\alpha\beta}} & 
\cO^{\oplus r}_{U_{\alpha\beta}} \otimes_{\cO_{U_{\alpha\beta}}} 
\cO_{U_{\alpha\beta}}[\epsilon] 
\ar[r]^-{1 \otimes \sigma_{\beta}^*} 
\ar[d]^-{ (\bar{\phi}_{\alpha}\otimes1)  (\bar{\phi}_{\beta}\otimes1)^{-1} 
(\mathrm{id}+ \epsilon b_{\alpha\beta})} &
\cO^{\oplus r}_{U_{\alpha\beta}} \otimes_{\cO_{U_{\alpha\beta}}} \cO_{U^{\epsilon}_{\alpha\beta}} 
\ar[r]^-{\phi_{\beta}^{-1}} \ar[d]^-{\phi_{\alpha}  \phi_{\beta}^{-1}} &
\cE_M|_{U^{\epsilon}_{\alpha\beta}} \ar[d]^-{\mathrm{id}} \\
\tilde{E}|_{U_{\alpha\beta}} \otimes_{\mathbb{C}} \mathbb{C}[\epsilon] 
\ar[r]^-{\bar{\phi}_{\alpha}\otimes 1} & 
\cO^{\oplus r}_{U_{\alpha\beta}} \otimes_{\cO_{U_{\alpha\beta}}}
 \cO_{U_{\alpha\beta}}[\epsilon] 
\ar[r]^-{1 \otimes \sigma_{\alpha}^*} &
\cO^{\oplus r}_{U_{\alpha\beta}} \otimes_{\cO_{U_{\alpha\beta}}} \cO_{U^{\epsilon}_{\alpha\beta}} 
\ar[r]^-{\phi_{\alpha}^{-1}} &
\cE_M|_{U^{\epsilon}_{\alpha\beta}}\rlap{.}
}
\end{equation*}
Here $b_{\alpha\beta} \colon \cO^{\oplus r}_{U_{\alpha\beta}} \ra 
\cO^{\oplus r}_{U_{\alpha\beta}}$ 
is a differential operator of degree $\leq 1$ 
satisfying $\epsilon
b_{\alpha\beta}(fa) = d_{\alpha\beta}(f)a 
+\epsilon f b_{\alpha\beta}(a)$ for 
$f \in \cO_{U_{\alpha\beta}}$ and $a \in \cO_{U_{\alpha\beta}}^{\oplus r}$.
If we denote the connection matrix of $\nabla_M^{\cE}|_{U_{\alpha}^{\epsilon}}$ 
via the frame $\phi_{\alpha}$ by
$A^{\epsilon}_{\alpha}(f^{\epsilon}_{\alpha})^{-1}d f^{\epsilon}_{\alpha}$,
then the connection matrix 
of $\tilde{\nabla}|_{U_{\alpha}} \otimes_{\mathbb{C}} \mathbb{C}[\epsilon]$ becomes
$(1\otimes\sigma_{\alpha}^*)^{-1}(A^{\epsilon}_{\alpha}(f^{\epsilon}_{\alpha})^{-1}
d f^{\epsilon}_{\alpha})(1\otimes\sigma_{\alpha}^*)$.
The patching condition 
for the connections $\tilde{\nabla}|_{U_{\alpha}} \otimes_{\mathbb{C}} \mathbb{C}[\epsilon]$ 
becomes
\begin{equation}\label{Pat. Cond. conn}
\begin{aligned}
&\bar{\phi}_{\beta} \bar{\phi}_{\alpha}^{-1} d (\bar{\phi}_{\alpha} \bar{\phi}_{\beta}^{-1}) 
\bar{\phi}_{\beta} \bar{\phi}_{\alpha}^{-1}
+\bar{\phi}_{\beta} \bar{\phi}_{\alpha}^{-1} (1\otimes\sigma_{\alpha}^*)^{-1}
\left(A^{\epsilon}_{\alpha}\frac{d f^{\epsilon}_{\alpha}}{f^{\epsilon}_{\alpha}}\right)
(1\otimes\sigma_{\alpha}^*)\\
&= (1\otimes\sigma_{\beta}^*)^{-1}
\left(A^{\epsilon}_{\beta}\frac{d f^{\epsilon}_{\beta}}{f^{\epsilon}_{\beta}}\right) 
(1\otimes\sigma_{\beta}^*)
\bar{\phi}_{\beta} \bar{\phi}_{\alpha}^{-1}.
\end{aligned}
\end{equation}
Here we denote $\bar{\phi}_{\alpha}\otimes 1$ simply by $\bar{\phi}_{\alpha}$.
Moreover the patching condition for $\nabla_M^{\cE}|_{U_{\alpha}^{\epsilon}}$ is the equality
\begin{equation}\label{Pat. Cond. horiz}
(\phi_{\alpha}  \phi_{\beta}^{-1}) A^{\epsilon}_{\beta}\frac{d f^{\epsilon}_{\beta}}{f^{\epsilon}_{\beta}}
=d(\phi_{\alpha}  \phi_{\beta}^{-1})
+A^{\epsilon}_{\alpha}
\frac{d f^{\epsilon}_{\alpha}}{f^{\epsilon}_{\alpha}}(\phi_{\alpha}  \phi_{\beta}^{-1}).
\end{equation}
For $a \in \cO^{\oplus r}_{U_{\alpha\beta}} 
\otimes_{\cO_{U_{\alpha\beta}}} \cO_{U^{\epsilon}_{\alpha\beta}} $,
we have
\begin{equation*}
\begin{aligned}
&d (\phi_{\alpha} \phi_{\beta}^{-1}(a))\\
&= 
d \left( (1 \otimes \sigma_{\alpha}^*) \bar{\phi}_{\alpha}\bar{\phi}_{\beta}^{-1}  
(1+ \epsilon b_{\alpha\beta}) (1 \otimes \sigma_{\beta}^*)^{-1} (a) \right) \\
&=  (1 \otimes \sigma_{\alpha}^*) d \left( \bar{\phi}_{\alpha} \bar{\phi}_{\beta}^{-1} \right)
 (1+ \epsilon b_{\alpha\beta}) (1 \otimes \sigma_{\beta}^*)^{-1} (a)
 +(1 \otimes \sigma_{\alpha}^*) \bar{\phi}_{\alpha} \bar{\phi}_{\beta}^{-1} 
 d \Big((1+ \epsilon b_{\alpha\beta}) (1 \otimes \sigma_{\beta}^*)^{-1}(a) \Big)  \\
&=  (1 \otimes \sigma_{\alpha}^*) d \left( \bar{\phi}_{\alpha} \bar{\phi}_{\beta}^{-1} \right) 
\bar{\phi}_{\beta} \bar{\phi}_{\alpha}^{-1} 
(1 \otimes \sigma_{\alpha}^*)^{-1} \phi_{\alpha} \phi_{\beta}^{-1} (a)
 +(1 \otimes \sigma_{\alpha}^*) \bar{\phi}_{\alpha} \bar{\phi}_{\beta}^{-1} 
 d \Big((1+ \epsilon b_{\alpha\beta}) (1 \otimes \sigma_{\beta}^*)^{-1}(a) \Big) .
\end{aligned}
\end{equation*}
Set $a':= (1 \otimes \sigma_{\beta}^*)^{-1}(a)$.
Rephrasing the above equality and substituting (\ref{Pat. Cond. conn}), 
(\ref{Pat. Cond. horiz}), we have
\begin{equation}\label{2019.5.20.16.52}
\begin{aligned}
&d ((1+ \epsilon b_{\alpha\beta})(a')) \\
&=\bar{\phi}_{\beta} \bar{\phi}_{\alpha}^{-1} 
(1 \otimes \sigma_{\alpha}^*)^{-1} d \Big((\phi_{\alpha}  \phi_{\beta}^{-1}) 
(1 \otimes \sigma_{\beta}^*) (a')\Big)
 - \bar{\phi}_{\beta} \bar{\phi}_{\alpha}^{-1} d \left( \bar{\phi}_{\alpha} \bar{\phi}_{\beta}^{-1} \right) 
\bar{\phi}_{\beta} \bar{\phi}_{\alpha}^{-1} 
(1 \otimes \sigma_{\alpha}^*)^{-1} \phi_{\alpha} \phi_{\beta}^{-1} (1 \otimes \sigma_{\beta}^*)(a') \\
&=\bar{\phi}_{\beta} \bar{\phi}_{\alpha}^{-1} (1 \otimes \sigma_{\alpha}^*)^{-1} 
\phi_{\alpha}  \phi_{\beta}^{-1} A^{\epsilon}_{\beta}\frac{d f^{\epsilon}_{\beta}}{f^{\epsilon}_{\beta}} 
(1 \otimes \sigma_{\beta}^*)(a') 
 -\bar{\phi}_{\beta} \bar{\phi}_{\alpha}^{-1} (1 \otimes \sigma_{\alpha}^*)^{-1} 
A^{\epsilon}_{\alpha}\frac{d f^{\epsilon}_{\alpha}}{f^{\epsilon}_{\alpha}}
(\phi_{\alpha}  \phi_{\beta}^{-1}) (1 \otimes \sigma_{\beta}^*) (a') \\
&\qquad
+\bar{\phi}_{\beta} \bar{\phi}_{\alpha}^{-1} 
(1 \otimes \sigma_{\alpha}^*)^{-1} (\phi_{\alpha}  \phi_{\beta}^{-1}) 
d \Big((1 \otimes \sigma_{\beta}^*) (a') \Big)
 - \bar{\phi}_{\beta} \bar{\phi}_{\alpha}^{-1} d \left( \bar{\phi}_{\alpha} \bar{\phi}_{\beta}^{-1} \right) 
\bar{\phi}_{\beta} \bar{\phi}_{\alpha}^{-1} 
(1 \otimes \sigma_{\alpha}^*)^{-1} \phi_{\alpha} \phi_{\beta}^{-1} (1 \otimes \sigma_{\beta}^*) (a') \\
&=\bar{\phi}_{\beta} \bar{\phi}_{\alpha}^{-1} (1 \otimes \sigma_{\alpha}^*)^{-1} 
\phi_{\alpha}  \phi_{\beta}^{-1} A^{\epsilon}_{\beta}\frac{d f^{\epsilon}_{\beta}}{f^{\epsilon}_{\beta}}
 (1 \otimes \sigma_{\beta}^*) (a')
 - (1\otimes\sigma_{\beta}^*)^{-1}
 \left(A^{\epsilon}_{\beta}\frac{d f^{\epsilon}_{\beta}}{f^{\epsilon}_{\beta}}\right) 
 (1\otimes\sigma_{\beta}^*)
 \bar{\phi}_{\beta} \bar{\phi}_{\alpha}^{-1}
(1 \otimes \sigma_{\alpha}^*)^{-1} \phi_{\alpha} \phi_{\beta}^{-1} (1 \otimes \sigma_{\beta}^*) 
(a')\\
&\qquad +\bar{\phi}_{\beta} \bar{\phi}_{\alpha}^{-1} 
(1 \otimes \sigma_{\alpha}^*)^{-1} (\phi_{\alpha}  \phi_{\beta}^{-1}) 
(1 \otimes \sigma_{\beta}^*) (d a') \\
&= (1+ \epsilon b_{\alpha\beta})
 \left( (1\otimes\sigma_{\beta}^*)^{-1} 
 A^{\epsilon}_{\beta}\frac{d f^{\epsilon}_{\beta}}{f^{\epsilon}_{\beta}} 
 (1\otimes\sigma_{\beta}^*) \right)(a')
-\left(  (1\otimes\sigma_{\beta}^*)^{-1}A^{\epsilon}_{\beta}
\frac{d f^{\epsilon}_{\beta}}{f^{\epsilon}_{\beta}} (1\otimes\sigma_{\beta}^*) \right)
 (1+ \epsilon b_{\alpha\beta}) (a')
+ (1+ \epsilon b_{\alpha\beta})(da').
\end{aligned}
\end{equation}
In the last line, $b_{\alpha\beta}$ acts on the elements of 
$(\mathcal{E}nd(\mathcal{O}_{U_{\alpha\beta}}^{\oplus r}) df_{\beta}/f_{\beta})
\oplus \mathcal{E}nd(\mathcal{O}_{U_{\alpha\beta}}^{\oplus r})  d \epsilon$.
This action is defined as follows.
Let $\tilde{d}_{\alpha\beta} \colon \mathcal{O}_{U_{\alpha\beta}} 
\rightarrow \mathcal{O}_{U_{\alpha\beta}}$
be the homomorphism of sheaves such that the composition
$\epsilon\circ \tilde{d}_{\alpha\beta}\colon \mathcal{O}_{U_{\alpha\beta}} 
\rightarrow \mathcal{O}_{U_{\alpha\beta}} \rightarrow 
\epsilon \otimes\mathcal{O}_{U_{\alpha\beta}}$
coincides with $d_{\alpha\beta}$.
By the equality \eqref{2020.5.22.21.58},
we have $b_{\alpha\beta}=\tilde{d}_{\alpha\beta}+g'_{\alpha\beta}$.
Let $X_0 df_{\beta}/f_{\beta} +X_1 d \epsilon$
be an element of
$(\mathcal{E}nd(\mathcal{O}_{U_{\alpha\beta}}^{\oplus r}) df_{\beta}/f_{\beta})
\oplus \mathcal{E}nd(\mathcal{O}_{U_{\alpha\beta}}^{\oplus r})  d \epsilon$.
Let $d_{f_{\beta}}$ be the relative exterior derivative
$\mathcal{E}nd(\mathcal{O}_{U_{\beta}}^{\oplus r})
\rightarrow
\mathcal{E}nd(\mathcal{O}_{U_{\beta}}^{\oplus r}) df_{\beta}$.
The operator $b_{\alpha\beta}$ acts on 
$X_0 df_{\beta}/f_{\beta} +X_1 d \epsilon$
as follows:
\begin{equation*}
\begin{aligned}
b_{\alpha\beta} (X_0 \frac{df_{\beta}}{f_{\beta}} +X_1 d \epsilon)
&=
(\tilde{d}_{\alpha\beta}  + g'_{\alpha\beta} ) (X_0 \frac{df_{\beta}}{f_{\beta}} 
+X_1 d \epsilon)\\
&= d_{f_{\beta}}\left\langle   \tilde{d}_{\alpha\beta} , X_0 \frac{ df_{\beta}}{f_{\beta}}
\right\rangle 
+ (g_{\alpha\beta}' \otimes 1) X_0 \frac{df_{\beta}}{f_{\beta} }
 +\left(
X_0\left\langle   \tilde{d}_{\alpha\beta} , \frac{ df_{\beta}}{f_{\beta}}
\right\rangle+ g_{\alpha\beta}'X_1 
\right)d\epsilon.
\end{aligned}
\end{equation*}
This action is natural since $d_{\alpha\beta}= (\sigma_{\alpha}^*)^{-1}  \sigma_{\beta}^*-1$ and
\begin{equation}\label{2019.5.21.11.36}
\begin{aligned}
(\sigma_{\alpha}^*)^{-1}  \sigma_{\beta}^*  \left(
(\sigma_{\beta}^*)^{-1} \frac{d f^{\epsilon}_{\beta}}{f^{\epsilon}_{\beta}}\right)
&=(\sigma_{\alpha}^*)^{-1}  \sigma_{\beta}^*
\left( \frac{df_{\beta}}{f_{\beta}} \right)
= \frac{d ( f_{\beta}+d_{\alpha\beta} (f_{\beta}))}{f_{\beta}+d_{\alpha\beta} (f_{\beta})} \\
&= \frac{df_{\beta}}{f_{\beta}} +
\frac{ (\tilde{d}_{\alpha\beta} (f_{\beta}))}{f_{\beta}} d\epsilon
+\epsilon\frac{  d(\tilde{d}_{\alpha\beta} (f_{\beta}))}{f_{\beta}}
- \epsilon\frac{  \tilde{d}_{\alpha\beta} (f_{\beta}) d f_{\beta}}{(f_{\beta})^2}\\
&= \frac{df_{\beta}}{f_{\beta}} + 
 \left\langle   \tilde{d}_{\alpha\beta} ,  \frac{ df_{\beta}}{f_{\beta}}
\right\rangle d\epsilon+
\epsilon d_{f_{\beta}}\left\langle   \tilde{d}_{\alpha\beta} , \frac{ df_{\beta}}{f_{\beta}}
\right\rangle ,
\end{aligned}
\end{equation}
where $d_{f_{\beta}}$ is the relative exterior derivative
$\mathcal{O}_{U_{\beta}}
\rightarrow \mathcal{O}_{U_{\beta}} df_{\beta}$.
Here the last line of \eqref{2019.5.21.11.36} follows from
 the equalities $\epsilon^2=0$ and $2 \epsilon d\epsilon= d(\epsilon^2)=0$.

We compute the $d \epsilon$ terms in the first line 
and the last line of the equations \eqref{2019.5.20.16.52}. 
By the definition of $f^{\epsilon}_{\beta}$, 
we have  $(\sigma_{\beta}^*)^{-1} (f^{\epsilon}_{\beta} )= f_{\beta} $.
We define $a'_0 \in \mathcal{O}_{U_{\alpha\beta}}^{\oplus r}$ 
and $a'_1 \in \mathcal{O}_{U_{\alpha\beta}}^{\oplus r}$ 
by $a' = a'_0+\epsilon a'_1$.
The $d \epsilon$ term in $d ((1+ \epsilon b_{\alpha\beta})(a'))$ is 
$(b_{\alpha\beta}(a'_0)+ a_1')d\epsilon$.
The $d \epsilon$ term in $(1+ \epsilon b_{\alpha\beta})(da')$ is 
$(\tilde{d}_{\alpha\beta} (a'_0) + a_1')d\epsilon$
by direct computation as in the equalities \eqref{2019.5.21.11.36}.
By comparing the $d \epsilon$ term in the first line of the equations \eqref{2019.5.20.16.52}
with the $d \epsilon$ term in the last line of the equations \eqref{2019.5.20.16.52},
we have
\begin{equation*}
b_{\alpha\beta}(a_0')d \epsilon +a_1'd\epsilon=
\left< \tilde{d}_{\alpha\beta} ,  
A_{\beta}\frac{d f_{\beta}}{f_{\beta}}(a_0') \right>d \epsilon
+   \tilde{d}_{\alpha\beta} (a'_0) d \epsilon 
+a_1'd\epsilon,
\end{equation*}
where $A_{\beta} f_{\beta}^{-1}d f_{\beta}$ is the connection matrix of $\tilde{\nabla}$ 
via the frame $\bar{\phi}_{\beta}$.
Then we obtain
\begin{equation}\label{horiz u}
u_{\alpha\beta}= \bar{\phi}_{\beta}^{-1}\circ \left( d_{\alpha\beta} + 
\left< d_{\alpha\beta} ,  A_{\beta}\frac{d f_{\beta}}{f_{\beta}} \right> \right)
\circ \bar{\phi}_{\beta}.
\end{equation}
On the other hand, we take the relative connection $\overline{\nabla}_M^{\cE}$ on $\cE_M$ 
associated to $\nabla_M^{\cE}$.
Since $d_{\epsilon}(A^{\epsilon}_{\alpha})=0$, we have
\begin{equation}\label{horiz v}
v_{\alpha}=  (\varphi_{\alpha}^{-1}\otimes \mathrm{id}) \circ 
\overline{\nabla}_M^{\cE} |_{U^{\epsilon}_{\alpha}} \circ \varphi_{\alpha} 
- \tilde{\nabla}|_{U_{\alpha}\times \Spec \mathbb{C}[\epsilon]} =0.
\end{equation}
So for $v \in H^0(U,\Theta_T)$, we obtain the element of $\bH^1(\cG^{\bullet}_M)$ 
given by \eqref{horiz u} and \eqref{horiz v}.
This correspondence from $v$ to this element of $\bH^1(\cG^{\bullet}_M)$ is just 
the morphism $D$ defined in \eqref{2020.10.11.17.50}.
On the other hand, by the local description of $\iota(\tilde{\nabla})$ in Section \ref{SS AA},
we have the following equality:
$$
\iota(\tilde{\nabla}) ([\{d_{\alpha\beta}\}]) 
=\bar{\phi}_{\beta}^{-1}\circ \left( d_{\alpha\beta} + 
\left< d_{\alpha\beta} ,  A_{\beta}\frac{d f_{\beta}}{f_{\beta}} \right> \right)
\circ \bar{\phi}_{\beta}.
$$
This equality means that the morphism (\ref{map which means IMD}) coincides with
the morphism
$D \colon  (p_1 \circ \varpi)^* (\Theta_T) \ra 
\Theta_{M_{\cC/T}^{\balpha}(\tilde{\bt}, r,d)/N^{(n)}_r(d)}$
defined in \eqref{2020.10.11.17.50}.
\end{proof}

\section{Hamiltonian description}\label{S Hamil}

In this section, we give a Hamiltonian description of the vector field determined 
by the isomonodromic deformation. 
We fix $\boldsymbol{\nu} \in N^{(n)}_r(d)$.
In Section \ref{SS const nabla0}, 
we take an affine open covering $\{ M\} $ of
$M_{\cC/T}^{\balpha}(\tilde{\bt},r,d)_{\boldsymbol{\nu}}$.
For each $M$, 
we construct an initial connection $\nabla_0$ by the following idea.
First, for the underlying vector bundle
on $\cC \times_T M$ induced by a universal family 
on $\cC \times_T M^{\balpha}_{\cC/T}(\tilde{\bt},r,d)$,
we give an injective morphism (of locally free sheaves) from some fixed vector bundle 
having same rank.
(This construction is not canonical).
Such a vector bundle with an injective morphism is treated in \cite[Section 2]{Hurt}.
Second, we construct a connection on the fixed vector bundle.
By the injective morphism and the connection on the fixed vector bundle, 
we have a connection on the underlying vector bundle 
on $M$ induced by a universal family on $M^{\balpha}_{\cC/T}(\tilde{\bt},r,d)$.
This connection is an initial connection $\nabla_0$.
In Section \ref{SS alge split nabla 0}, we construct vector fields 
on each $M$ associated to an initial connection $\nabla_0$.
Note that an initial connection $\nabla_0$ has poles 
along divisors on $\cC_{M} \setminus D(\tilde{\bt})_{M}$.
Then, for construction of the vector fields,
we need some condition of deformations of $n$-pointed curves 
on neighborhoods of the poles along divisors on $\cC_{M} \setminus D(\tilde{\bt})_{M}$.
We use these vector fields instead of vector fields associated to time variables.
In Section \ref{SS Hamil}, we describe the main theorem.
First, we give a 2-form $\omega$ on $M^{\balpha}_{\cC/T}(\tilde{\bt},r,d)_{\boldsymbol{\nu}}$
such that the kernel $\mathrm{Ker}(\omega)$ induces the vector fields determined 
by the isomonodromic deformations and $\omega$ is the symplectic form fiberwise.
Second, we define Hamiltonian functions on each $M$.
Finally, if we take good coordinates on each $M$,
we obtain a Hamiltonian description of the vector field on each $M$ induced 
by the isomonodromic deformation.

\subsection{Construction of initial connections $\tilde{\nabla}_0^{\sigma_M}$}\label{SS const nabla0}

Let $T$ be a smooth algebraic scheme which is an \'etale covering 
of the moduli stack $\cM_{g,n}$ of $n$-pointed smooth projective curves of genus $g$ 
over $\mathbb{C}$ 
and take a universal family $(\cC,\tilde{t}_1, \ldots,\tilde{t}_n)$ over $T$.
We take a $T$-ample line bundle $\cO_{\cC}(1)$ on $\cC$.
Let $m$ be an integer sufficiently large.
We take an exact sequence 
\begin{equation*}
0 \lra \cO_{\cC} (-m) \lra \cO_{\cC} \lra \cO_{\cC}/ \cO_{\cC} (-m) \lra 0.
\end{equation*}
From this section,
we fix $\boldsymbol{\nu} \in N^{(n)}_r(d)$
and we denote by
\begin{equation*}
\varpi_{\boldsymbol{\nu}}  \colon 
M_{\cC/T}^{\balpha}(\tilde{\bt},r,d)_{\boldsymbol{\nu}} 
\longrightarrow T\times \{ \boldsymbol{\nu} \}
\end{equation*}
the fiber product of 
$M_{\cC/T}^{\balpha}(\tilde{\bt},r,d) \ra T \times N^{(n)}_r(d)$ 
and 
$T\times \{ \boldsymbol{\nu} \}  \ra T \times N^{(n)}_r(d)$.
Let $\tilde{E}$ be the underlying vector bundle of the universal family 
of $M_{\cC/T}^{\balpha}(\tilde{\bt},r,d)_{\boldsymbol{\nu}}$.
Let $M$ be 
an affine open subset of $M^{\balpha}_{\cC/T}(\tilde{\bt},r,d)_{\boldsymbol{\nu}}$.
We put $\cC_M:= \cC \times_TM$.
Let $D(m)$ be the Cartier divisor of $\cC_{M}$ 
such that for any $x \in M$,
\begin{equation*}
\begin{aligned}
D(m)_{\cO_{\cC_{x}}}&= 
\sum_{p \in \cC_{x}} \mathrm{length} \left(\mathrm{Coker}\left(\cO_{\cC_{M}} (-m) 
\ra \cO_{\cC_{M}}\right)_{p}\right)[p] .
\end{aligned}
\end{equation*}
We assume that $D(m)_{\cO_{\cC_x}}$ 
consists of distinct points for any $x \in M$.
Let 
$\pi_{M_{\cC/T}^{\balpha}(\tilde{\bt} ,r,d)_{\boldsymbol{\nu}}} \colon 
\cC \times_T M_{\cC/T}^{\balpha}(\tilde{\bt} ,r,d)_{\boldsymbol{\nu}}
 \rightarrow M_{\cC/T}^{\balpha}(\tilde{\bt} ,r,d)_{\boldsymbol{\nu}}$ 
and 
$\pi_{M} \colon \cC_M \rightarrow M$ 
be the projections, respectively.

\begin{Prop}\label{Prop nabla0Wr}
\textit{There exists an affine open
 covering $\{ M \}$ of $M^{\balpha}_{\cC/T}(\tilde{\bt},r,d)_{\boldsymbol{\nu}}$ 
 such that for each $M$,}
\begin{itemize}
\item[(1)] \textit{$(\pi_{M_{\cC/T}^{\balpha}(\tilde{\bt} ,r,d)_{\boldsymbol{\nu}}})_*
(\tilde{E}(m))|_M$
is a free sheaf on $M$, and}
\item[(2)]\textit{there is a subbundle inclusion 
$\mathcal{O}_M^{\oplus r} \hookrightarrow 
(\pi_{M_{\cC/T}^{\balpha}(\tilde{\bt} ,r,d)_{\boldsymbol{\nu}}})_*(\tilde{E}(m))|_M$
such that for each $x \in M$, the composition 
\begin{equation}\label{2020.5.6.22.21}
\mathcal{O}_M^{\oplus r} \otimes \mathcal{O}_{\cC_M} \otimes k(x)
 \xrightarrow{\ \subset \ } 
(\pi_{M_{\cC/T}^{\balpha}(\tilde{\bt} ,r,d)_{\boldsymbol{\nu}}})_*(\tilde{E}(m))|_{M}
\otimes \mathcal{O}_{\cC_M} \otimes k(x)
\longrightarrow \tilde{E}(m)\otimes k(x)
\end{equation}
is an injection and the cokernel of this injection gives a reduced divisor
whose support is disjoint from 
the supports of $D(\tilde{\bt})_{\cO_{\cC_{x}}}$ and $D(m)_{\cO_{\cC_{x}}}$.}
\end{itemize}
\textit{We denote by $\sigma_M\colon \mathcal{O}_{\cC_M}^{\oplus r} 
\rightarrow \tilde{E}(m)|_{\cC_M}$ 
the composition
$$ \mathcal{O}_{\cC_M}^{\oplus r} \cong
\mathcal{O}_M^{\oplus r} \otimes \mathcal{O}_{\cC_M}
\xrightarrow{\ \subset \ } 
(\pi_{M_{\cC/T}^{\balpha}(\tilde{\bt} ,r,d)_{\boldsymbol{\nu}}})_*(\tilde{E}(m))|_M 
\otimes \mathcal{O}_{\cC_M}
\longrightarrow
\tilde{E}(m)|_{\cC_M}.
$$
}
\end{Prop}

\begin{proof}
Since $m$ is a sufficiently large integer, 
$(\pi_{M_{\cC/T}^{\balpha}(\tilde{\bt} ,r,d)_{\boldsymbol{\nu}}})_*(\tilde{E}(m))$
is a locally free sheaf on $M_{\cC/T}^{\balpha}(\tilde{\bt} ,r,d)_{\boldsymbol{\nu}}$.
So we can take an affine open
 covering $\{ M \}$ of $M^{\balpha}_{\cC/T}(\tilde{\bt},r,d)_{\boldsymbol{\nu}}$
 such that $(\pi_{M_{\cC/T}^{\balpha}(\tilde{\bt} ,r,d)_{\boldsymbol{\nu}}})_*
 (\tilde{E}(m))|_M$
 is a free sheaf for each $M$.
We take a point $x \in M$.
We denote by $E_x$ the restriction $\tilde{E}|_{\cC_x}$ of the vector bundle $\tilde{E}$.
Let $y$ be a point on $\cC_x$.
Set
$$
\begin{aligned}
Y_y :=&\  \left\{ 
(\bar{s}_1,\ldots,\bar{s}_r) \in (E_x|_{2y})^{\oplus r} \ \middle| \
\text{ $\bar{s}_1 \wedge \cdots \wedge \bar{s}_{r}$
does not generate $\bigwedge^r E_x|_{2y}$}
 \right\} \text{ and}\\
 Z_y :=& \ \left\{ 
(\bar{s}_1,\ldots,\bar{s}_r) \in Y_y \ \middle| \
\text{ $\bar{s}_1 \wedge \cdots \wedge \bar{s}_{r}=0$
in $\bigwedge^r E_x|_{2y}$}
 \right\}.
\end{aligned}
$$
Then we have $Z_y \subsetneq Y_y \subsetneq (E_x|_{2y})^{\oplus r}$
and $\dim Z_y \leq \dim (E_x|_{2y})^{\oplus r}-2$.
If we choose an integer $m$ sufficiently large, we have an
exact sequence 
$$
0 \longrightarrow H^0(E_x(-2 y) (m) )
\longrightarrow H^0(E_x(m)) \longrightarrow E|_{2y}
\longrightarrow 0
$$ 
for any $y \in \cC_x$.
So the restriction map
$$
\begin{aligned}
\varphi_y \colon H^0(E_x(m))^{\oplus r} &\longrightarrow (E|_{2y})^{\oplus r} \\
(s_1,\ldots,s_r) &\longmapsto (\bar{s}_1,\ldots,\bar{s}_r)
\end{aligned}
$$
is surjective.
Fibers of $\varphi_y$ are isomorphic to $H^0(E_x (-2y )(m))^{\oplus r}$.
Then we get the inequality 
$\dim \varphi_y^{-1} (Z_y) \leq \dim H^0(E_x(m))^{\oplus r}  -2$.
This inequality implies that $  H^0(E_x(m))^{\oplus r} \setminus
 \bigcup_{y \in \cC_x} \varphi_y^{-1} (Z_y) \neq \emptyset$.
 We take a basis $\{ \boldsymbol{e}_1, \ldots,\boldsymbol{e}_r\}$
of the free sheaf $(\pi_{M_{\cC/T}^{\balpha}(\tilde{\bt} ,r,d)})_*(\tilde{E}(m))|_M$
such that 
$$\{ (\boldsymbol{e}_1)_x, \ldots, (\boldsymbol{e}_r)_x\}
\in H^0(E_x(m))^{\oplus r} \setminus
 \bigcup_{y \in \cC_x} \varphi_y^{-1} (Z_y)$$ for the point $x\in M$.
Then
we may assume that
the subbundle inclusion 
$\mathcal{O}_M^{\oplus r} \hookrightarrow 
(\pi_{M_{\cC/T}^{\balpha}(\tilde{\bt} ,r,d)_{\boldsymbol{\nu}}})_*(\tilde{E}(m))|_M$
 defined by $\{ \boldsymbol{e}_1, \ldots,\boldsymbol{e}_r\}$
satisfies the condition (2) in the statement of this porposition
by taking a refined affine open covering $\{ M \}$.
\end{proof}

For each family $\pi_M\colon \mathcal{C}_M \rightarrow M$ of curves, 
we assume that $\mathrm{Supp}(D(m))$ is disjoint from $\mathrm{Supp}(D(\tilde{\bt}))$,
since $\cO_{\cC}(1)$ is a $T$-ample line bundle.
Let $D(\sigma_M)$ be the Cartier divisor of $\cC_{M}$ 
such that for any $x \in M$,
\begin{equation*}
\begin{aligned}
D(\sigma_M)_{\cO_{\cC_{x}}}
&= \sum_{p\in \cC_{x}} \mathrm{length}\left(\mathrm{Coker}\left(\sigma_M\right)_{p}\right)[p].
\end{aligned}
\end{equation*}
Let $d_m$ be the relative connection induced by $0 \ra \cO_{\cC_{M}}^{\oplus r} (-m) 
\ra \cO_{\cC_{M}}^{\oplus r} $
and the relative exterior derivative $d_{\cC_{M}/M}$ on $\cO_{\cC_{M}}^{\oplus r}$,
that is,
the diagram 
\begin{equation*}
\xymatrix{
\cO_{\cC_{M}}^{\oplus r} (-m) \ar[r]^-{d_m} \ar[d]  
& \cO_{\cC_{M}}^{\oplus r} (-m) \otimes \Omega^1_{\cC_{M}/M }(D(m)) \ar[d] \\
\cO_{\cC_{M}}^{\oplus r}  \ar[r]^-{d_{\cC_{M}/M}} 
& \cO_{\cC_{M}}^{\oplus r}   \otimes \Omega^1_{\cC_{M}/M }(D(m))
}
\end{equation*}
is commutative.

\begin{Def}\label{definition nabla0}
For each affine open subset $M \subset M^{\balpha}_{\cC/T}(\tilde{\bt},r,d)_{\boldsymbol{\nu}}$ 
of Proposition \ref{Prop nabla0Wr},
we fix $\sigma_{M} \colon  \mathcal{O}_{\cC_M}^{\oplus r} 
\rightarrow \tilde{E}(m)|_{\cC_M}$
as in this proposition.
We define a \textit{relative initial connection}
\begin{equation*}
\tilde{\nabla}_{0}^{\sigma_M} \colon  \tilde{E}_{\cC_M}   
\lra \tilde{E}_{\cC_M}   \otimes \Omega^1_{\cC_{M}/M}(D(m)+D(\sigma_M))
\end{equation*}
by the relative connection induced by 
$0 \ra \cO_{\cC_{M}}^{\oplus r} (-m) \xrightarrow{\sigma_M} \tilde{E}|_{\cC_{M}}$ 
and the relative connection $d_m$ on $\cO_{\cC_{M}}^{\oplus r} (-m)$.
\end{Def}

We will construct a parabolic structure of 
$(\tilde{E}_{\cC_M}  ,\tilde{\nabla}_{0}^{\sigma_M})$ as follows.
Let $\{\tilde{\boldsymbol{p}}_1,\ldots,\tilde{\boldsymbol{p}}_{N_1}\}$ be 
the support of the divisor $D(m)$ 
and $\{\tilde{\boldsymbol{p}}'_1,\ldots,\tilde{\boldsymbol{p}}'_{N_2}\}$ be 
the support of the divisor $D(\sigma_M)$. 
Here we put
$N_1 := \deg D(m)|_{\cC_x}$ and $N_2:=\deg D(\sigma_M)|_{\cC_x}$
for each $x \in M$.
Let $ \{ (U_{\alpha}, g_{\alpha}) \}_{\alpha}$ 
and $\{ (U_{\alpha}, h_{\alpha}) \}_{\alpha}$ be the Cartier divisors
corresponding to $D(m)$ and $D(\sigma_M)$, respectively.
Here $\{U_{\alpha}\}_{\alpha}$ is an affine open covering of $\cC_M$.
We take a basis $\boldsymbol{e}_1\otimes g_{\alpha},\ldots,
\boldsymbol{e}_r\otimes g_{\alpha} $ 
of $\cO_{\cC_{M}}^{\oplus r} (-m)|_{U_{\alpha}}$.
We put $\boldsymbol{s}_{\alpha,j}:= 
\sigma_{M}|_{U_{\alpha}}(\boldsymbol{e}_j\otimes g_{\alpha})$ for $j=1,\ldots,r$.
First, we consider a filtration on $\tilde{\boldsymbol{p}}_{i'}$
($i'=1,2,\ldots,N_1$).
Since $\{\tilde{\boldsymbol{p}}_1,\ldots,\tilde{\boldsymbol{p}}_{N_1}\}$ 
and $\{\tilde{\boldsymbol{p}}'_1,\ldots,\tilde{\boldsymbol{p}}'_{N_2}\}$
are disjoint, 
for each $\tilde{\boldsymbol{p}}_{i'}$, 
we can define a filtration $(l^{\sigma_M}_*)^{(i')}$ 
by subbundles $\tilde{E}|_{\tilde{\boldsymbol{p}}_{i'}}=
(l^{\sigma_M}_0)^{(i')} \supset (l^{\sigma_M}_1)^{(i')}
\supset \cdots \supset (l^{\sigma_M}_r)^{(i')}=0$
by using the basis 
$\{ \boldsymbol{s}_{\alpha,1}|_{\tilde{\boldsymbol{p}}_{i'}},
 \ldots, \boldsymbol{s}_{\alpha,r}|_{\tilde{\boldsymbol{p}}_{i'}}\}$
 where $\tilde{\boldsymbol{p}}_{i'} \in  U_{\alpha} $.
By the definition of $\tilde{\nabla}_{0}^{\sigma_M}$,
the residue matrix of $\nabla_0^{\sigma_M}$ along $\tilde{\boldsymbol{p}}_{i'} $
associated to the trivialization given by 
$\{ \boldsymbol{s}_{\alpha,j}\}_{1\le j \le r}$ is 
\begin{equation}\label{2020.5.11.13.13}
\begin{pmatrix}
1 & 0  & \cdots &0 \\
 0 & 1  & \cdots & 0\\
\vdots & \vdots & \ddots &\vdots   \\
0 & 0& \cdots & 1 
\end{pmatrix}.
\end{equation}
Then we have that
 the filtration $(l^{\sigma_M}_*)^{(i')}$ is compatible with $\tilde{\nabla}_{0}^{\sigma_M}$.
 Second, we consider a filtration on $\tilde{\boldsymbol{p}}'_{i''}$
 ($i''=1,2,\ldots,N_2$).
For each $i'' $,
we take $\alpha$ such that $\tilde{\boldsymbol{p}}'_{i''} \in U_{\alpha}$.
We change the order of 
$( \boldsymbol{s}_{\alpha,1} , \boldsymbol{s}_{\alpha,2},\ldots
,\boldsymbol{s}_{\alpha,r})$ as follows.
For each $i''$,
we set 
$ \boldsymbol{s}_{\alpha,i'',k} := \boldsymbol{s}_{\alpha,j_{(i'',k)}}$
(where $k=1,2,\ldots,r$ and $\{ j_{(i'',1)} ,\ldots,j_{(i'',r)}\} =\{ 1,2,\ldots,r\}$)
such that 
$(\boldsymbol{s}_{\alpha,i'',2})|_{\tilde{\boldsymbol{p}}'_{i''}},
\ldots, (\boldsymbol{s}_{\alpha,i'',r})|_{\tilde{\boldsymbol{p}}'_{i''}}$
are linearly independent. 
Let $a_{\alpha,i'',2},\ldots,a_{\alpha,i'',r}$ be elements 
of $\pi_M^*(\mathcal{O}_M)|_{U_{\alpha}}$
such that  
$ \boldsymbol{s}_{\alpha,i'',1}|_{\tilde{\boldsymbol{p}}'_{i''}}
= a_{\alpha,i',2}\cdot \boldsymbol{s}_{\alpha,i'',2}|_{\tilde{\boldsymbol{p}}'_{i''}}
+\cdots+ a_{\alpha,i',r}\cdot \boldsymbol{s}_{\alpha,i'',r}|_{\tilde{\boldsymbol{p}}'_{i''}}$.
We define $\hat{\boldsymbol{s}}_{\alpha,i'',1}$
as 
$$\boldsymbol{s}_{\alpha,i'',1}= h_{\alpha} \cdot \hat{\boldsymbol{s}}_{\alpha,i'',1}
+ a_{\alpha,i'',2}\cdot\boldsymbol{s}_{\alpha,i'',2}
+\cdots+ a_{\alpha,i'',r}\cdot\boldsymbol{s}_{\alpha,i'',r}.$$
For each $\tilde{\boldsymbol{p}}'_{i''}$, 
we can define a filtration $(l'^{\sigma_M}_*)^{(i'')}$ 
by subbundles $\tilde{E}|_{\tilde{\boldsymbol{p}}'_{i''}}=
(l'^{\sigma_M}_0)^{(i'')} \supset (l'^{\sigma_M}_1)^{(i'')}
\supset \cdots \supset (l'^{\sigma_M}_r)^{(i'')}=0$
by using the basis $\{ \hat{\boldsymbol{s}}_{\alpha,i'',1}|_{\tilde{\boldsymbol{p}}'_{i''}},
 \boldsymbol{s}_{\alpha,i'',2}|_{\tilde{\boldsymbol{p}}'_{i''}},
 \ldots, \boldsymbol{s}_{\alpha,i'',r}|_{\tilde{\boldsymbol{p}}'_{i''}}\}.$
Next we check 
the filtration $(l^{\sigma_M}_*)^{(i'')}$ is compatible with $\tilde{\nabla}_{0}^{\sigma_M}$.
We set
\begin{equation}\label{2019.8.3.7.59}
T_{\alpha,i''}:=
\begin{pmatrix}
h_{\alpha} & 0  & 0& \cdots &0 \\
 a_{\alpha,i'',2} & 1 & 0 & \cdots & 0\\
a_{\alpha,i'',3} & 0& 1 & \cdots & 0 \\ 
\cdots & \cdots &\cdots &\cdots &\cdots  \\
a_{\alpha,i'',r} & 0& 0& \cdots & 1 
\end{pmatrix},
\end{equation}
which gives a commutative diagram
\begin{equation}\label{2020.5.6.13.07}
\xymatrix{
 \cO^{\oplus r}_{\cC_M}(-m) |_{U_{\alpha}}
  \ar[d]^-{\sigma_{M}|_{U_{\alpha}}}  
& \cO_{U_{\alpha}}^{\oplus r} \ar[d]^-{T_{\alpha,i''}} \ar[l]_-{g_{\alpha}}\\
\tilde{E}|_{U_{\alpha}} \ar[r]^-{\cong} 
& \cO_{U_{\alpha}}^{\oplus r} \rlap{,}
}
\end{equation}
where the morphism $\tilde{E}|_{U_{\alpha}}
\rightarrow \cO_{U_{\alpha}}^{\oplus r}$ is defined by 
$\hat{\boldsymbol{s}}_{\alpha,i''}:=\{ \hat{\boldsymbol{s}}_{\alpha,i'',1},
 \boldsymbol{s}_{\alpha,i'',2},
 \ldots, \boldsymbol{s}_{\alpha,i'',r}\}.$
Here $U_{\alpha}$ shrinks so that 
$\tilde{\boldsymbol{p}}'_{i''} \in U_{\alpha}$ 
and the morphism from $\tilde{E}|_{U_{\alpha}}$
to $\cO_{U_{\alpha}}^{\oplus r}$ in the diagram \eqref{2020.5.6.13.07} is an isomorphism. 
This isomorphism gives a trivialization of $\tilde{E}|_{U_{\alpha}}$.
The residue matrix of $\nabla_0^{\sigma_M}$ along $\tilde{\boldsymbol{p}}'_{i''} $
associated to this trivialization is 
\begin{equation}\label{2020.5.11.13.14}
\begin{pmatrix}
-\res_{h_{\alpha}=0}(h_{\alpha}^{-1} dh_{\alpha}) & 0  & \cdots &0 \\
-\res_{h_{\alpha}=0}(h_{\alpha}^{-1} d a_{\alpha,i'',2})  & 0  & \cdots & 0\\
\vdots & \vdots & \ddots &\vdots   \\
-\res_{h_{\alpha}=0}(h_{\alpha}^{-1} d a_{\alpha,i'',r}) & 0& \cdots & 0 
\end{pmatrix}.
\end{equation}
Then we have that 
 the filtration $(l'^{\sigma_M}_*)^{(i'')}$ is compatible with $\tilde{\nabla}_{0}^{\sigma_M}$.

\begin{Def}\label{2020.5.6.13.33}
For each affine open subset $M \subset M^{\balpha}_{\cC/T}(\tilde{\bt},r,d)_{\boldsymbol{\nu}}$ 
of Proposition \ref{Prop nabla0Wr},
we fix $\sigma_{M} \colon  \mathcal{O}_{\cC_M}^{\oplus r} 
\rightarrow \tilde{E}(m)|_{\cC_M}$
as in this proposition.
Let $(\tilde{E}_{\cC_M} , \tilde{\nabla}_{0}^{\sigma_M})$
be the relative initial connection defined in Definition \ref{definition nabla0}.
Let $ \boldsymbol{l}^{\sigma_M}:= \{ (l^{\sigma_M}_*)^{(i')}, (l'^{\sigma_M}_*)^{(i'')} 
\}_{1 \leq i' \leq N_1, 
1\leq i'' \leq N_2}$ be the filtrations defined as above.
We call $(\tilde{E}_{\cC_M} , \tilde{\nabla}_{0}^{\sigma_M},
\boldsymbol{l}^{\sigma_M})$
a \textit{relative initial parabolic connection} on $\pi_M \colon \cC_M\rightarrow M$.
\end{Def}

\begin{Def}\label{2020.11.7.17.19}
For $\alpha$ such that $\tilde{\boldsymbol{p}}_{i'} \in U_{\alpha}$
and $\tilde{\boldsymbol{p}}'_{i''} \notin U_{\alpha}$,
the sections
$\{ \boldsymbol{s}_{\alpha,i'',k}\}_{1\le k \le r}$
 give a trivialization of $E|_{U_{\alpha}}$.
For $\alpha$ such that $\tilde{\boldsymbol{p}}'_{i''} \in U_{\alpha}$,
the sections
$\{ \hat{\boldsymbol{s}}_{\alpha,i'',1},
 \boldsymbol{s}_{\alpha,i'',k}\}_{2\le k \le r}$
 give a trivialization of $E|_{U_{\alpha}}$.
We say trivializations $(U_{\alpha},
\bar{\phi}_{\alpha} \colon \tilde{E}|_{U_{\alpha}}
\rightarrow \mathcal{O}^{\oplus r}_{U_{\alpha}})$ 
are \textit{compatible with $\sigma_M$}
if the trivialization $\bar{\phi}_{\alpha}$ 
(for $\alpha$ such that $\tilde{\boldsymbol{p}}_{i'} \in U_{\alpha}$
or $\tilde{\boldsymbol{p}}'_{i''} \in U_{\alpha}$)
coincides with
the trivialization given by the sections 
$\{ \boldsymbol{s}_{\alpha,i'',k}\}_{1\le k \le r}$ or 
the trivialization given by the sections 
$\{ \hat{\boldsymbol{s}}_{\alpha,i'',1},
 \boldsymbol{s}_{\alpha,i'',k}\}_{2\le k \le r}$.
\end{Def}

Now we have two families of parabolic connections parametrized 
by $M\subset M_{\cC/T}^{\balpha}(\tilde{\bt} ,r,d)_{\boldsymbol{\nu}} $: 
$$
(\tilde{E}_{\cC_M} , \tilde{\nabla}_{\cC_M},
\{ l_{*}^{(i)}\}_{1 \leq i \leq n})
\quad \text{and} \quad
(\tilde{E}_{\cC_M} , \tilde{\nabla}_{0}^{\sigma_M},
\boldsymbol{l}^{\sigma_M}).
$$
The former is induced by a universal family
of the moduli space 
$M_{\cC/T}^{\balpha}(\tilde{\bt} ,r,d)_{\boldsymbol{\nu}}$.
The latter is defined in Definition \ref{2020.5.6.13.33}.
By these two families, 
we have two morphisms from $M$ to some moduli spaces.
The first morphism is by taking difference of the connections:
$\tilde{\nabla}_{\cC_M} - \tilde{\nabla}_{0}^{\sigma_M}$.
Note that the difference of connections 
is a \textit{Higgs field} on the vector bundle.
We call a vector bundle with parabolic structures a \textit{parabolic bundle}.
We call a parabolic bundle with a Higgs field 
(which is compatible with the parabolic structures) 
a \textit{parabolic Higgs bundle}.
Then we have a morphism to a moduli space of parabolic Higgs bundles,
roughly speaking.
The second morphism is by taking $\tilde{\nabla}_{0}^{\sigma_M}$.
Then we have a morphism to another moduli space of parabolic connections,
roughly speaking.
We will see a rough sketch of construction of these morphisms.
We set $N:=N_1+ N_2$. 
Let $\mathcal{M}_{g, n + N}$ be the
moduli stack of
smooth projective curves of genus $g$ with (ordered) points of degree $n+N$.
First for each $(C,  \mathrm{Supp}(D(t) + D(p)) ) \in \mathcal{M}_{g, n + N}$,
we consider 
parabolic Higgs bundles of rank $r$ and of degree $d$ on $C$ 
with poles on $D(t)+ D(p)$.
Let $\varpi_H'\colon M_H' \rightarrow \mathcal{M}_{g, n + N}$ be
the moduli space of such parabolic Higgs bundles.
We define a morphism from $M$ to $M'_H$ as follows:
\begin{equation}\label{2020.5.7.21.58}
\begin{aligned}
h_{\nabla - \nabla_0} \colon M &\longrightarrow M'_H \\
x  &\longmapsto ((\cC_x, \mathrm{Supp}(D(\tilde{\bt})_{\cC_x} +D( \tilde{\boldsymbol{p}})_{\cC_x} ) ), 
(\tilde{E}_{\cC_x}, 
(\tilde{\nabla}- \tilde{\nabla}_{0}^{\sigma_M})|_{\cC_x},
\boldsymbol{l}^{\sigma_M}|_{\cC_x} \cup \{ l_{*}^{(i)}|_{\cC_x}\}_{i}
 )).
\end{aligned}
\end{equation}
Second, for each $(C,\mathrm{Supp} (D(t)+D(p))) \in \mathcal{M}_{g, n + N}$,
we consider 
parabolic connections of rank $r$ and of degree $d$ on $C$ 
with poles on $D(t)+ D(p)$.
Let $\varpi'\colon M' \rightarrow \mathcal{M}_{g, n + N}$ be
the moduli space of such parabolic connections.
We define a morphism from $M$ to $M'$ as follows:
\begin{equation}\label{2020.5.7.21.59}
\begin{aligned}
h_{\nabla_0} \colon M &\longrightarrow M' \\
x  &\longmapsto ((\cC_x, \mathrm{Supp}(D(\tilde{\bt})_{\cC_x} +D( \tilde{\boldsymbol{p}})_{\cC_x} ) ), 
(\tilde{E}_{\cC_x}, 
\tilde{\nabla}_{0}^{\sigma_M}|_{\cC_x},
\boldsymbol{l}^{\sigma_M}|_{\cC_x} \cup \{ l_{*}^{(i)}|_{\cC_x}\}_{i} )).
\end{aligned}
\end{equation}

\subsection{Algebraic vector fields 
associated to $\tilde{\nabla}_0^{\sigma_M}$}\label{SS alge split nabla 0}

We take an affine open covering 
$\{M \}$ of $M^{\balpha}_{\cC/T}(\tilde{\bt},r,d)_{\boldsymbol{\nu}}$ 
as in Proposition \ref{Prop nabla0Wr}.
We denote by the same notations $D(\tilde{\bt})$ and $(\tilde{E},\tilde{\nabla},\{ \tilde{l}_j^{(i)}\})$ 
the pull-backs of the divisor $D(\tilde{\bt})$ 
and a universal family $(\tilde{E},\tilde{\nabla},\{ \tilde{l}_j^{(i)}\})$
under the morphism $\cC_M \ra \cC \times_TM^{\balpha}_{\cC/T}(\tilde{\bt},r,d)_{\boldsymbol{\nu}}$, 
respectively.
We fix an injection $\sigma_{M} \colon \mathcal{O}_{\cC_M}^{\oplus r} 
\rightarrow \tilde{E}(m)$
for each $M$ as in Proposition \ref{Prop nabla0Wr}.
We put $D( \tilde{\boldsymbol{p}} ) := D(m) + D(\sigma_M)$.
In this section, we show the following.

\begin{Prop}\label{prop algebra split 2nd}
\textit{
Let $\tilde{\nabla}_0^{\sigma_M}$ be the relative initial connection of Definition} 
\ref{definition nabla0}.
\textit{
For $\mu \in H^1 (\cC_M, \Theta_{C_M/M}(-D(\tilde{\bt})))$,
we take a lift $\hat{\mu} \in H^1(\cC_M, \Theta_{\cC_M/M}( - D(\tilde{\bt}) 
-  D( \tilde{\boldsymbol{p}} )))$ as} (\ref{mulift}) \textit{below}.
\textit{
Then we can construct an algebraic vector field 
associated to $\tilde{\nabla}_0^{\sigma_M}$ on $M$}  
(\textit{Lemma} \ref{Lemma alg vect field nabla0} \textit{below}).
\textit{This vector field is described by
$[ ( \{ 
u_{\alpha\beta}^{\hat{\mu}\nabla_0}\}, 
\{  -v_{\alpha}^{\hat{\mu}\nabla_0}\} ) ] \in \bH^1(\cG^{\bullet}_M)$ 
in Lemma} \ref{Lemma alg vect field nabla0}.
\end{Prop}

Let $\{U_{\alpha} \}_{\alpha}$ be an affine open covering of $\cC_M$
such that $\sharp\{ i \mid \hat{t}_i |_{\cC_M} \cap U_{\alpha} \neq \emptyset \} \le 1$ for any $\alpha$ and 
$\sharp\{ \alpha \mid \hat{t}_i |_{\cC_M} \cap U_{\alpha} \neq \emptyset \} \le 1$ for any $i$.
Here $\{ \hat{t}_i \}$ is the set of the supports 
of the Cartier divisor $D(\tilde{\bt})+D( \tilde{\boldsymbol{p}} )$.
Set
$$I_{D( \tilde{\boldsymbol{p}} )} := \{ \alpha \mid 
U_{\alpha} \cap \mathrm{Supp}(D( \tilde{\boldsymbol{p}} )) \neq \emptyset\}.$$
We denote by $D(\tilde{\bt})=\{(U_{\alpha}, f_{\alpha}) \}_{\alpha}$, 
$D(m)= \{ (U_{\alpha}, g_{\alpha}) \}_{\alpha}$ 
and $D(\sigma_M)= \{ (U_{\alpha}, h_{\alpha}) \}_{\alpha}$ the Cartier divisors.
For any $\alpha$, we assume that there exists a trivialization 
$\bar{\phi}_{\alpha} \colon  \tilde{E}|_{U_{\alpha}} \xrightarrow{\sim} \cO^{\oplus r}_{U_{\alpha}}$
of $\tilde{E}$ over $U_{\alpha}$.

We take $\mu  \in H^1 (\cC_M, \Theta_{\cC_M/M}(-D(\tilde{\bt})))$.
Since
\begin{equation*}
H^1(\cC_M, \Theta_{\cC_M/M}( - D(\tilde{\bt}) - D( \tilde{\boldsymbol{p}} ))) 
\lra H^1 (\cC_M, \Theta_{\cC_M/M}(-D(\tilde{\bt}))) \lra 0,
\end{equation*}
we choose a lift 
\begin{equation}\label{mulift}
\hat{\mu} = \left[ \{ \hat{d}_{\alpha\beta} \} \right] 
\in H^1(\cC_M, \Theta_{\cC_M/M}( - D(\tilde{\bt})-  D( \tilde{\boldsymbol{p}} )))
\end{equation}
of $\mu$ as follows.
First we take the isomonodromic lift $M[\epsilon] 
\rightarrow M$ associated to $\mu$ as in Section \ref{SS IMD by c c}.
Here we set $M[\epsilon] := M \times \mathrm{Spec}\, \mathbb{C} [\epsilon]$.
We denote by $\mathrm{im}_{M}^{\mu}\colon M[\epsilon] \rightarrow M$ 
this isomonodromic lift.
Let $\cC_{\epsilon}$ be the fiber product $\cC_M \times_M M[\epsilon]$
with respect to the projection $\cC_M \rightarrow M$ and 
the morphism $\mathrm{im}_{M}^{\mu}\colon M[\epsilon] \rightarrow M$.
The lift $\mathrm{im}_{M}^{\mu}$ induces a morphism
\begin{equation}\label{2020.5.3.15.52}
\mathrm{IM}_{M}^{\mu} \colon \cC_{\epsilon} \longrightarrow \cC_M.
\end{equation}
Second we take the pull-back 
$(\mathrm{IM}_{M}^{\mu})^*(D(\tilde{\bt})+D( \tilde{\boldsymbol{p}} ))$,
which is a Cartier divisor on $\cC_{\epsilon}$.
We consider the pair 
$(\cC_{\epsilon},\mathrm{Supp}((\mathrm{IM}_{M}^{\mu})^*(D(\tilde{\bt})+D( \tilde{\boldsymbol{p}} ))) )$,
which is a flat family of curves with (ordered) points.
This flat family gives a morphism $M[\epsilon] \rightarrow \mathcal{M}_{g, n+ N}$
to the moduli stack $\mathcal{M}_{g, n+ N }$ of 
smooth projective curves of genus $g$ with points of degree $n+N$.
The morphism $M[\epsilon] \rightarrow \mathcal{M}_{g, n+ N}$ 
determines a class of 
$H^1(\cC_M, \Theta_{\cC_M/M}( - D(\tilde{\bt})-  D( \tilde{\boldsymbol{p}} )))$.
We denote by $\hat{\mu}$ this class, which is a lift of $\mu$.

Now we define an algebraic vector field on $M$ 
by the relative initial connection $(\tilde{E},\tilde{\nabla}_0^{\sigma_M})$ 
and the lift $\hat{\mu}=[ \{ \hat{d}_{\alpha\beta} \} ] $ of $\mu$.
Let $\cG_M^0$ and $\cG_M^1$ be 
the pull-backs of $\cG^0$ and $\cG^1$
by the morphism 
$\cC_M \rightarrow \cC \times_T M^{\balpha}_{\cC/T}(\tilde{\bt},r,d)_{\boldsymbol{\nu}}$,
respectively.
We have the morphism $\nabla_{\cG_M^{\bullet}} \colon \cG_M^0 \rightarrow \cG_M^1$
induced by $\nabla_{\cG^{\bullet}} \colon \cG^0 \rightarrow \cG^1$.
Let $\{ (U_{\alpha}, z_{\alpha}) \}_{\alpha}$ be the Cartier divisor 
$D(\tilde{\bt}) + D( \tilde{\boldsymbol{p}} )$.
Let $\tilde{A}^0_{\alpha}  z_{\alpha}^{-1} dz_{\alpha}$ 
and $\tilde{A}_{\alpha}  f_{\alpha}^{-1} df_{\alpha}$ 
be connection matrices of $\tilde{\nabla}_0^{\sigma_M}$ 
and $\tilde{\nabla}$ on $U_{\alpha}$ associated to the trivialization $\bar{\phi}_{\alpha}$, 
respectively.
First, we set
\begin{equation}\label{2019.7.31.11.05}
u_{\alpha\beta}^{\hat{\mu}\nabla_0} := 
\bar{\phi}_{\beta}^{-1}\circ \left( \hat{d}_{\alpha\beta} + 
\left< \hat{d}_{\alpha\beta} ,  \tilde{A}^0_{\beta} \frac{d z_{\beta}}{z_{\beta}} \right> \right)
\circ \bar{\phi}_{\beta} \in \cG_M^0 (U_{\alpha\beta}),
\end{equation}
which satisfy the equality $u_{\beta\gamma}^{\hat{\mu}\nabla_0}
- u_{\alpha\gamma}^{\hat{\mu}\nabla_0} + u_{\alpha\beta}^{\hat{\mu}\nabla_0}=0$.
The idea of the definition of $u_{\alpha\beta}^{\hat{\mu}\nabla_0}$ is 
``isomonodromic deformation'' of $\tilde{\nabla}_0^{\sigma_M}$ associated to $\hat{\mu}$.

Now we construct a vector field on $M$ such that 
its image of the natural morphism
$\bH^1(\cG_M^{\bullet}) \rightarrow H^1(\cG_M^0) $
is $[\{ u_{\alpha\beta}^{\hat{\mu}\nabla_0} \}]$.
We will define $v_{\alpha}^{\hat{\mu}\nabla_0}$ later at \eqref{2019.7.31.11.14}
such that $[(\{ u_{\alpha\beta}^{\hat{\mu}\nabla_0} \}, \{ -v_{\alpha}^{\hat{\mu}\nabla_0} \})] 
\in \bH^1(\cG_M^{\bullet})$.
For $\{ \hat{d}_{\alpha\beta} \}_{\alpha\beta}$,
we put
\begin{equation}\label{2019.7.31.11.02}
\hat{u}_{\alpha\beta}^{\text{IMD}} := \bar{\phi}_{\beta}^{-1}\circ \left( \hat{d}_{\alpha\beta} + 
\left< \hat{d}_{\alpha\beta} ,  \tilde{A}_{\beta} \frac{d f_{\beta}}{f_{\beta}} \right> \right)
\circ \bar{\phi}_{\beta}  , \quad \hat{v}_{\alpha}^{\text{IMD}} :=0.
\end{equation}
By Proposition \ref{Prop explicit description IMD},
the class $[(\{\hat{u}_{\alpha\beta}^{\text{IMD}} \} ,\{ \hat{v}_{\alpha}^{\text{IMD}}\})]
\in \bH^1(\mathcal{G}_M^{\bullet})$ means 
the vector field determined by the isomonodromic deformations of $\tilde{\nabla}$
associated to $\mu$.
Hence $[(\{\hat{u}_{\alpha\beta}^{\text{IMD}} \} ,\{ \hat{v}_{\alpha}^{\text{IMD}}\})]$
gives the isomonodromic lift 
$\mathrm{im}_M^{\mu} \colon 
M[\epsilon] \rightarrow M$.
Let
$(E_{\epsilon},(\mathrm{IM}_{M}^{\mu})^*\tilde{\nabla}, 
\{ (l_{\epsilon})^{(i)}_j \}) $ 
be the pull-back of $(\tilde{E},\tilde{\nabla}, \{ \tilde{l}^{(i)}_j \})$
under the morphism $\mathrm{IM}_{M}^{\mu} \colon \cC_\epsilon \rightarrow \cC_M$,
which is induced by $\mathrm{im}_M^{\mu}$.
Let $D_{\epsilon}(\tilde{\bt})=\{ (U_{\alpha}^{\epsilon}, f_{\alpha}^{\epsilon}) \}_{\alpha}$ 
and
$D_{\epsilon}( \tilde{\boldsymbol{p}} )=
\{ (U_{\alpha}^{\epsilon}, z_{\alpha}^{\epsilon}) \}_{\alpha}$ be 
the pull-backs of the Cartier divisors 
$D(\tilde{\bt})=\{ (U_{\alpha}, f_{\alpha}) \}_{\alpha}$  
and $D( \tilde{\boldsymbol{p}} )=\{ (U_{\alpha}, z_{\alpha}) \}_{\alpha}$ 
by the morphism $\mathrm{IM}_{M}^{\mu} \colon \cC_\epsilon \ra \cC_M$, respectively.

\begin{Def}\label{2019.8.8.14.02}
Let $(\tilde{E},\tilde{\nabla}_0^{\sigma_M})$ be the relative initial connection 
defined in Definition \ref{definition nabla0}.
We define an infinitesimal deformation 
\begin{equation*}
\nabla_{0,\epsilon}^{\sigma_M,\text{IMD}}
 \colon E_{\epsilon} \lra E_{\epsilon} \otimes 
\Omega^1_{\cC_{\epsilon} /M[\epsilon]} 
\left(D_{\epsilon}( \tilde{\boldsymbol{p}} )\right)
\end{equation*}
of the relative initial connection $(\tilde{E},\tilde{\nabla}_0^{\sigma_M})$
by taking the pull-back of $(\tilde{E},\tilde{\nabla}_0^{\sigma_M})$ 
by $\mathrm{IM}_{M}^{\mu} \colon \cC_\epsilon \rightarrow \cC_M$.
\end{Def}

Let $\cC_{\epsilon} = \bigcup_{\alpha} U^{\epsilon}_{\alpha}$ be the open covering 
corresponding to the affine open covering $\{U_{\alpha} \}_{\alpha}$ of $\cC_M$.
Here,
the affine open covering $\{U_{\alpha} \}_{\alpha}$
 is defined after Proposition \ref{prop algebra split 2nd}.
We can take trivializations
$\{ \phi_{\alpha}\}_{\alpha}$ 
of $E_{\epsilon}$ 
on $U_{\alpha}^{\epsilon}$
such that
$\bar{\phi}_{\alpha}= \phi_{\alpha}$ ($\mathrm{mod}$ $\epsilon$)
and 
\begin{equation}\label{2020.5.14.23.06}
\varphi_{\alpha}^{-1}\circ\varphi_{\beta}-\mathrm{id}
= \hat{u}_{\alpha\beta}^{\mathrm{IMD}},
\end{equation}
where $\varphi_{\alpha}$ is defined as in the proof of 
Proposition \ref{Prop d t s of m by c c}.
We set
\begin{equation}\label{2019.7.31.11.14}
 v_{\alpha}^{\hat{\mu}\nabla_0} := (\varphi_{\alpha}^{-1} \otimes \mathrm{id}) 
 \circ (  \nabla_{0,\epsilon}^{\sigma_M,\text{IMD}} ) \circ \varphi_{\alpha}
-  \tilde{\nabla}_0^{\sigma_M}
\end{equation}
for any $\alpha$.

\begin{Lem}\label{Lemma alg vect field nabla0}
\textit{
For $u_{\alpha\beta}^{\hat{\mu}\nabla_0}$ defined as \eqref{2019.7.31.11.02}
and  $v_{\alpha}^{\hat{\mu}\nabla_0}$ defined as \eqref{2019.7.31.11.14},
we have the following equality
\begin{equation}\label{2019.5.28.12.16}
\tilde{\nabla} \circ u_{\alpha\beta}^{\hat{\mu}\nabla_0}  -u_{\alpha\beta}^{\hat{\mu}\nabla_0}  
\circ \tilde{\nabla}
= - (v_{\beta}^{\hat{\mu}\nabla_0}- v_{\alpha}^{\hat{\mu}\nabla_0}).
\end{equation}
Moreover $ v_{\alpha}^{\hat{\mu}\nabla_0}$  is an element of $\epsilon \otimes \cG_M^1$.
As a result, the pair
$\left[ \left( \left\{ 
u_{\alpha\beta}^{\hat{\mu}\nabla_0}
 \right\}, \left\{  
-v_{\alpha}^{\hat{\mu}\nabla_0}
\right\} \right) \right]$ is a class of $\bH^1(\cG^{\bullet}_M)$.
This class is independent of the choice of a representative of 
the lift $\hat{\mu}=[ \{ \hat{d}_{\alpha\beta} \} ] 
\in H^1(\cC_M, \Theta_{\cC_M/M}( - D(\tilde{\bt})-  D( \tilde{\boldsymbol{p}} )))$
of $\mu$.
}
\end{Lem}

\begin{proof}

First, we show the equality (\ref{2019.5.28.12.16}).
We can check the equalities
\begin{equation*}
\begin{aligned}
&(\tilde{\nabla} \circ u_{\alpha\beta}^{\hat{\mu}\nabla_0}  -u_{\alpha\beta}^{\hat{\mu}\nabla_0}  
\circ \tilde{\nabla})(a)\\
&=(\bar{\phi}_{\beta}^{-1}\otimes \mathrm{id})
 \circ\left(d+  \tilde{A}_{\beta}\frac{d f_{\beta}}{f_{\beta}}\right)\circ 
\left( \hat{d}_{\alpha\beta} + 
\left< \hat{d}_{\alpha\beta} ,  \tilde{A}^0_{\beta} \frac{d z_{\beta}}{z_{\beta}} \right> \right)
( \bar{\phi}_{\beta} (a)) \\
&\qquad -(\bar{\phi}_{\beta}^{-1}\otimes \mathrm{id})  \circ
\left( \hat{d}_{\alpha\beta} + 
\left< \hat{d}_{\alpha\beta} ,  \tilde{A}^0_{\beta} \frac{d z_{\beta}}{z_{\beta}} \right> \right) 
\circ \left(d+  \tilde{A}_{\beta}\frac{d f_{\beta}}{f_{\beta}}\right)
( \bar{\phi}_{\beta} (a))\\
&= (\bar{\phi}_{\beta}^{-1}\otimes \mathrm{id}) \circ \left(
d\left(\left< \hat{d}_{\alpha\beta} ,  
\tilde{A}^0_{\beta} \frac{d z_{\beta}}{z_{\beta}} -\tilde{A}_{\beta} \frac{d f_{\beta}}{f_{\beta}} \right>\right) 
( \bar{\phi}_{\beta} (a)) 
 +\left[ \tilde{A}_{\beta}\frac{d f_{\beta}}{f_{\beta}} ,
\left< \hat{d}_{\alpha\beta} ,  \tilde{A}^0_{\beta} \frac{d z_{\beta}}{z_{\beta}} \right>\right]
 ( \bar{\phi}_{\beta} (a)) \right)
\end{aligned}
\end{equation*}
and
\begin{equation*}
\begin{aligned}
&(\tilde{\nabla}_0^{\sigma_M} \circ \hat{u}_{\alpha\beta}^{\text{IMD}}
 -\hat{u}_{\alpha\beta}^{\text{IMD}} \circ \tilde{\nabla}_0^{\sigma_M})(a)\\
 &=(\bar{\phi}_{\beta}^{-1}\otimes \mathrm{id})
 \circ\left(d+ \tilde{A}^0_{\beta} \frac{d z_{\beta}}{z_{\beta}} \right)\circ 
\left( \hat{d}_{\alpha\beta} + 
\left< \hat{d}_{\alpha\beta} , \tilde{A}_{\beta}\frac{d f_{\beta}}{f_{\beta}} \right> \right)
( \bar{\phi}_{\beta} (a)) \\
&\qquad -(\bar{\phi}_{\beta}^{-1}\otimes \mathrm{id})  \circ
\left( \hat{d}_{\alpha\beta} + 
\left< \hat{d}_{\alpha\beta} ,  \tilde{A}_{\beta}\frac{d f_{\beta}}{f_{\beta}} \right> \right) 
\circ \left(d+  \tilde{A}^0_{\beta} \frac{d z_{\beta}}{z_{\beta}}\right)
( \bar{\phi}_{\beta} (a))\\
&= (\bar{\phi}_{\beta}^{-1}\otimes \mathrm{id}) \circ \left(
d\left(\left< \hat{d}_{\alpha\beta} ,  
\tilde{A}_{\beta} \frac{d f_{\beta}}{f_{\beta}}
- \tilde{A}^0_{\beta} \frac{d z_{\beta}}{z_{\beta}}\right>\right) 
( \bar{\phi}_{\beta} (a)) 
 +\left[  \tilde{A}^0_{\beta} \frac{d z_{\beta}}{z_{\beta}} ,
\left< \hat{d}_{\alpha\beta} , \tilde{A}_{\beta}\frac{d f_{\beta}}{f_{\beta}} \right>\right]
 ( \bar{\phi}_{\beta} (a)) \right)
\end{aligned}
\end{equation*}
for $a \in \tilde{E}|_{U_{\alpha\beta}}$.
Since $\cC_M \ra M$ is a family of smooth projective curves, 
we have
$\langle \hat{d}_{\alpha\beta},  d f_{\beta} \rangle dz_{\beta}
=\langle \hat{d}_{\alpha\beta},  d z_{\beta} \rangle df_{\beta}$.
By the equalities above, we have
\begin{equation*}
\tilde{\nabla} \circ u_{\alpha\beta}^{\hat{\mu}\nabla_0}  -u_{\alpha\beta}^{\hat{\mu}\nabla_0}  
\circ \tilde{\nabla}=-(
\tilde{\nabla}_0^{\sigma_M} \circ \hat{u}_{\alpha\beta}^{\text{IMD}}
 -\hat{u}_{\alpha\beta}^{\text{IMD}} \circ \tilde{\nabla}_0^{\sigma_M}).
\end{equation*}
On the other hand, since $v_{\alpha}^{\hat{\mu}\nabla_0}$ is 
defined by the infinitesimal
 deformation of the relative initial connection $\tilde{\nabla}_0^{\sigma_M}$
associated to $(\{\hat{u}_{\alpha\beta}^{\text{IMD}} \},\{ \hat{v}_{\alpha}^{\text{IMD}}\})$,
we can check the equality
\begin{equation}\label{2019.8.6.18.40}
\tilde{\nabla}_0^{\sigma_M} \circ \hat{u}_{\alpha\beta}^{\text{IMD}}
 -\hat{u}_{\alpha\beta}^{\text{IMD}} \circ \tilde{\nabla}_0^{\sigma_M}
= v_{\beta}^{\hat{\mu}\nabla_0}- v_{\alpha}^{\hat{\mu}\nabla_0}.
\end{equation}
Then we obtain the equality
\begin{equation*}
\tilde{\nabla} \circ u_{\alpha\beta}^{\hat{\mu}\nabla_0}  -u_{\alpha\beta}^{\hat{\mu}\nabla_0}  
\circ \tilde{\nabla}
= - (v_{\beta}^{\hat{\mu}\nabla_0}- v_{\alpha}^{\hat{\mu}\nabla_0}).
\end{equation*}

Next, we show that $ v_{\alpha}^{\hat{\mu}\nabla_0}\in 
 \epsilon \otimes \cE nd (\tilde{E}) \otimes \Omega^1_{\cC/T} (D(\tilde{\bt}))$.
For this purpose, 
we show that $\{ v_{\alpha}^{\hat{\mu}\nabla_0} \}$ has no pole 
on the supports of $D(m)$ and $D(\sigma_M)$.
Let $\sigma_{\alpha} \colon U_{\alpha}^{\epsilon} 
\rightarrow U_{\alpha} \times \mathrm{Spec}\, \mathbb{C} [\epsilon]$
be an isomorphism such that 
the isomorphism $\sigma_{\alpha}$ preserves the divisor 
$D_{\epsilon}(\tilde{\boldsymbol{p}})\cap U_{\alpha}^{\epsilon}$ and 
$(D(\tilde{\boldsymbol{p}})\cap U_{\alpha})\times \mathrm{Spec}\, \mathbb{C}[\epsilon]$
and the isomorphisms $\{\sigma_{\alpha}\}_{\alpha}$
give a representative of $\hat{\mu}$ 
by $\hat{d}_{\alpha\beta}=(\sigma_{\alpha}^{*})^{-1}\sigma_{\beta}^* -\mathrm{id}
\colon \cO_{U_{\alpha\beta}} \ra \epsilon \otimes \cO_{U_{\alpha\beta}}$.
We set $U_{I_{D(\tilde{\boldsymbol{p}})}}:= 
\bigcup_{\alpha \in I_{D(\tilde{\boldsymbol{p}})}} U_{\alpha}$.
We consider a subset 
$(\mathrm{IM}_M^{\mu})^{-1}(U_{I_{D(\tilde{\boldsymbol{p}})}} )
\subset \cC_{\epsilon}$ as the tuple of open subsets 
$\{ U_{\alpha} \times \mathrm{Spec}\, \mathbb{C} [\epsilon] 
\}_{\alpha \in I_{D(\tilde{\boldsymbol{p}})}}$
 and gluing data $ \{(\sigma_{\alpha}^{*})^{-1}\sigma_{\beta}^* 
 =\mathrm{id}+ \hat{d}_{\alpha\beta}
 \}_{\alpha,\beta \in I_{D(\tilde{\boldsymbol{p}})}}$:
$$
(\mathrm{IM}_M^{\mu})^{-1}(U_{I_{D(\tilde{\boldsymbol{p}})}} )=
(\{ U_{\alpha} \times \mathrm{Spec}\, \mathbb{C} [\epsilon] 
\}_{\alpha \in I_{D(\tilde{\boldsymbol{p}})}} ,
 \{ \mathrm{id}+ \hat{d}_{\alpha\beta}
 \}_{\alpha,\beta \in I_{D(\tilde{\boldsymbol{p}})}}) .
 $$
By the definition of the lift $\hat{\mu}$, we have an isomorphism 
$$
(\mathrm{IM}_M^{\mu})^{-1}(U_{I_{D(\tilde{\boldsymbol{p}})}} )
 \cong (\{ U_{\alpha} \times \mathrm{Spec}\, \mathbb{C} [\epsilon] 
\}_{\alpha \in I_{D(\tilde{\boldsymbol{p}})}} ,
 \{\mathrm{id} \}_{\alpha,\beta \in I_{D(\tilde{\boldsymbol{p}})}}) 
$$
such that this isomorphism preserves the divisor 
$D_{\epsilon}(\tilde{\boldsymbol{p}})\cap 
(\mathrm{IM}_M^{\mu})^{-1}(U_{I_{D(\tilde{\boldsymbol{p}})}} )$ and 
$(D(\tilde{\boldsymbol{p}})\cap U_{I_{D(\tilde{\boldsymbol{p}})}} )
\times \mathrm{Spec}\, \mathbb{C}[\epsilon]$.
That is, there exists a set $\{ \hat{d}_{\alpha} \}_{\alpha \in I_{D(\tilde{\boldsymbol{p}})}}$
such that $\hat{d}_{\alpha} \in \Theta_{\cC_M/M}(-D(\tilde{\bt})- D( \tilde{\boldsymbol{p}} ))(U_{\alpha})$
and $\hat{d}_{\alpha\beta} = \hat{d}_{\beta} - \hat{d}_{\alpha}$
for $\alpha,\beta \in I_{D(\tilde{\boldsymbol{p}})}$.
We consider a vector bundle 
$E_{\epsilon}|_{(\mathrm{IM}_M^{\mu})^{-1}(U_{I_{D(\tilde{\boldsymbol{p}})}} )}$ 
as a tuple of an open covering $\{ U_{\alpha} \times \mathrm{Spec}\, \mathbb{C} [\epsilon] 
\}_{\alpha \in I_{D(\tilde{\boldsymbol{p}})}}$
 and gluing data 
  $\{ g_{\alpha\beta} 
\}_{\alpha,\beta \in I_{D(\tilde{\boldsymbol{p}})}}$, where 
 \begin{equation*}
g_{\alpha\beta} := 
 (\bar{\phi}_{\alpha} \otimes 1 ) 
\circ (\mathrm{id}+\hat{u}_{\alpha\beta}^{\text{IMD}})
\circ  (\bar{\phi}_{\beta} \otimes 1 )^{-1}.
\end{equation*}
That is,
$E_{\epsilon}|_{(\mathrm{IM}_M^{\mu})^{-1}(U_{I_{D(\tilde{\boldsymbol{p}})}} )}
=(\{ U_{\alpha} \times \mathrm{Spec}\, \mathbb{C} [\epsilon] 
\}_{\alpha \in I_{D(\tilde{\boldsymbol{p}})}},
\{ g_{\alpha\beta} 
\}_{\alpha,\beta \in I_{D(\tilde{\boldsymbol{p}})}})$.
Set 
$$
u_{\alpha}^{\hat{d}_{\alpha} \nabla} := 
\bar{\phi}_{\alpha}^{-1}\circ \left( \hat{d}_{\alpha} + 
\left< \hat{d}_{\alpha} ,  \tilde{A}_{\alpha} \frac{d f_{\alpha}}{f_{\alpha}} \right> \right)
\circ \bar{\phi}_{\alpha}
$$
for $\alpha \in I_{D(\tilde{\boldsymbol{p}})}$.
Since 
$\hat{u}_{\alpha\beta}^{\text{IMD}} =  
u_{\beta}^{\hat{d}_{\beta} \nabla}
-u_{\alpha}^{\hat{d}_{\alpha} \nabla}$,
we have an isomorphism
$$
\phi_{I_{D(\tilde{\boldsymbol{p}})}} \colon  E_{\epsilon}|_{(\mathrm{IM}_M^{\mu})^{-1}
(U_{I_{D(\tilde{\boldsymbol{p}})}} )}
\xrightarrow{\ \cong \ } (\{ U_{\alpha} \times \mathrm{Spec}\, \mathbb{C} [\epsilon] 
\}_{\alpha \in I_{D(\tilde{\boldsymbol{p}})}},
\{  (\bar{\phi}_{\alpha} \otimes 1 ) 
\circ  (\bar{\phi}_{\beta} \otimes 1 )^{-1}
\}_{\alpha,\beta \in I_{D(\tilde{\boldsymbol{p}})}}).
$$
That is, $\phi_{I_{D(\tilde{\boldsymbol{p}})}}$
is an isomorphism from $E_{\epsilon}|_{(\mathrm{IM}_M^{\mu})^{-1}
(U_{I_{D(\tilde{\boldsymbol{p}})}} )}$
to the trivial deformation of 
$\tilde{E}|_{U_{I_{D(\tilde{\boldsymbol{p}})}}}$.
Now we consider a connection 
\begin{equation}\label{2020.10.21.16.47}
(\phi_{I_{D(\tilde{\boldsymbol{p}})}}^{-1})^*
(\nabla_{0,\epsilon}^{\sigma_M,\text{IMD}}|_{(\mathrm{IM}_M^{\mu})^{-1}
(U_{I_{D(\tilde{\boldsymbol{p}})}} )}).
\end{equation}
Let $\tilde{A}_{\alpha}^0 \frac{dz_{\alpha}}{z_{\alpha}}$
be the connection matrix of $\tilde{\nabla}_0^{\sigma_M}$ 
associated to the trivialization $\bar{\phi}_{\alpha}$.
The connection $\nabla_{0,\epsilon}^{\sigma_M,\text{IMD}}$
is defined by the injection $\sigma_M
\colon \cO_{\cC_M} (-m) \hookrightarrow \tilde{E}_{\cC_M}$ and 
the injection $\cO_{\cC_M} (-m) \hookrightarrow \cO_{\cC_M}$,
which are defined in Section \ref{SS const nabla0}.
We restrict the injections $\sigma_M
\colon \cO^{\oplus r}_{\cC_M} (-m) \hookrightarrow \tilde{E}_{\cC_M}$
and $\cO_{\cC_M} (-m) \hookrightarrow \cO_{\cC_M}$
 to the open set $U_{I_{D(\tilde{\boldsymbol{p}})}} $.
The morphism $\mathrm{IM}_{M}^{\mu} \colon \cC_{\epsilon} \rightarrow \cC_M$ 
and the restricted injections
induce injections over $(\mathrm{IM}_M^{\mu})^{-1}
(U_{I_{D(\tilde{\boldsymbol{p}})}} )$.
The images of these induced injections 
under the isomorphism $\phi_{I_{D(\tilde{\boldsymbol{p}})}}$
are independent of $\epsilon$.
Then we have that 
the connection matrix of the connection \eqref{2020.10.21.16.47}
on $U_{\alpha} \times \mathrm{Spec}\, \mathbb{C} [\epsilon]$ 
by the trivialization $(\phi_{I_{D(\tilde{\boldsymbol{p}})}}^{-1})^* (\phi_{\alpha})$ is just 
$\tilde{A}_{\alpha}^0 \frac{dz_{\alpha}}{z_{\alpha}}$.
That is, the $\epsilon$-term of this connection matrix of the connection \eqref{2020.10.21.16.47}
vanishes.
Since this $\epsilon$-term vanishes,  
we have  
$$
v_{\alpha}^{\hat{\mu}\nabla_0} = 
\tilde{\nabla}_0^{\sigma_M} \circ u_{\alpha}^{\hat{d}_{\alpha} \nabla}
 -u_{\alpha}^{\hat{d}_{\alpha} \nabla} \circ \tilde{\nabla}_0^{\sigma_M}.
$$
We may check this equality 
as in \eqref{2020.10.22.13.36} and \eqref{2020.10.22.13.37}.
Since $\hat{d}_{\alpha}$ vanishes at the support of $D(\tilde{\boldsymbol{p}})$,
$u_{\alpha}^{\hat{d}_{\alpha} \nabla}$
also vanishes at the support of $D(\tilde{\boldsymbol{p}})$.
This fact means that $\{ v_{\alpha}^{\hat{\mu}\nabla_0} \}$ has no pole 
on the supports of $D(m)$ and $D(\sigma_M)$.
Since $v_{\alpha}^{\hat{\mu}\nabla_0}$ has no poles along $D( \tilde{\boldsymbol{t}} ) $,
we have the compatibility of $v_{\alpha}^{\hat{\mu}\nabla_0}$
with the parabolic structures $\{ \tilde{l}^{(i)}_j \}$.
Then we have $ v_{\alpha}^{\hat{\mu}\nabla_0} \in \epsilon \otimes \cG_M^1$.
So $[(\{ u_{\alpha\beta}^{\hat{\mu}\nabla_0}\}, \{  -   v_{\alpha}^{\hat{\mu}\nabla_0}\})]$ 
is an element of $\bH^1(\cG^{\bullet}_M)$.

Let $d_{\alpha}$ be an element of 
$\Theta_{\cC_M/M}(-D(\tilde{\bt})- D( \tilde{\boldsymbol{p}} ))(U_{\alpha})$.
Set $\hat{d}'_{\alpha\beta}:= \hat{d}_{\alpha\beta} +d_{\beta}-d_{\alpha}$ and
\begin{equation*}
u_{\alpha}^{d_{\alpha}\nabla_0} := \bar{\phi}_{\alpha}^{-1} \circ \left( d_{\alpha} 
+ \left< d_{\alpha} ,  \tilde{A}^0_{\alpha} \frac{d z_{\alpha}}{z_{\alpha}} \right>
\right) \circ \bar{\phi}_{\alpha}, \text{ and }
u_{\alpha}^{d_{\alpha}\nabla} := \bar{\phi}_{\alpha}^{-1} \circ \left( d_{\alpha} 
+ \left< d_{\alpha} ,  \tilde{A}_{\alpha} \frac{d f_{\alpha}}{f_{\alpha}} \right>
\right) \circ \bar{\phi}_{\alpha}.
\end{equation*}
For the representative
$\{\hat{d}'_{\alpha\beta}\}$,
we define 
${}'u_{\alpha\beta}^{\hat{\mu}\nabla_0}$ and
$(\{ {}'\hat{u}_{\alpha\beta}^{\text{IMD}} \}, \{ {}'\hat{v}_{\alpha}^{\text{IMD}}  \})$ 
as (\ref{2019.7.31.11.05}) and (\ref{2019.7.31.11.02}), respectively.
We can define the infinitesimal deformation 
of the relative initial connection $\tilde{\nabla}_0^{\sigma_M}$
along the vector field
 $(\{ {}'\hat{u}_{\alpha\beta}^{\text{IMD}} \}, \{ {}'\hat{v}_{\alpha}^{\text{IMD}}  \})$ 
as Definition \ref{2019.8.8.14.02}.
By this infinitesimal deformation,
we can define
${}' v_{\alpha}^{\hat{\mu}\nabla_0}$
as (\ref{2019.7.31.11.14}).
We can check the following equalities:
\begin{equation*}
\begin{aligned}
{}'u_{\alpha\beta}^{\hat{\mu}\nabla_0}&=  u_{\alpha\beta}^{\hat{\mu}\nabla_0} 
+\left(u_{\beta}^{d_{\beta}\nabla_0}-u_{\alpha}^{d_{\alpha}\nabla_0} \right) \text{ and} \\
{}' v_{\alpha}^{\hat{\mu}\nabla_0} &=
v_{\alpha}^{\hat{\mu}\nabla_0}
 + \left(\tilde{\nabla}_0^{\sigma_M} \circ u_{\alpha}^{d_{\alpha}\nabla}
 -u_{\alpha}^{d_{\alpha}\nabla} \circ \tilde{\nabla}_0^{\sigma_M} \right)
= v_{\alpha}^{\hat{\mu}\nabla_0}
 - \left(\tilde{\nabla} \circ u_{\alpha}^{d_{\alpha}\nabla_0}
 -u_{\alpha}^{d_{\alpha}\nabla_0} \circ \tilde{\nabla} \right).
\end{aligned}
\end{equation*}
Note that $d_{\alpha} + \left< d_{\alpha} ,  \tilde{A}^0_{\alpha} \frac{d z_{\alpha}}{z_{\alpha}} \right>$ 
has no pole at the support of $D( \tilde{\boldsymbol{p}} )$.
Moreover for any $a \in \tilde{E}|_{U_{\alpha}}$,
\begin{equation*}
\begin{aligned}
&\left(\tilde{\nabla} \circ u_{\alpha}^{d_{\alpha}\nabla_0}
 -u_{\alpha}^{d_{\alpha}\nabla_0} \circ \tilde{\nabla} \right)(a)\\
& =
 (\bar{\phi}_{\alpha}^{-1}\otimes \mathrm{id}) \circ \left(
d\left(\left< \hat{d}_{\alpha} ,   
\tilde{A}^0_{\alpha} \frac{d z_{\alpha}}{z_{\alpha}}
-\tilde{A}_{\alpha} \frac{d f_{\alpha}}{f_{\alpha}}\right>\right) 
( \bar{\phi}_{\alpha} (a)) 
 +\left[  \tilde{A}_{\alpha}\frac{d f_{\alpha}}{f_{\alpha}} ,
\left< \hat{d}_{\alpha} , \tilde{A}^0_{\alpha} \frac{d z_{\alpha}}{z_{\alpha}} \right>\right]
 ( \bar{\phi}_{\alpha} (a)) \right)
\end{aligned}
\end{equation*}
has no pole at the support of $D( \tilde{\boldsymbol{p}} )$.
Then the class $[(\{ {}'u_{\alpha\beta}^{\hat{\mu}\nabla_0}\}, \{  -  {}' v_{\alpha}^{\hat{\mu}\nabla_0}\})]$
coincides with
the class $[(\{ u_{\alpha\beta}^{\hat{\mu}\nabla_0}\}, \{  -   v_{\alpha}^{\hat{\mu}\nabla_0}\})]$
in $\bH^1 (\cG^{\bullet}_{M})$.
So the class $[(\{ u_{\alpha\beta}^{\hat{\mu}\nabla_0}\}, \{  -   v_{\alpha}^{\hat{\mu}\nabla_0}\})]$ 
is independent of the choice of a representative of 
the class $[ \{ \hat{d}_{\alpha\beta} \} ] 
\in H^1(\cC_M, \Theta_{\cC_M/M}( - D(\tilde{\bt})-  D( \tilde{\boldsymbol{p}} )))$.
\end{proof}

Now we consider meaning of the vector field 
$[ ( \{ 
u_{\alpha\beta}^{\hat{\mu}\nabla_0}
 \}_{\alpha\beta}, \{  
-v_{\alpha}^{\hat{\mu}\nabla_0}
\}_{\alpha} ) ]\in \bH^1(\cG^{\bullet}_M)$
by using the morphisms 
$$h_{\nabla-\nabla_0}\colon M \lra M'_H 
\quad \text{and} \quad 
 h_{\nabla_0}\colon M \lra M'$$
(see \eqref{2020.5.7.21.58} and \eqref{2020.5.7.21.59}).
In some sense,
we define tangent sheaves $\Theta_{M'_H}$ and $\Theta_{M'}$ of 
the moduli spaces $M'_H$ and $M'$, respectively.
Let $v_1 \in h^*_{\nabla-\nabla_0} \Theta_{M'_H}$ be the image
of the vector field of isomonodromic deformations of $\tilde{\nabla}$
associated to $\mu$ on $M$
under the morphism 
$\Theta_M \rightarrow h_{\nabla-\nabla_0}^*\Theta_{M'_H}$.
First, we consider 
the infinitesimal deformation of $\tilde{\nabla} - \tilde{\nabla}_0^{\sigma_M}$
on $M[\epsilon]$ corresponding to $v_1$.
Recall that $\{v_{\alpha}^{\hat{\mu}\nabla_0}\}_{\alpha}$
is defined by $\nabla_{0,\epsilon}^{\sigma_M,\text{IMD}}$,
which is the pull-back of $\nabla^{\sigma_M}_0$ under 
the morphism $\mathrm{IM}^{\mu}_M\colon \cC_{\epsilon} \rightarrow \cC_M$.
Then the infinitesimal deformation of $\tilde{\nabla} - \tilde{\nabla}_0^{\sigma_M}$
on $M[\epsilon]$ corresponding to $v_1$ 
is 
$(\mathrm{IM}^{\mu}_M)^* \tilde{\nabla}
- \nabla_{0,\epsilon}^{\sigma_M,\text{IMD}}$.
By
$\hat{v}_{\alpha}^{\text{IMD}} =0$ in \eqref{2019.7.31.11.02} and 
the definition \eqref{2019.7.31.11.14} of $v_{\alpha}^{\hat{\mu}\nabla_0}$, 
we have the following equality
\begin{equation*}
 (\varphi_{\alpha}^{-1} \otimes \mathrm{id}) 
 \circ ( (\mathrm{IM}^{\mu}_M)^* \tilde{\nabla}
  -   \nabla_{0,\epsilon}^{\sigma_M,\text{IMD}} )|_{U^{\epsilon}_{\alpha}}
   \circ \varphi_{\alpha}
-  (\tilde{\nabla} - \tilde{\nabla}_0^{\sigma_M})
= - v_{\alpha}^{\hat{\mu}\nabla_0}. 
\end{equation*}
Roughly speaking, the pair
$( \{ 
\hat{u}_{\alpha\beta}^{\mathrm{IMD}}
 \}_{\alpha\beta}, \{  
-v_{\alpha}^{\hat{\mu}\nabla_0}
\}_{\alpha} )$
means the infinitesimal deformation of 
$(\tilde{E},\tilde{\nabla}-\tilde{\nabla}_0^{\sigma}, 
\boldsymbol{l}^{\sigma_M} \cup\{ \tilde{l}_*^{(i)} \} )$
corresponding to
 the pushforward $v_1$.
Second, we will define $v_2 \in h^*_{\nabla-\nabla_0} \Theta_{M'_H}$ 
later at \eqref{2020.5.9.13.25}
and we consider the difference $v_1-v_2$ (see \eqref{2020.5.9.13.41} below).
We have the following equality
\begin{equation*}
\hat{u}_{\alpha\beta}^{\text{IMD}} -u_{\alpha\beta}^{\hat{\mu}\nabla_0}
= \bar{\phi}_{\beta}^{-1}\circ \left(  
\left< \hat{d}_{\alpha\beta} ,  \tilde{A}_{\beta} \frac{d f_{\beta}}{f_{\beta}}
-\tilde{A}^0_{\beta} \frac{d z_{\beta}}{z_{\beta}} \right> \right)
\circ \bar{\phi}_{\beta} .
\end{equation*}
Here $\tilde{A}_{\beta} \frac{d f_{\beta}}{f_{\beta}}$ and 
$\tilde{A}^0_{\beta} \frac{d z_{\beta}}{z_{\beta}}$
are the connection matrices of $\tilde{\nabla}$ and $\tilde{\nabla}_0^{\sigma_M}$
associated to the trivialization $\bar{\phi}_{\beta}$, respectively.
Then we have
$[ \tilde{\nabla} - \tilde{\nabla}_0^{\sigma_M},
\hat{u}_{\alpha\beta}^{\text{IMD}} -u_{\alpha\beta}^{\hat{\mu}\nabla_0}]=0.$
Here $[\cdot ,\cdot]$ is the commutator.
This equality means that the pair
\begin{equation}\label{2020.5.9.13.25}
( \{ \hat{u}_{\alpha\beta}^{\text{IMD}} -u_{\alpha\beta}^{\hat{\mu}\nabla_0}
 \}_{\alpha\beta}, \{  0\}_{\alpha} ) 
 \end{equation}
gives an infinitesimal deformation of 
$(\tilde{E} , \tilde{\nabla}-  \tilde{\nabla}_{0}^{\sigma_M},
\boldsymbol{l}^{\sigma_M} \cup\{ \tilde{l}_*^{(i)} \})$ 
parametrized by $M[\epsilon]$.
We denote by $v_2$ the corresponding vector field.
We consider the difference $v_1 - v_2$, which is the image of 
the vector field of isomonodromic deformations of $\tilde{\nabla}$ 
associated to $\mu$
under the following composition:
\begin{equation*}
\xymatrix{
\Theta_M \ar[r] \ar[dr] & h^*_{\nabla-\nabla_0} \Theta_{M'_H}  \ar[d]^{\text{adding $- v_2$}}\\
&h^*_{\nabla-\nabla_0} \Theta_{M'_H}\rlap{.}
}
\end{equation*} 
The difference $v_1 - v_2$ is described by
\begin{equation}\label{2020.5.9.13.41}
( \{ u_{\alpha\beta}^{\hat{\mu}\nabla_0} \}_{\alpha\beta}, 
\{  -v_{\alpha}^{\hat{\mu}\nabla_0}\}_{\alpha} ),
\end{equation}
which gives an infinitesimal deformation of 
$(\tilde{E} , \tilde{\nabla}-  \tilde{\nabla}_{0}^{\sigma_M},
\boldsymbol{l}^{\sigma_M} \cup\{ \tilde{l}_*^{(i)} \})$ 
parametrized by $M[\epsilon]$.
Third, 
we construct an infinitesimal deformation of 
$\tilde{\nabla}=\tilde{\nabla}_0^{\sigma_M} + (\tilde{\nabla}- \tilde{\nabla}_0^{\sigma_M})$
by the infinitesimal deformation of $\tilde{\nabla}_0^{\sigma_M}$
corresponding to the isomonodromic deformation of $\tilde{\nabla}_0^{\sigma_M}$
and the infinitesimal deformation of $\tilde{\nabla}- \tilde{\nabla}_0^{\sigma_M}$
corresponding to $v_1-v_2$.
For this purpose, we consider the following commutative diagram:
\begin{equation}\label{2020.11.7.11.28}
\xymatrix{
M' \times_{\text{par-Bun}} M'_H \ar[r]^-{\text{Add}} &  M'  \\
 & M \rlap{,}
   \ar[u]_-{h_{\nabla} } \ar[lu]^-{(h_{\nabla_0},h_{\nabla-\nabla_0}) \quad } 
}
\end{equation}
where $\text{par-Bun}$ is the moduli space of parabolic bundles (in some sense),
\begin{equation*}
\begin{aligned}
(h_{\nabla_0},h_{\nabla-\nabla_0}) \colon 
(E, \tilde{\nabla} , \{ \tilde{l}_*^{(i)} \} ; 
 \tilde{\nabla}_0^{\sigma_M} , \boldsymbol{l}^{\sigma_M})
 & \longmapsto 
 ((\tilde{E} ,  
\boldsymbol{l}^{\sigma_M} \cup\{ \tilde{l}_*^{(i)} \}),
\tilde{\nabla}_{0}^{\sigma_M}, \tilde{\nabla}-  \tilde{\nabla}_{0}^{\sigma_M}), \\
h_{\nabla} \colon (E, \tilde{\nabla} , \{ \tilde{l}_*^{(i)} \} ; 
 \tilde{\nabla}_0^{\sigma_M} , \boldsymbol{l}^{\sigma_M})
 &\longmapsto (E, \tilde{\nabla} , \boldsymbol{l}^{\sigma_M} \cup \{ \tilde{l}_*^{(i)} \} ) ,
\end{aligned}
\end{equation*}
and
\begin{equation*}
\begin{aligned}
\text{Add} ((\tilde{E} ,  
\boldsymbol{l}^{\sigma_M} \cup\{ \tilde{l}_*^{(i)} \}),
\tilde{\nabla}_{0}^{\sigma_M}, \tilde{\nabla}-  \tilde{\nabla}_{0}^{\sigma_M})
&=(E, \tilde{\nabla}_{0}^{\sigma_M}+(\tilde{\nabla}-  \tilde{\nabla}_{0}^{\sigma_M}) , 
\boldsymbol{l}^{\sigma_M}  \cup \{ \tilde{l}_*^{(i)} \} )\\
&=(E, \tilde{\nabla} , \boldsymbol{l}^{\sigma_M} \cup \{ \tilde{l}_*^{(i)} \} ).
\end{aligned}
\end{equation*}
We have a lift $M[\epsilon] \rightarrow M' \times_{\mathrm{par-Bun}} M'_H$
 of $(h_{\nabla_0},h_{\nabla-\nabla_0})
\colon M \rightarrow M' \times_{\mathrm{par-Bun}} M'_H$
by the infinitesimal deformation of $\tilde{\nabla}_0^{\sigma_M}$
corresponding to the isomonodromic deformation of $\tilde{\nabla}_0^{\sigma_M}$
and the infinitesimal deformation of $\tilde{\nabla}- \tilde{\nabla}_0^{\sigma_M}$
corresponding to $v_1-v_2$.
By the diagram \eqref{2020.11.7.11.28},
we have a lift $M[\epsilon] \rightarrow M' $
 of $h_{\nabla} \colon M \rightarrow M' $.
 The desired infinitesimal deformation of 
$\tilde{\nabla}=\tilde{\nabla}_0^{\sigma_M} + (\tilde{\nabla}- \tilde{\nabla}_0^{\sigma_M})$
just corresponds to this lift $M[\epsilon] \rightarrow M' $.
We describe this infinitesimal deformation of 
$\tilde{\nabla}=\tilde{\nabla}_0^{\sigma_M} + (\tilde{\nabla}- \tilde{\nabla}_0^{\sigma_M})$
by the \v{C}ech cohomology as follows 
(The conclusion is \eqref{2020.11.7.11.53} bellow).
A description of the infinitesimal deformation of $\tilde{\nabla}- \tilde{\nabla}_0^{\sigma_M}$
corresponding to $v_1-v_2$ by the \v{C}ech cohomology 
is \eqref{2020.5.9.13.41}.
We consider the infinitesimal deformation 
of $\tilde{\nabla}_0^{\sigma_M}$
corresponding to the isomonodromic deformation 
of $\tilde{\nabla}_0^{\sigma_M}$. 
We will give a description of this infinitesimal deformation 
of $\tilde{\nabla}_0^{\sigma_M}$ 
by the \v{C}ech cohomology (\eqref{2020.5.9.16.28} bellow).
As in Section \ref{SS IMD by c c}, we may define 
the vector field of the isomonodromic deformation for 
$\varpi' \colon M' \rightarrow \mathcal{M}_{g, n + N}$.
That is, we define a splitting $(\varpi')^* \Theta_{\mathcal{M}_{g, n + N}} \rightarrow \Theta_{M'}$.
By $h_{\nabla_0}$, 
we have a morphism 
$h_{\nabla_0}^* (\varpi')^* \Theta_{\mathcal{M}_{g, n + N}} \rightarrow h_{\nabla_0}^* \Theta_{M'}$.
The lift $\hat{\mu}$ of $\mu$ gives a section of 
$h_{\nabla_0}^* (\varpi')^* \Theta_{\mathcal{M}_{g, n + N}}$.
We take the image of $\hat{\mu}$ by this
morphism $h_{\nabla_0}^* (\varpi')^* \Theta_{\mathcal{M}_{g, n + N}} 
\rightarrow h_{\nabla_0}^* \Theta_{M'}$.
This image gives
an infinitesimal deformation of $\tilde{\nabla}_0^{\sigma_M}$ 
parametrized by $M[\epsilon]$.
By this infinitesimal deformation of $\tilde{\nabla}_0^{\sigma_M}$,
we have a pair 
\begin{equation}\label{2020.5.9.16.28}
( \{ u_{\alpha\beta}^{\hat{\mu}\nabla_0} \}_{\alpha\beta}, \{  0 \}_{\alpha} )
\end{equation}
as Proposition \ref{Prop explicit description IMD}.
The first components of the pairs \eqref{2020.5.9.13.41} and 
\eqref{2020.5.9.16.28} mean infinitesimal deformations 
of the underlying parabolic bundles.
The pairs \eqref{2020.5.9.13.41} and \eqref{2020.5.9.16.28}
have same first components.
Then we can add the infinitesimal 
deformation of $\tilde{\nabla}_0^{\sigma_M}$ 
corresponding to the pair \eqref{2020.5.9.16.28}
to the infinitesimal 
deformation of $\tilde{\nabla} - \tilde{\nabla}_0^{\sigma_M}$ 
corresponding to the pair \eqref{2020.5.9.13.41}.
By this addition, 
we have a pair 
\begin{equation}\label{2020.11.7.11.53}
( \{ u_{\alpha\beta}^{\hat{\mu}\nabla_0} \}_{\alpha\beta}, 
\{  -v_{\alpha}^{\hat{\mu}\nabla_0}\}_{\alpha} ).
\end{equation}
This pair corresponds to the lift $M[\epsilon] \rightarrow M'$.
By this infinitesimal deformation of $\tilde{\nabla}$ and 
forgetting the parabolic structure $\boldsymbol{l}^{\sigma_M}$,
 we have a morphism $M[\epsilon] \rightarrow M$.
Here this morphism is given by the universal property of the moduli space $M$.
Finally, this morphism $M[\epsilon] \rightarrow M$
corresponds to
the vector field 
$$[ ( \{ 
u_{\alpha\beta}^{\hat{\mu}\nabla_0} \}_{\alpha\beta}, \{  
-v_{\alpha}^{\hat{\mu}\nabla_0}\}_{\alpha} ) ]\in \bH^1(\cG^{\bullet}_M).$$

\subsection{Trivializations 
of the infinitesimal deformation of $\tilde{E}$
given by the infinitesimal deformation of $\sigma_M$}

We take any $v \in \Theta_M$.
For the vector field $v$, we have a morphism $f_v\colon M[\epsilon] \rightarrow M$.
Put $\cC_{\epsilon}^v:= \cC_M \times_M M[\epsilon]$ given by
the family $\cC_M \rightarrow M$ 
and $f_v\colon M[\epsilon] \rightarrow M$.
We denote by $F_v \colon \cC_{\epsilon}^v \rightarrow \cC_M$
the natural morphism from $\cC_{\epsilon}^v$ to $\cC_M$.
For $\mathrm{symb}_1(v) \in H^1 (\cC_M, \Theta_{C_M/M}(-D(\tilde{\bt})))$,
we take the lift of $\mathrm{symb}_1(v)$
in $ H^1(\cC_M, \Theta_{\cC_M/M}( - D(\tilde{\bt}) 
- D( \tilde{\boldsymbol{p}} )))$ defined by the family
$(\cC_{\epsilon}^v, \mathrm{Supp}(F_v^* (D(\tilde{\bt}) + D( \tilde{\boldsymbol{p}} ))))$.
Let $\{U_{\alpha}\}_{\alpha}$ be an affine open covering of $\cC_M$.
Set $U_{\alpha}^{\epsilon}:=F_v^{-1}( U_{\alpha})$.
Let $\sigma_{\alpha} \colon  U_{\alpha}^{\epsilon} \rightarrow U_{\alpha}
\times \mathrm{Spec}\, \mathbb{C}[\epsilon]$ be 
an isomorphism corresponding to the lift.
The infinitesimal deformation
$(F_v^*\tilde{E}, F_v^* \tilde{\nabla}_0^{\sigma_M} )$
is defined by taking the pull-back of $(\tilde{E} , \tilde{\nabla}_0^{\sigma_M}) $
under the morphism $F_v \colon \cC_{\epsilon}^v \rightarrow \cC_M$.
Let $\sigma_M \colon  \mathcal{O}^{\oplus r}_{\cC_M} \rightarrow \tilde{E}(m)|_{\cC_M}$
be the injection in Proposition \ref{Prop nabla0Wr}.
Let $F_v^*\sigma_{M}$ be the pull-back 
of $\sigma_M \in 
\mathrm{Hom} (
\mathcal{O}^{\oplus r}_{\cC_M} (-m), \tilde{E}|_{\cC_M})$
under the morphism $F_v \colon \cC_{\epsilon}^v \rightarrow \cC_M$:
\begin{equation}\label{2019.7.30.15.48}
F_v^*\sigma_{M} \colon \ \cO_{\cC^v_{\epsilon}}^{\oplus r}(-m) 
\lra F_v^*\tilde{E}.
\end{equation}
In this section we will define trivializations of 
$F_v^*\tilde{E}$
which are compatible with $F_v^* \sigma_{M}$.
We will use these trivializations in the proof of Lemma \ref{2020.5.24.14.16}.
Let $D_{\epsilon}(m)=\{ (U_{\alpha}^{\epsilon}, g_{\alpha}^{\epsilon}) \}_{\alpha}$ and
$D_{\epsilon}( \sigma )=
\{ (U_{\alpha}^{\epsilon}, h_{\alpha}^{\epsilon}) \}_{\alpha}$ be 
the pull-backs of the Cartier divisors 
$D(m)=\{ (U_{\alpha}, g_{\alpha}) \}_{\alpha}$  
and $D( \sigma )=\{ (U_{\alpha}, h_{\alpha}) \}_{\alpha}$ 
by the morphism $F_v \colon \cC_{\epsilon}^v \rightarrow \cC_M$, respectively.

We define trivializations $\{(U_{\alpha}^{\epsilon}, \psi_{\alpha})\}_{\alpha}$
of $F_v^*\tilde{E}$ as follows.
First, we consider an affine open subset $U_{\alpha}$
such that the intersection of $U_{\alpha}$ and the support of $D(\tilde{\boldsymbol{p}})$
is empty.
On such an affine open subset $U_{\alpha}^{\epsilon}$,
the natural injection 
$\mathcal{O}^{\oplus r}_{U_{\alpha}^{\epsilon}}(-m)
\rightarrow \mathcal{O}^{\oplus r}_{U_{\alpha}^{\epsilon}} $ 
and $F_v^*\sigma_{M} \colon 
\mathcal{O}^{\oplus r}_{U_{\alpha}^{\epsilon}}(-m)
\rightarrow F_v^*\tilde{E}|_{U_{\alpha}^{\epsilon}}$ are isomorphisms.
Then the morphisms
\begin{equation*}
\mathcal{O}^{\oplus r}_{U_{\alpha}^{\epsilon}} 
\xleftarrow[\cong]{\ \supset \ }
\mathcal{O}^{\oplus r}_{U_{\alpha}^{\epsilon}}(-m)
\xrightarrow[\cong]{\ F_v^*\sigma_{M}\ }  F_v^*\tilde{E}|_{U_{\alpha}^{\epsilon}}
\end{equation*}
give a trivialization of $F_v^*\tilde{E}|_{U_{\alpha}^{\epsilon}}$.
We denote by $(U_{\alpha}^{\epsilon}, \psi_{\alpha})$
this trivialization of $F_v^*\tilde{E}|_{U_{\alpha}^{\epsilon}}$.
Second, 
we consider an affine open set $U_{\alpha}^{\epsilon}$
such that 
$U_{\alpha}^{\epsilon} \cap \mathrm{Supp}(D_{\epsilon}(m)) \neq \emptyset$
and
$U_{\alpha}^{\epsilon} \cap \mathrm{Supp}(D_{\epsilon}(\sigma_M)) = \emptyset$.
On such an affine open subset $U_{\alpha}^{\epsilon}$,
the composition 
\begin{equation*}
\mathcal{O}^{\oplus r}_{U_{\alpha}^{\epsilon}} 
\xrightarrow[\cong]{\ g_{\alpha}^{\epsilon} \ }
\mathcal{O}^{\oplus r}_{U_{\alpha}^{\epsilon}}(-m)
\xrightarrow[\cong]{\ F_v^*\sigma_{M}\ } F_v^*\tilde{E}|_{U_{\alpha}^{\epsilon}}
\end{equation*}
is an isomorphism.
This isomorphism 
gives a trivialization of $F_v^*\tilde{E}|_{U_{\alpha}^{\epsilon}}$.
We denote by $(U_{\alpha}^{\epsilon}, \psi_{\alpha})$
this trivialization of $F_v^*\tilde{E}|_{U_{\alpha}^{\epsilon}}$.
Third, we consider an affine open set $U_{\alpha}^{\epsilon}$
where $U_{\alpha}^{\epsilon} \cap \mathrm{Supp}(D_{\epsilon}(m)) = \emptyset$
and
$U_{\alpha}^{\epsilon} \cap \mathrm{Supp}(D_{\epsilon}(\sigma_M)) \neq \emptyset$.
Set $(\tilde{\boldsymbol{p}}'_{i''})^{\epsilon}:= 
\tilde{\boldsymbol{p}}'_{i''} \times_{M}\Spec \cO_{M}[\epsilon]$.
Here $\{\tilde{\boldsymbol{p}}'_1,\ldots,\tilde{\boldsymbol{p}}'_{N_2}\}$ are 
the components of
the support of the divisor $D(\sigma_M)$,
where $N_2:=\deg D(\sigma_M)|_{\cC_x}$
for each $x \in M$.
We put $\boldsymbol{s}_{\alpha,j}^{\epsilon}:= 
(F_v^*\sigma_{M})|_{U_{\alpha}^{\epsilon}}
(\boldsymbol{e}_j\otimes g^{\epsilon}_{\alpha})$ for $j=1,\ldots,r$.
We change the order of 
$( \boldsymbol{s}^{\epsilon}_{\alpha,1} , \boldsymbol{s}^{\epsilon}_{\alpha,2},\ldots
,\boldsymbol{s}^{\epsilon}_{\alpha,r})$ as follows.
For each $i''$,
we set 
$ \boldsymbol{s}^{\epsilon}_{\alpha,i'',k} := \boldsymbol{s}^{\epsilon}_{\alpha,j_{(i'',k)}}$
where $k=1,2,\ldots,r$ and $\{ j_{(i'',1)} ,\ldots,j_{(i'',r)}\} =\{ 1,2,\ldots,r\}$.
We assume that 
$\boldsymbol{s}^{\epsilon}_{\alpha,i'',2}
|_{(\tilde{\boldsymbol{p}}'_{i''})^{\epsilon}},
\ldots, \boldsymbol{s}^{\epsilon}_{\alpha,i'',r}
|_{(\tilde{\boldsymbol{p}}'_{i''})^{\epsilon}}$
are linearly independent. 
Let $a^{\epsilon}_{\alpha,i'',2},\ldots,a^{\epsilon}_{\alpha,i'',r}$ be elements 
of $\pi_{M[\epsilon]}^*(\mathcal{O}_{M[\epsilon]})|_{U^{\epsilon}_{\alpha}}$
such that  
$$ \boldsymbol{s}^{\epsilon}_{\alpha,i'',1}
|_{(\tilde{\boldsymbol{p}}'_{i''})^{\epsilon}}
= a_{\alpha,i',2}^{\epsilon}\cdot \boldsymbol{s}^{\epsilon}_{\alpha,i'',2}
|_{(\tilde{\boldsymbol{p}}'_{i''})^{\epsilon}}
+\cdots+ 
a^{\epsilon}_{\alpha,i',r} \cdot \boldsymbol{s}^{\epsilon}_{\alpha,i'',r}
|_{(\tilde{\boldsymbol{p}}'_{i''})^{\epsilon}}.$$
Here $\pi_{M[\epsilon]} \colon \cC^v_{\epsilon} \rightarrow M[\epsilon]$
is the projection.
We define $\hat{\boldsymbol{s}}^{\epsilon}_{\alpha,i'',1}$
as 
$$\boldsymbol{s}^{\epsilon}_{\alpha,i'',1}
= h^{\epsilon}_{\alpha} \cdot \hat{\boldsymbol{s}}^{\epsilon}_{\alpha,i'',1}
+ a^{\epsilon}_{\alpha,i'',2}\cdot\boldsymbol{s}_{\alpha,i'',2}^{\epsilon}
+\cdots+ a^{\epsilon}_{\alpha,i'',r}\cdot\boldsymbol{s}^{\epsilon}_{\alpha,i'',r}.$$
The tuple $\hat{\boldsymbol{s}}_{\alpha,i''}^{\epsilon}:=
(\hat{\boldsymbol{s}}^{\epsilon}_{\alpha,i'',1},
\boldsymbol{s}_{\alpha,i'',2}^{\epsilon},\ldots,
\boldsymbol{s}^{\epsilon}_{\alpha,i'',r})$ gives 
a trivialization of $F_v^*\tilde{E}|_{U_{\alpha}^{\epsilon}}$.
We denote by $(U_{\alpha}^{\epsilon}, \psi_{\alpha})$
this trivialization of $F_v^*\tilde{E}|_{U_{\alpha}^{\epsilon}}$.

\begin{Def}
We say trivializations $\{(U_{\alpha}^{\epsilon}, \psi_{\alpha})\}_{\alpha}$
of $F_v^*\tilde{E}$
are \textit{compatible with $F_v^*\sigma_{M}$}
if $\{(U_{\alpha}^{\epsilon}, \psi_{\alpha})\}_{\alpha}$ are the trivializations 
constructed above for any $\alpha$.
\end{Def}

We set $\bar{\psi}_{\alpha}:= \psi_{\alpha}$ ($\mathrm{mod}\ \epsilon$).
For any $\alpha$,
we define $\varphi_{\alpha}' \colon \tilde{E}|_{U_{\alpha}}
\rightarrow F_v^*\tilde{E}|_{U^{\epsilon}_{\alpha}}$
by the composition
\begin{equation}\label{2020.5.15.12.53}
\varphi_{\alpha}' \colon \tilde{E}|_{U_{\alpha}} \otimes_{\mathbb{C}} \mathbb{C}[\epsilon] 
\xrightarrow{\bar{\psi}_{\alpha}\otimes 1}  
\cO^{\oplus r}_{U_{\alpha}} \otimes_{\mathbb{C}} \mathbb{C}[\epsilon] 
= 
\cO^{\oplus r}_{U_{\alpha}} \otimes_{\cO_{U_{\alpha}}} (\cO_{U_{\alpha}} 
\otimes_{\mathbb{C}} \mathbb{C}[\epsilon] )
\xrightarrow{1 \otimes \sigma_{\alpha}^*} 
\cO^{\oplus r}_{U_{\alpha}} \otimes_{\cO_{U_{\alpha}}} \cO_{U^{\epsilon}_{\alpha}} 
\xrightarrow{\psi_{\alpha}^{-1}} 
F_v^*\tilde{E}|_{U^{\epsilon}_{\alpha}}
\end{equation} 
as \eqref{isom varphi}.

\begin{Lem}\label{2020.5.15.11.56}
\textit{
For the isomorphisms $\sigma_{\alpha}$ $(\alpha \in I_{D( \tilde{\boldsymbol{p}} )})$,
we assume that $(\sigma_{\alpha}^*)^{-1}g_{\alpha}^{\epsilon} = g_{\alpha}$
and $(\sigma_{\alpha}^*)^{-1}h_{\alpha}^{\epsilon} = h_{\alpha}$.
Moreover assume that
the trivializations $\{(U_{\alpha}^{\epsilon}, \psi_{\alpha})\}_{\alpha}$
of $F_v^*\tilde{E}$
are compatible with $F_v^*\sigma_{M}$.
Then}
\begin{equation}\label{2020.5.16.13.00}
( (\varphi'_{\alpha})^{-1} \otimes \mathrm{id} ) 
\circ (  F_v^* \tilde{\nabla}_0^{\sigma_M} ) 
\circ   \varphi'_{\alpha} -  \tilde{\nabla}_0^{\sigma_M}=0
\end{equation}
\textit{for any $\alpha$.}
\end{Lem}

\begin{proof}
First, we consider an affine open subset $U_{\alpha}$
such that the intersection of $U_{\alpha}$ and the support of $D(\tilde{\boldsymbol{p}})$
is empty.
On such an affine open subset $U_{\alpha}^{\epsilon}$,
we have the following commutative diagram
\begin{equation*}
\xymatrix{
\mathcal{O}^{\oplus r}_{U_{\alpha}^{\epsilon}} 
\ar[d]^-{d_{\cC_{\epsilon}^v/M[\epsilon]}|_{U_{\alpha}^{\epsilon}}  } &
\ar[l]^{\cong}_{\supset}
\mathcal{O}^{\oplus r}_{U_{\alpha}^{\epsilon}}(-m) 
\ar[r]_{\cong}^{F_v^*\sigma_{M}} 
\ar[d]^-{F_v^* d_m|_{U_{\alpha}^{\epsilon}}} &
  F_v^*\tilde{E}|_{U_{\alpha}^{\epsilon}}
  \ar[d]^-{F_v^* \tilde{\nabla}_0^{\sigma_M}|_{U^{\epsilon}_{\alpha}}} \\
  \mathcal{O}^{\oplus r}_{U_{\alpha}^{\epsilon}} 
  \otimes \Omega^1_{U^{\epsilon}_{\alpha}/M[\epsilon]} (D_{\epsilon}(\tilde{\boldsymbol{p}})) &
\ar[l]^-{\cong}_-{\supset}
\mathcal{O}^{\oplus r}_{U_{\alpha}^{\epsilon}}(-m) 
\otimes \Omega^1_{U^{\epsilon}_{\alpha}/M[\epsilon]} (D_{\epsilon}(\tilde{\boldsymbol{p}}))
\ar[r]_{\cong}^{F_v^*\sigma_{M}} &
  F_v^*\tilde{E}|_{U_{\alpha}^{\epsilon}}
  \otimes \Omega^1_{U^{\epsilon}_{\alpha}/M[\epsilon]} (D_{\epsilon}(\tilde{\boldsymbol{p}}))
  \rlap{.}
}
\end{equation*}
By this commutative diagram and the definition of 
the trivialization $(U_{\alpha}^{\epsilon}, \psi_{\alpha})$,
we have 
$$(\psi_{\alpha}\otimes \mathrm{id} )\circ(F_v^*\tilde{\nabla}^{\sigma_{M}}_0|_{U^{\epsilon}_{\alpha}}
 )\circ \psi_{\alpha}^{-1}=
d_{\cC_{\epsilon}^v/M[\epsilon]}|_{U_{\alpha}^{\epsilon}}. $$
By this equality, we may check the equality \eqref{2020.5.16.13.00}.

Second, 
we consider an affine open set $U_{\alpha}^{\epsilon}$
such that 
$U_{\alpha}^{\epsilon} \cap \mathrm{Supp}(D_{\epsilon}(m)) \neq \emptyset$
and
$U_{\alpha}^{\epsilon} \cap \mathrm{Supp}(D_{\epsilon}(\sigma_M)) = \emptyset$.
By the definition of the trivialization $(U_{\alpha}^{\epsilon}, \psi_{\alpha})$,
we may check that 
\begin{equation*}\label{dmdescrip}
(\psi_{\alpha}\otimes \mathrm{id} )\circ(F_v^*\sigma_{M}|_{U^{\epsilon}_{\alpha}}
 )\circ \psi_{\alpha}^{-1}
= d_{\cC_{\epsilon}^v/M[\epsilon]}|_{U_{\alpha}^{\epsilon}}+ 
 \mathrm{diag}(\frac{d g_{\alpha}^{\epsilon}}{g_{\alpha}^{\epsilon}} ,
 \ldots, \frac{d g_{\alpha}^{\epsilon}}{g_{\alpha}^{\epsilon}}).
\end{equation*}
By this equality and $(\sigma_{\alpha}^*)^{-1}g_{\alpha}^{\epsilon} = g_{\alpha}$,
we may check the equality \eqref{2020.5.16.13.00}.

Third, we consider an affine open set $U_{\alpha}^{\epsilon}$
where $U_{\alpha}^{\epsilon} \cap \mathrm{Supp}(D_{\epsilon}(m)) = \emptyset$
and
$U_{\alpha}^{\epsilon} \cap \mathrm{Supp}(D_{\epsilon}(\sigma_M)) \neq \emptyset$.
On such an affine open subset $U_{\alpha}^{\epsilon}$,
we have the following commutative diagram
\begin{equation*}
\xymatrix{
\mathcal{O}^{\oplus r}_{U_{\alpha}^{\epsilon}} (-m)
\ar[d]^-{F_v^* d_m|_{U_{\alpha}^{\epsilon}}}  &
\ar[l]^{\cong}_{g_{\alpha}^{\epsilon}}
\mathcal{O}^{\oplus r}_{U_{\alpha}^{\epsilon}}
\ar[r]^{T^{\epsilon}_{\alpha,i''}}  &
\mathcal{O}^{\oplus r}_{U_{\alpha}^{\epsilon}}
  \ar[d]_-{(\psi_{\alpha}\otimes \mathrm{id} ) \circ
(  F_v^* \tilde{\nabla}_0^{\sigma_M}) \circ \psi_{\alpha}^{-1}
   |_{U^{\epsilon}_{\alpha}}} \\
  \mathcal{O}^{\oplus r}_{U_{\alpha}^{\epsilon}}(-m) 
  \otimes \Omega^1_{U^{\epsilon}_{\alpha}/M[\epsilon]} (D_{\epsilon}(\tilde{\boldsymbol{p}})) &
\ar[l]^-{\cong}_-{g_{\alpha}^{\epsilon}}
\mathcal{O}^{\oplus r}_{U_{\alpha}^{\epsilon}} 
\otimes \Omega^1_{U^{\epsilon}_{\alpha}/M[\epsilon]} (D_{\epsilon}(\tilde{\boldsymbol{p}}))
\ar[r]^{T^{\epsilon}_{\alpha,i''}} &
\mathcal{O}^{\oplus r}_{U_{\alpha}^{\epsilon}}
  \otimes \Omega^1_{U^{\epsilon}_{\alpha}/M[\epsilon]} (D_{\epsilon}(\tilde{\boldsymbol{p}}))
}
\end{equation*}
(see the diagram \eqref{2020.5.6.13.07}).
Here we set
\begin{equation*}
T^{\epsilon}_{\alpha,i''}:=
\begin{pmatrix}
h^{\epsilon}_{\alpha} & 0  & 0& \cdots &0 \\
 a^{\epsilon}_{\alpha,i'',2} & 1 & 0 & \cdots & 0\\
a^{\epsilon}_{\alpha,i'',3} & 0& 1 & \cdots & 0 \\ 
\cdots & \cdots &\cdots &\cdots &\cdots  \\
a^{\epsilon}_{\alpha,i'',r} & 0& 0& \cdots & 1 
\end{pmatrix}.
\end{equation*}
By the commutative diagram,
we have
\begin{equation*}
\left(\psi_{\alpha} \otimes \mathrm{id}\right) 
\circ ( F_v^* \tilde{\nabla}_0^{\sigma_M} ) 
\circ \psi_{\alpha}^{-1} 
= d_{\cC_{\epsilon}^v/M[\epsilon]}|_{U_{\alpha}^{\epsilon}}+
T_{\alpha}^{\epsilon} d ((T_{\alpha}^{\epsilon})^{-1})
+
\mathrm{diag}(\frac{d g_{\alpha}^{\epsilon}}{g_{\alpha}^{\epsilon}},
 \ldots, \frac{d g_{\alpha}^{\epsilon}}{g_{\alpha}^{\epsilon}}).
\end{equation*}
By this equality, $(\sigma_{\alpha}^*)^{-1}g_{\alpha}^{\epsilon} = g_{\alpha}$, 
and $(\sigma_{\alpha}^*)^{-1}h_{\alpha}^{\epsilon} = h_{\alpha}$,
we may check the equality \eqref{2020.5.16.13.00}.
\end{proof}

\subsection{Hamiltonians on the moduli spaces}\label{SS Hamil}

Let $M_{\text{aff}}$ be an affine open set of 
$M_{\cC/T}^{\balpha}(\tilde{\bt} ,r,d)_{\boldsymbol{\nu}}$.
We define an algebraic splitting $\eta \colon \bH^1(\cG^{\bullet}_{M_{\text{aff}}}) 
\ra \bH^1(\cF^{\bullet}_{M_{\text{aff}}})$ of
the tangent map $\bH^1(\cF^{\bullet}_{M_{\text{aff}}}) 
\ra \bH^1(\cG^{\bullet}_{M_{\text{aff}}})$ by
\begin{equation*}
\eta \colon \bH^1(\cG^{\bullet}_{M_{\text{aff}}}) 
\lra \bH^1(\cF^{\bullet}_{M_{\text{aff}}}); \quad
[(\{ u_{\alpha\beta} \}, \{ v_{\alpha} \} )] 
\longmapsto  [(\{ \eta(u_{\alpha\beta}) \}, \{  v_{\alpha} \} )].
\end{equation*}
Here we set
\begin{equation}\label{2020.11.3.19.35}
\eta(s) := s - \iota(\tilde{\nabla}) \circ \mathrm{symb}_1(s) , \quad s \in \cG_{M_{\text{aff}}}^0 .
\end{equation}
First, we define a lift  
\begin{equation}\label{2019.7.28.14.00}
\omega \colon \bold{R}^1 (\pi_{M_{\cC/T}^{\balpha}(\tilde{\bt} ,r,d)_{\boldsymbol{\nu}}})_*
(\cG^{\bullet}) 
\otimes \bold{R}^1 (\pi_{M_{\cC/T}^{\balpha}(\tilde{\bt} ,r,d)_{\boldsymbol{\nu}}})_*
(\cG^{\bullet}) \lra 
\cO_{M_{\cC/T}^{\balpha}(\tilde{\bt} ,r,d)}
\end{equation}
of the symplectic form defined in
Section \ref{S2 symp form} as follows.
We define a pairing
\begin{equation}\label{2020.5.11.14.00}
\begin{aligned}
\omega_{M_{\text{aff}}} \colon \bH^1(\cC_{M_{\text{aff}}}, \cG_{M_{\text{aff}}}^{\bullet}) 
\otimes \bH^1(\cC_{M_{\text{aff}}}, \cG_{M_{\text{aff}}}^{\bullet}) 
&\lra \bH^2(\cC_{M_{\text{aff}}} , \Omega_{\cC_{M_{\text{aff}}}/{M_{\text{aff}}}}^{\bullet}) 
\cong H^0(\cO_{M_{\text{aff}}})\\
v \otimes v' &\longmapsto \omega_{M_{\text{aff}}}(v,v'),
\end{aligned}
\end{equation}
where we put 
\begin{equation*}
\omega_{M_{\text{aff}}}(v,v'):=
[(\{  \Tr( \eta( u_{\alpha\beta}) \circ \eta( u_{\beta\gamma}')) \}, 
-\{ \Tr (\eta( u_{\alpha\beta}) \circ  v_{\beta}' ) - \Tr (v_{\alpha} \circ \eta( u'_{\alpha\beta})) \} )]
\end{equation*}
for $v=[(\{ u_{\alpha\beta} \}, \{ v_{\alpha} \})]$
and $v'=[(\{ u_{\alpha\beta}' \}, \{ v_{\alpha}' \})]$.
Here we consider in \v{C}ech cohomology with respect to an affine open covering 
$\{ U_{\alpha} \}$ 
of $\cC \times_T{M_{\text{aff}}}$, 
$\{ u_{\alpha\beta} \} \in C^1(\cG_{M_{\text{aff}}}^0)$, 
$\{ v_{\alpha} \} \in C^0(\cG^1_{M_{\text{aff}}})$.
This pairing induces a lift $\omega$.

\begin{Prop}\label{Kernel is IMD}
\textit{
The kernel $\mathrm{Ker}(\omega_{M_{\text{aff}}})$
of $\omega_{M_{\text{aff}}}\colon \bH^1(\cG_{M_{\text{aff}}}^{\bullet}) 
\ra \mathcal{H}om_{\cO_{M_{\text{aff}}}} (\bH^1(\cG^{\bullet}_{M_{\text{aff}}}), 
\cO_{M_{\text{aff}}})$
 induces the vector fields on $M_{\text{aff}}$ determined by the isomonodromic deformations.
}
\end{Prop}

\begin{proof}

The algebraic splitting \eqref{map which means IMD} for 
$\bH^1(\cG^{\bullet}_{M_{\text{aff}}}) \rightarrow 
H^1(\Theta_{\cC\times_TM_{\text{aff}}/M_{\text{aff}}} 
(-D(\tilde{\bt})_{M_{\text{aff}}}  ))$
gives an isomorphism 
$$
\begin{aligned}
H^1(\Theta_{\cC\times_TM_{\text{aff}}/M_{\text{aff}}} 
(-D(\tilde{\bt})_{M_{\text{aff}}}  ))\oplus \bH^1(\cF^{\bullet}_{M_{\text{aff}}})
&\xrightarrow{\, \cong \, }
\bH^1(\cG^{\bullet}_{M_{\text{aff}}})   \\
 ( [\{ d_{\alpha\beta} \}]  ,
 [(\{ u_{\alpha\beta} \}, \{  v_{\alpha} \} )])
 &\longmapsto [ \{  \iota(\tilde{\nabla}) (d_{\alpha\beta}) + u_{\alpha \beta}\} , \{ v_{\alpha} \} ].
\end{aligned}
$$
By this isomorphism, we can define a composition 
\begin{equation}\label{2020.11.7.16.08}
H^1(\Theta_{\cC\times_TM_{\text{aff}}/M_{\text{aff}}} 
(-D(\tilde{\bt})_{M_{\text{aff}}}  )) \oplus \bH^1(\cG^{\bullet}_{M_{\text{aff}}}) 
\lra
 \bH^1( \cG_{M_{\text{aff}}}^{\bullet}) 
\otimes \bH^1( \cG_{M_{\text{aff}}}^{\bullet}) 
\xrightarrow{\omega_{M_{\text{aff}}}} 
 \bH^2( \Omega_{\cC_{M_{\text{aff}}}/{M_{\text{aff}}}}^{\bullet}) 
\end{equation}
The image of $[\{d_{\alpha\beta} \}] \in H^1(\Theta_{\cC\times_TM_{\text{aff}}/M_{\text{aff}}} 
(-D(\tilde{\bt})_{M_{\text{aff}}}  ))$ is 
$[ \{  \iota(\tilde{\nabla}) (d_{\alpha\beta}) \} , \{ 0 \} ]$.
Moreover, we have 
$$
\begin{aligned}
 \eta( [ \{  \iota(\tilde{\nabla}) (d_{\alpha\beta}) \} , \{ 0 \} ])
 &= [ \{  \iota(\tilde{\nabla}) (d_{\alpha\beta}) - \iota(\tilde{\nabla}) (d_{\alpha\beta}) \} , \{ 0 \} ] \\
 &=[ \{  0\} , \{ 0 \} ].
\end{aligned}
$$
By \eqref{2020.5.11.14.00}, the composition \eqref{2020.11.7.16.08} is 
the zero morphism. 
Then the image of $H^1(\Theta_{\cC\times_TM_{\text{aff}}/M_{\text{aff}}} 
(-D(\tilde{\bt})_{M_{\text{aff}}}  ))$ in $\bH^1(\cG^{\bullet}_{M_{\text{aff}}})$
is contained in
$\mathrm{Ker}(\omega_{M_{\text{aff}}})$.
On the other hand, the composition 
$$
\begin{aligned}
\bH^1(\cF^{\bullet}_{M_{\text{aff}}}) \oplus \bH^1(\cF^{\bullet}_{M_{\text{aff}}}) 
\lra
 \bH^1( \cG_{M_{\text{aff}}}^{\bullet}) 
\otimes \bH^1( \cG_{M_{\text{aff}}}^{\bullet}) 
\xrightarrow{\omega_{M_{\text{aff}}}} 
& \bH^2( \Omega_{\cC_{M_{\text{aff}}}/{M_{\text{aff}}}}^{\bullet}) \\
&\cong H^0(\cO_{M_{\text{aff}}}).
\end{aligned}
$$
is just the symplectic structure \eqref{2020.11.7.15.48} in Section \ref{S2 symp form}.
In particular, this pairing is nondegenerate. 
Then the image of $H^1(\Theta_{\cC\times_TM_{\text{aff}}/M_{\text{aff}}} 
(-D(\tilde{\bt})_{M_{\text{aff}}}  ))$ in $\bH^1(\cG^{\bullet}_{M_{\text{aff}}})$
coincides with
$\mathrm{Ker}(\omega_{M_{\text{aff}}})$.
By Proposition \ref{Prop explicit description IMD},
this image of $H^1(\Theta_{\cC\times_TM_{\text{aff}}/M_{\text{aff}}} 
(-D(\tilde{\bt})_{M_{\text{aff}}}  ))$
means the vector fields on $M_{\text{aff}}$ determined by the isomonodromic deformations.
\end{proof}

Second, we define Hamiltonian functions for
each affine open subset $M\subset M_{\cC/T}^{\balpha}(\tilde{\bt} ,r,d)_{\boldsymbol{\nu}}$
in Proposition \ref{Prop nabla0Wr}
as follows.
We take an affine open sebset $U\subset T$
and we take
$\mu_1 ,\ldots , \mu_{3g-3+n} \in  \Theta_{T}(U)$
such that $\mu_1 ,\ldots , \mu_{3g-3+n}$ give a trivialization 
of $\Theta_T|_U$. 
By taking a refined covering $\{ M\}$, we may define a morphism
$\varpi_M \colon M \rightarrow U$.
The sections $\mu_1 ,\ldots , \mu_{3g-3+n}$ 
induce classes of 
$H^1(\cC_M, \Theta_{\cC_M/M}( - D(\tilde{\bt})))$
by $\varpi_M \colon M \rightarrow U$.
We also denote by $\mu_1 ,\ldots , \mu_{3g-3+n}$
the corresponding classes in $H^1(\cC_M, \Theta_{\cC_M/M}( - D(\tilde{\bt})))$.
We take a representative $\{ (d_k)_{\alpha\beta} \}_{\alpha\beta}$
 of the class $\mu_k$:
\begin{equation}\label{2020.5.24.14.44}
\mu_k = \left[ \{ (d_k)_{\alpha\beta} \} \right] 
\in H^1(\cC_M, \Theta_{\cC_M/M}( - D(\tilde{\bt})))
\end{equation}
for each $k$ ($k=1,\ldots,3g-3+n$).
Let $\mathrm{im}_M^{\mu_k} \colon M[\epsilon] \rightarrow M$
be the isomonodromic lift associated to $\mu_k$
and $\mathrm{IM}_{M}^{\mu_k} \colon \cC^k_{\epsilon} \rightarrow \cC_M$
be the induced morphism by $\mathrm{im}_M^{\mu_k} \colon M[\epsilon] \rightarrow M$.
We take a representative $ \{ (\hat{d}_k)_{\alpha\beta} \}_{\alpha\beta} $ of the lift 
\begin{equation*}
\hat{\mu}_k =[\{(\hat{d}_k)_{\alpha\beta}\}_{\alpha\beta}] 
\in H^1(\cC_M, \Theta_{\cC_M/M}( - D(\tilde{\bt})-  D( \tilde{\boldsymbol{p}} ))) 
\end{equation*}
defined as \eqref{mulift}.

\begin{Def}\label{2020.5.22.11.40}
We describe
a representative $ \{ (\hat{d}_k)_{\alpha\beta} \}_{\alpha\beta} $ 
as
$(\hat{d}_k)_{\alpha\beta} = (d_k)_{\alpha\beta}-(\delta_k)_{\alpha}+(\delta_k)_{\beta},$
where 
$$(\delta_k)_{\alpha} = \tilde{\delta}_{k,\alpha} \frac{\partial}{\partial z_{\alpha}}
 \in  \Theta_{\cC_M/M}(-D(\tilde{\bt}))(U_{\alpha})$$ 
for any $\alpha$
and $(\delta_k)_{\alpha}=0$ for $\alpha\notin I_{D(\tilde{\boldsymbol{p}})}$.
Here $\tilde{\delta}_{k,\alpha}
 \in \pi_M^*( \mathcal{O}_M)|_{U_{\alpha}}$.
By changing of coordinates $\{ z_{\alpha}\}$
($\alpha\in I_{D(\tilde{\boldsymbol{p}})}$),
we may assume that 
\begin{equation}\label{2020.5.22.11.42}
f^*_{v} ( \tilde{\delta}_{k,\alpha})
-f^*_{0} ( \tilde{\delta}_{k,\alpha} ) =0
\end{equation}
for any relative vector field $v \in \Theta_{M/T}$.
Here $f_v \colon M[\epsilon] \rightarrow M$ is induced by $v \in \Theta_{M/T}$.
Remark that $f_0 \colon M[\epsilon] \rightarrow M$ is the trivial projection.
Note that $T$ is the base space of the family
$(\cC, \tilde{\bt}) \rightarrow T$ of 
$n$-pointed curves.
We say such coordinates $\{ z_{\alpha}\}_{\alpha}$ are \textit{refined}.
Now we assume that $\{ (\delta_k)_{\alpha}\}_{\alpha \in  I_{D(\tilde{\boldsymbol{p}})}}$ 
are described by using the refined coordinates.
\end{Def}

Remark that 
we avoid the refined coordinate to describe the representative 
$(d_k)_{\alpha\beta}$. 
We describe $(d_k)_{\alpha\beta}$ as 
$(d_k)_{\alpha\beta}= \tilde{d}_{k,\alpha\beta} \frac{\partial}{\partial f_{\beta}}$
(where $\tilde{d}_{k,\alpha\beta} \in   \mathcal{O}_{\cC_M}|_{U_{\alpha\beta}}$)
so that 
\begin{equation}\label{2020.6.3.10.57}
F^*_{U_{\alpha\beta}, v} ( \tilde{d}_{k,\alpha\beta} )-
F^*_{U_{\alpha\beta}, 0} ( \tilde{d}_{k,\alpha\beta} )=0
\end{equation}
for any relative vector field $v \in \Theta_{M/T}$.
Here $F_{U_{\alpha\beta}, v} \colon U_{\alpha\beta}  \times_M M[\epsilon]
\rightarrow U_{\alpha\beta}$
is defined by 
the projection $U_{\alpha\beta} \rightarrow M$
and $f_v \colon M[\epsilon] \rightarrow M$ for $v \in \Theta_{M/T}$.
For the representative 
$\{(\hat{d}_k)_{\alpha\beta}\}_{\alpha\beta}$ of $\hat{\mu}_k$,
we take a lift
$$
[\{(\hat{d}_k)_{\alpha\beta}\}_{\alpha\beta}] 
\in H^1(\cC_M, \Theta_{\cC_M/M}( - 2D(\tilde{\bt})- 2 D( \tilde{\boldsymbol{p}} ))) .
$$
We also denote by $\hat{\mu}_k$ this lift
of $\hat{\mu}_k
\in H^1(\cC_M, \Theta_{\cC_M/M}( - D(\tilde{\bt})-  D( \tilde{\boldsymbol{p}} )))$.

\begin{Rem}
In Definition \ref{Definition of Hamil} below, we will define Hamiltonians.
In this definition,
we will couple $(\hat{d}_k)_{\alpha\beta}$ on the square of 
$\tilde{\nabla}-\tilde{\nabla}_0^{\sigma_M}$ (see \eqref{2020.5.22.10.45} below).
This square of 
$\tilde{\nabla}-\tilde{\nabla}_0^{\sigma_M}$
has poles of order 2 at $D(\tilde{\bt})+ D( \tilde{\boldsymbol{p}} )$.
Then the Hamiltonians will depend on the choice of a lift in 
$H^1(\cC_M, \Theta_{\cC_M/M}( - 2D(\tilde{\bt})- 2 D( \tilde{\boldsymbol{p}} )))$
of $\hat{\mu}_k
\in H^1(\cC_M, \Theta_{\cC_M/M}( - D(\tilde{\bt})-  D( \tilde{\boldsymbol{p}} )))$.
Now we take a lift of $\hat{\mu}_k
\in H^1(\cC_M, \Theta_{\cC_M/M}( - D(\tilde{\bt})-  D( \tilde{\boldsymbol{p}} )))$
so that the equality \eqref{2020.11.7.17.11} (below) satisfies.
\end{Rem}

We take a trivialization 
$\bar{\phi}_{\alpha} \colon \tilde{E}|_{U_{\alpha}}\ra \cO^{\oplus r}_{U_{\alpha}}$.
We assume that 
the trivializations $(U_{\alpha}, \bar{\phi}_{\alpha})_{\alpha}$
are compatible with $\sigma_M$ (see Definition \ref{2020.11.7.17.19}).
Let
$\tilde{A}^0_{\alpha}  z_{\alpha}^{-1} dz_{\alpha}$ 
and $\tilde{A}_{\alpha}  f_{\alpha}^{-1} df_{\alpha}$ 
be connection matrices of $\tilde{\nabla}_0^{\sigma_M}$ 
and $\tilde{\nabla}$ on $U_{\alpha}$ via the trivialization $\bar{\phi}_{\alpha}$, respectively.
Set 
\begin{equation}\label{2020.5.22.10.45}
\begin{aligned}
(H_k)_{\alpha\beta} :=&\,
\bar{\phi}^{-1}_{\beta} \circ 
\left\langle (\hat{d}_k)_{\alpha\beta}, 
\left( \tilde{A}_{\beta}  \frac{df_{\beta}}{f_{\beta}} -
\tilde{A}^0_{\beta}  \frac{dz_{\beta}}{z_{\beta}} \right)^2 \right\rangle 
\circ \bar{\phi}_{\beta} 
\end{aligned}
\end{equation}
which is an element 
of $(\mathcal{E} nd (\tilde{E}) \otimes \Omega^1_{\cC_{M}/M})(U_{\alpha\beta})$
since 
$\sharp\{ i \mid \hat{t}_i |_{\cC_M} \cap U_{\alpha} \neq \emptyset \} \le 1$ for any $\alpha$ and 
$\sharp\{ \alpha \mid \hat{t}_i |_{\cC_M} \cap U_{\alpha} \neq \emptyset \} \le 1$ for any $i$.
Here $\{ \hat{t}_i \}$ is the set of the supports 
of the Cartier divisor $D(\tilde{\bt})+D( \tilde{\boldsymbol{p}} )$.

\begin{Def}\label{Definition of Hamil}
We define \textit{Hamiltonian functions} $H_k$ ($k=1,\ldots , 3g-3+n$) on $M$ as
\begin{equation*}
H_k = \frac{1}{2} [\{ \Tr ((H_k)_{\alpha\beta}) \}]
\in H^1(\Omega^1_{\cC_{M}/M}) \cong
H^0(\mathcal{O}_M)
\end{equation*}
for the lifts $\hat{\mu}_1 ,\ldots , \hat{\mu}_{3g-3+n}
\in H^1(\cC_M, \Theta_{\cC_M/M}( - 2D(\tilde{\bt})-  2D( \tilde{\boldsymbol{p}} )))$.
\end{Def}

\subsection{Calculation of the Hamiltonians}

By Proposition \ref{prop algebra split 2nd}, 
for the lifts $\hat{\mu}_1 ,\ldots , \hat{\mu}_{3g-3+n}$,
we have the vector fields $v_{\hat{\mu}_k}$ on $M$:
\begin{equation*}
v_{\hat{\mu}_k} :=
\left[ \left( \left\{ 
u_{\alpha\beta}^{\hat{\mu}_k\nabla_0}
 \right\}, \left\{  
-v_{\alpha}^{\hat{\mu}_k\nabla_0}
\right\} \right) \right] \in \bH^1(\cG^{\bullet}_M).
\end{equation*}

Let $v'=[ (\{ u'_{\alpha\beta}\} , \{v_{\alpha}' \} ) ]$ be an element of $\bH^1(\cG_M^{\bullet})$.
We put $v'_{\eta}:=[ ( \{\eta( u'_{\alpha\beta})\} , \{v_{\alpha}' \} ) ] \in \bH^1(\cF_M^{\bullet})$.
Let $f_{v'_{\eta}} \colon M[\epsilon] \rightarrow M$
be the morphism induced by $v'_{\eta}\in \bH^1(\cF_M^{\bullet})$.
Remark that $\bH^1(\cF_M^{\bullet}) \cong
 \Theta_{M_{\cC/T}^{\balpha}(\tilde{\bt} ,r,d)_{\boldsymbol{\nu}}/T}(M)$.
That is, the family $\cC_M$ of curves is constant along the direction 
of the vector field $v'_{\eta}$.
We denote by 
$$F_{v'_{\eta}} \colon \cC_M \times_M M[\epsilon]
\longrightarrow \cC_M$$
the morphism induced by $f_{v'_{\eta}} \colon M[\epsilon] \rightarrow M$.
We define an infinitesimal deformation of $\tilde{\nabla}_0^{\sigma_M}$ 
on $F_{v'_{\eta}}^* \tilde{E}$ 
by taking the pull-back of $\tilde{\nabla}_0^{\sigma_M}$
by $F_{v'_{\eta}} \colon\cC_M \times_M M[\epsilon]
\rightarrow \cC_M$:
\begin{equation*}
F_{v'_{\eta}}^* \tilde{\nabla}_0^{\sigma_M} \colon F_{v'_{\eta}}^* \tilde{E}
\longrightarrow \tilde{F}_{v'_{\eta}}^* \tilde{E} \otimes
 \Omega^1_{\cC_M \times_M M[\epsilon] /M[\epsilon]}
 (F_{v'_{\eta}}^*D(m) + F_{v'_{\eta}}^*D(\sigma_M)).
\end{equation*}
Let 
\begin{equation*}
\varphi_{\alpha} (v'_{\eta}) \colon  
F_{v'_{\eta}}^* \tilde{E}|_{U_{\alpha}\times \Spec \mathcal{O}_M[\epsilon] }
\xrightarrow[\sim]{\phi_{\alpha}(v'_{\eta})} 
\cO^{\oplus r}_{U_{\alpha}\times \Spec \mathcal{O}_M[\epsilon]} 
\xrightarrow{\bar{\phi}_{\alpha}^{-1}\otimes 1} \tilde{E}|_{U_{\alpha}} 
\otimes \cO_M[\epsilon]
\end{equation*}
be an isomorphism as in Section \ref{Inf defor fixed curve}. 
We put 
\begin{equation}\label{2020.5.22.10.50}
(v_{\alpha}^{0})' := (\varphi_{\alpha} (v'_{\eta}) \otimes \mathrm{id}) \circ   
(F_{v'_{\eta}}^* \tilde{\nabla}_0^{\sigma_M} )  
\circ (\varphi_{\alpha}(v'_{\eta}))^{-1} -  \tilde{\nabla}_0^{\sigma_M},
\end{equation}
which is an element of 
$\epsilon \otimes  \cE nd(\tilde{E}) \otimes \Omega^1_{\cC_M/M}(D(m) + D(\sigma_M))$.
The collection $(\{\eta( u'_{\alpha\beta})\} , \{ (v_{\alpha}^{0})' \} )$ satisfies
\begin{equation}\label{cocycle cond. nabla0 1}
\tilde{\nabla}_0^{\sigma_M}\circ (\eta( u'_{\alpha\beta})) 
- (\eta( u'_{\alpha\beta})) \circ \tilde{\nabla}_0^{\sigma_M} 
= (v_{\beta}^{0})' - (v_{\alpha}^{0})'
\end{equation}
as with (\ref{2019.8.6.18.40}).

\begin{Lem}\label{2020.5.24.14.16}
We take any $v'=[ (\{ u'_{\alpha\beta}\} , \{v_{\alpha}' \} ) ]\in \bH^1(\cG_M^{\bullet})$.
For $v'$, we define $(v_{\beta}^0)'$ by \eqref{2020.5.22.10.50}.
Then we have the following equality:
\begin{equation}\label{2020.5.19.17.57}
\omega(v_{\hat{\mu}_k}, [ (\{ \eta( u'_{\alpha\beta} )\} , \{v_{\alpha}' \}) ] ) 
= \left[ \left( \{ 0 \}_{\alpha\beta\gamma},
\left\{ \Tr ( \left( -\eta ( u_{\alpha\beta}^{\hat{\mu}_k\nabla_0} ) \right) 
\circ ( v_{\beta}' - (v_{\beta}^0)')) \right\}_{\alpha\beta} \right) \right]
\in \bH^2 (\Omega^{\bullet}_{\cC_M/M}).
\end{equation}
\end{Lem}

\begin{proof}
We consider the 1-form $\omega_M (v_{\hat{\mu}_k}, \cdot) 
\in \cH om_{\cO_M}( \bH^1(\cG_M^{\bullet}),\cO_M)$.
For $v'$, we compute 
\begin{equation*}
\omega_M(v_{\hat{\mu}_k} ,v' )=
[(\{  \Tr( \eta ( u_{\alpha\beta}^{\hat{\mu}_k\nabla_0} )
 \circ \eta( u'_{\beta\gamma})) \}, 
-\{ \Tr (\eta ( u_{\alpha\beta}^{\hat{\mu}_k\nabla_0} ) \circ  v_{\beta}' ) 
- \Tr (-v_{\alpha}^{\hat{\mu}_k\nabla_0} \circ \eta( u'_{\alpha\beta}) ) \} )],
\end{equation*}
which is a class of $\bH^2(\cC_M, \Omega_{\cC_M/M}^{\bullet})$.
We have
\begin{equation}\label{equation 2 in proof of Theorem}
\begin{aligned}
&- \Tr ( \eta ( u_{\alpha\beta}^{\hat{\mu}_k\nabla_0} ) \circ  v_{\beta}' ) 
+ \Tr (-v_{\alpha}^{\hat{\mu}_k\nabla_0} \circ \eta( u'_{\alpha\beta})) \\
&=  \Tr ( \left( -\eta ( u_{\alpha\beta}^{\hat{\mu}_k\nabla_0} ) \right) 
\circ ( v_{\beta}'- (v_{\beta}^0)'))
 + \left(  \Tr ( - \eta ( u_{\alpha\beta}^{\hat{\mu}_k\nabla_0} )\circ (v_{\beta}^0)' )
 - \Tr ( v_{\alpha}^{\hat{\mu}_k\nabla_0}  \circ \eta( u'_{\alpha\beta})) \right).
\end{aligned}
\end{equation}
We consider the second term of the right hand side of (\ref{equation 2 in proof of Theorem}).
We set
\begin{equation}\label{2020.5.20.23.24}
\eta_0(s) := s - \iota(\tilde{\nabla}^{\sigma_M}_0) \circ \mathrm{symb}_1(s).
\end{equation}
As in the proof of Lemma \ref{2020.5.24.14.29}, we may check the equality
\begin{equation*}
\tilde{\nabla}_0^{\sigma_M}\circ u_{\alpha\beta}^{\hat{\mu}_k\nabla_0} 
- u_{\alpha\beta}^{\hat{\mu}_k\nabla_0}  \circ \tilde{\nabla}_0^{\sigma_M}
= \tilde{\nabla}_0^{\sigma_M}\circ \eta_0(u_{\alpha\beta}^{\hat{\mu}_k\nabla_0} ) 
- \eta_0(u_{\alpha\beta}^{\hat{\mu}_k\nabla_0} ) \circ \tilde{\nabla}_0^{\sigma_M}.
\end{equation*}
Since $\eta_0(u_{\alpha\beta}^{\hat{\mu}_k\nabla_0} )=0$,
we have $\tilde{\nabla}_0^{\sigma_M}\circ u_{\alpha\beta}^{\hat{\mu}_k\nabla_0} 
- u_{\alpha\beta}^{\hat{\mu}_k\nabla_0}  \circ \tilde{\nabla}_0^{\sigma_M}=0$.
By this equality
and the equality \eqref{2019.8.6.18.40},
we have the following equalities:
\begin{equation}\label{cocycle cond. nabla0 2}
\begin{aligned}
\tilde{\nabla}_0^{\sigma_M}\circ (-\eta(u_{\alpha\beta}^{\hat{\mu}_k\nabla_0}) ) 
- (-\eta(u_{\alpha\beta}^{\hat{\mu}_k\nabla_0}) ) \circ \tilde{\nabla}_0^{\sigma_M} 
&= \tilde{\nabla}_0^{\sigma_M} \circ 
( \hat{u}_{\alpha\beta}^{\text{IMD}}- u_{\alpha\beta}^{\hat{\mu}_k\nabla_0}  )
 -( \hat{u}_{\alpha\beta}^{\text{IMD}}- u_{\alpha\beta}^{\hat{\mu}_k\nabla_0}  )
  \circ \tilde{\nabla}_0^{\sigma_M}\\ 
&= \tilde{\nabla}_0^{\sigma_M} \circ \hat{u}_{\alpha\beta}^{\text{IMD}}
 -\hat{u}_{\alpha\beta}^{\text{IMD}} \circ \tilde{\nabla}_0^{\sigma_M}\\
&= v_{\beta}^{\hat{\mu}_k\nabla_0} - v_{\alpha}^{\tilde{\mu}_k\nabla_0}.
\end{aligned}
\end{equation}
Set 
$\omega^0_{\alpha \beta}:=
 \Tr ( - \eta ( u_{\alpha\beta}^{\hat{\mu}_k\nabla_0} )\circ (v_{\beta}^0)' )
 - \Tr ( v_{\alpha}^{\hat{\mu}_k\nabla_0}  \circ \eta( u'_{\alpha\beta})).$
By the cocycle conditions of $\eta(u_{\alpha\beta}^{\hat{\mu}_k\nabla_0})$ and 
$\eta(u'_{\alpha\beta})$,
and by the equalities (\ref{cocycle cond. nabla0 1}) and (\ref{cocycle cond. nabla0 2}), 
we have the following equalities:
\begin{equation*}
\begin{aligned}
\omega^0_{\beta\gamma}-\omega^0_{\alpha \gamma}+\omega^0_{\alpha \beta}
&= \Tr ( - \eta ( u_{\alpha\beta}^{\hat{\mu}_k\nabla_0} )\circ 
 ((v_{\beta}^0)'- (v_{\gamma}^0)') )
 + \Tr ( (v_{\alpha}^{\hat{\mu}_k\nabla_0} -v_{\beta}^{\hat{\mu}_k\nabla_0} )
 \circ \eta( u'_{\beta\gamma})) \\
 &=  - \Tr \left( - \eta ( u_{\alpha\beta}^{\hat{\mu}_k\nabla_0} )\circ 
 \left( \tilde{\nabla}_0^{\sigma_M}\circ (\eta( u'_{\beta\gamma})) 
- (\eta( u'_{\beta\gamma})) \circ \tilde{\nabla}_0^{\sigma_M} \right)   \right) \\
&\qquad
 + \Tr \left( \left( \tilde{\nabla}_0^{\sigma_M}\circ (\eta(u_{\alpha\beta}^{\hat{\mu}_k\nabla_0}) ) 
- (\eta(u_{\alpha\beta}^{\hat{\mu}_k\nabla_0}) ) \circ \tilde{\nabla}_0^{\sigma_M}  \right)
 \circ \eta( u'_{\beta\gamma}) \right).
 \end{aligned}
\end{equation*}
Let $\tilde{A}^0_{\beta} \frac{d z_{\beta}}{z_{\beta}}$
be the connection matrices of $\tilde{\nabla}_0^{\sigma_M}$
associated to the trivialization $\bar{\phi}_{\beta}$.
That is $\tilde{\nabla}_0^{\sigma_M}|_{U_{\beta}} =
(\bar{\phi}_{\beta}^{-1} \otimes \mathrm{id})  
\circ (d+\tilde{A}^0_{\beta}  \frac{d z_{\beta}}{z_{\beta}})\circ \bar{\phi}_{\beta} $.
By this local description of $\tilde{\nabla}_0^{\sigma_M}$, we have the following equalities:
\begin{equation}\label{2020.11.3.17.41}
\begin{aligned}
 \omega^0_{\beta\gamma}-\omega^0_{\alpha \gamma}+\omega^0_{\alpha \beta}
 &=  \Tr \left(  \eta ( u_{\alpha\beta}^{\hat{\mu}_k\nabla_0} )\circ 
 \left( d (\eta( u'_{\beta\gamma}))+
  [  (\bar{\phi}_{\beta}^{-1} \otimes \mathrm{id}) \circ \tilde{A}^0_{\beta}
  \frac{d z_{\beta}}{z_{\beta}} \circ \bar{\phi}_{\beta} ,
   \eta( u'_{\beta\gamma}) ] \right)   \right) \\
&\qquad
 + \Tr \left( \left( d (\eta(u_{\alpha\beta}^{\hat{\mu}_k\nabla_0}) ) + 
 [  (\bar{\phi}_{\beta}^{-1} \otimes \mathrm{id}) \circ  \tilde{A}^0_{\beta}
 \frac{d z_{\beta}}{z_{\beta}} \circ \bar{\phi}_{\beta} ,
  \eta(u_{\alpha\beta}^{\hat{\mu}_k\nabla_0})  ] \right)
 \circ \eta( u'_{\beta\gamma}) \right)\\
 &= -d  \left( \Tr( \eta( -u_{\alpha\beta}^{\hat{\mu}_k\nabla_0}) \circ \eta( u_{\beta\gamma}'))
 \right).
 \end{aligned}
\end{equation}
Then 
the pair $(- \{ \Tr( \eta( -u_{\alpha\beta}^{\hat{\mu}_k\nabla_0}) \circ \eta( u_{\beta\gamma}')) \} , 
\{ \omega_{\alpha\beta}^0 \} )$
gives a class
\begin{equation}\label{nabla0 var symp}
[ (- \{ \Tr( \eta( -u_{\alpha\beta}^{\hat{\mu}_k\nabla_0}) \circ \eta( u_{\beta\gamma}')) 
\}_{\alpha\beta\gamma} , 
\{ \omega_{\alpha\beta}^0 \}_{\alpha\beta} )]
\end{equation}
of $\bH^2(\cC_M, \Omega_{\cC_M/M}^{\bullet})$.
If the class \eqref{nabla0 var symp} vanishes, 
then we obtain the equality \eqref{2020.5.19.17.57}
by the equality \eqref{equation 2 in proof of Theorem}.

We claim that the class (\ref{nabla0 var symp}) 
of $\bH^2(\cC_M, \Omega_{\cC_M/M}^{\bullet})$ vanishes.
To show this claim,
we consider replacement of 
$(\{ \eta( -u_{\alpha\beta}^{\hat{\mu}_k\nabla_0}) \} , 
\{ v_{\alpha}^{\hat{\mu}_k\nabla_0} \} )$ 
in the class \eqref{nabla0 var symp}.
We recall the definition of $v_{\alpha}^{\hat{\mu}_k\nabla_0}$.
We set $E_{\epsilon}:= (\mathrm{IM}_{M}^{\mu_k})^*\tilde{E}$ 
and $\nabla_{0,\epsilon}^{\sigma_M,\mathrm{IMD}}:=
(\mathrm{IM}_{M}^{\mu_k})^* \tilde{\nabla}_0^{\sigma_M}$.
If we take trivializations of $E_{\epsilon}$ which satisfies \eqref{2020.5.14.23.06},
then we can define 
 $v_{\alpha}^{\hat{\mu}_k\nabla_0}$
by \eqref{2019.7.31.11.14}.
The tuple $\{v_{\alpha}^{\hat{\mu}_k\nabla_0}\}_{\alpha}$
satisfies the equality \eqref{2019.8.6.18.40}.
We apply Lemma \ref{2020.5.15.11.56}
to $\mathrm{IM}_{M}^{\mu_k} \colon 
\cC_{\epsilon}^{k} \rightarrow \cC_M$ and the lift $\hat{\mu}_k
\in H^1(\cC_M, \Theta_{\cC_M/M}( - D(\tilde{\bt})-  D( \tilde{\boldsymbol{p}} )))$
as follows.
Let $\varphi'_{\alpha}$ be the composition \eqref{2020.5.15.12.53}.
By the definition \eqref{2019.7.31.11.14} and the equality \eqref{2020.5.16.13.00}, 
we have
\begin{equation*}
\begin{aligned}
 v_{\alpha}^{\hat{\mu}_k\nabla_0} &= (\varphi_{\alpha}^{-1} \otimes \mathrm{id}) 
 \circ (  \nabla_{0,\epsilon}^{\sigma_M,\text{IMD}} ) \circ \varphi_{\alpha}
-  \tilde{\nabla}_0^{\sigma_M} \\
&= (\varphi_{\alpha}^{-1} \otimes \mathrm{id}) 
\circ 
( (\varphi'_{\alpha})^{-1} \otimes \mathrm{id} )^{-1}
 \circ
( (\varphi'_{\alpha})^{-1} \otimes \mathrm{id} ) 
\circ 
  (  \nabla_{0,\epsilon}^{\sigma_M,\text{IMD}} ) \circ
  \varphi'_{\alpha}
  \circ (\varphi'_{\alpha})^{-1}
  \circ \varphi_{\alpha}
-  \tilde{\nabla}_0^{\sigma_M}\\
&= (\varphi_{\alpha}^{-1} \otimes \mathrm{id}) 
\circ 
( (\varphi'_{\alpha})^{-1} \otimes \mathrm{id} )^{-1}
 \circ
 \tilde{\nabla}_0^{\sigma_M}
  \circ (\varphi'_{\alpha})^{-1}
  \circ \varphi_{\alpha}
-  \tilde{\nabla}_0^{\sigma_M}.
\end{aligned}
\end{equation*}
We define $\chi_{\alpha} \in \mathcal{A}_{\tilde{E}}(U_{\alpha})$ (for any $\alpha$)
as
$(\varphi'_{\alpha})^{-1}
  \circ \varphi_{\alpha}= \mathrm{id}+\epsilon \chi_{\alpha}$. 
  Remark that 
  $ \varphi_{\alpha}^{-1}\circ
  \varphi'_{\alpha}= \mathrm{id}-\epsilon \chi_{\alpha}$.
Then we have the following equality:
\begin{equation}\label{2020.6.2.21.42}
\begin{aligned}
 v_{\alpha}^{\hat{\mu}_k\nabla_0} &= (\varphi_{\alpha}^{-1} \otimes \mathrm{id}) 
\circ 
( (\varphi'_{\alpha})^{-1} \otimes \mathrm{id} )^{-1}
 \circ
 \tilde{\nabla}_0^{\sigma_M}
  \circ (\varphi'_{\alpha})^{-1}
  \circ \varphi_{\alpha}
-  \tilde{\nabla}_0^{\sigma_M} \\
&= ( ((\mathrm{id}-\epsilon \chi_{\alpha}) \otimes \mathrm{id} )
 \circ
 \tilde{\nabla}_0^{\sigma_M}
  \circ (\mathrm{id}+\epsilon \chi_{\alpha})
-  \tilde{\nabla}_0^{\sigma_M}\\
&= \tilde{\nabla}_0^{\sigma_M} \circ \chi_{\alpha}
- \chi_{\alpha} \circ \tilde{\nabla}_0^{\sigma_M} 
\end{aligned}
\end{equation}
for any $\alpha$.
The symbol $\mathrm{symb}_1(\chi_{\alpha})$ 
of $\chi_{\alpha}\in  \mathcal{A}_{\tilde{E}}(U_{\alpha})$
has zero 
along the components of $D(\tilde{\boldsymbol{p}})$,
since the class
 $[\{ \mathrm{symb}_1 (\varphi_{\alpha}^{-1}\circ \varphi_{\beta} -\mathrm{id})\}_{\alpha\beta}]$
coincides with the class
$[\{\mathrm{symb}_1 ((\varphi'_{\alpha})^{-1}\circ \varphi'_{\beta} -\mathrm{id})\}_{\alpha\beta}]$
in $H^1(\cC_M, \Theta_{\cC_M/M}( - D(\tilde{\bt})-  D( \tilde{\boldsymbol{p}} )))$.
Then $\eta_0(\chi_{\alpha})$ has no pole along the components of $D(\tilde{\boldsymbol{p}})$.
As in the proof of Lemma \ref{2020.5.24.14.29}, we may check the equality
\begin{equation}\label{2020.11.3.17.35}
\begin{aligned}
\tilde{\nabla}_0^{\sigma_M} \circ  \eta_0(\chi_{\alpha})
- \eta_0( \chi_{\alpha}  )\circ \tilde{\nabla}_0^{\sigma_M} 
&=  \tilde{\nabla}_0^{\sigma_M} \circ \chi_{\alpha}
- \chi_{\alpha} \circ \tilde{\nabla}_0^{\sigma_M} \\
&=v_{\alpha}^{\hat{\mu}_k\nabla_0}
\end{aligned}
\end{equation}
by the equality \eqref{2020.6.2.21.42}.
Then we may replace $(\{ \eta( -u_{\alpha\beta}^{\hat{\mu}_k\nabla_0}) \} , 
\{ v_{\alpha}^{\hat{\mu}_k\nabla_0} \} )$
for 
\begin{equation}\label{2020.5.15.23.02}
(\{ \eta( -u_{\alpha\beta}^{\hat{\mu}_k\nabla_0}) -\eta_0(\chi_{\beta})
 + \eta_0(\chi_{\alpha}) \}_{\alpha\beta} , 
\{ 0 \}_{\alpha} ).
\end{equation}
We set 
$$
x_{\alpha\beta}^{\nabla_0}:=\eta( -u_{\alpha\beta}^{\hat{\mu}_k\nabla_0}) -\eta_0(\chi_{\beta})
 + \eta_0(\chi_{\alpha}).
$$
Second, we consider the pair
 $(\{ \eta( u_{\alpha\beta}') \} , 
\{ (v_{\alpha}^0)' \} )$.
We take a lift 
\begin{equation*}
\hat{\mu}(0) 
\in H^1(\cC_M, \Theta_{\cC_M/M}( - D(\tilde{\bt})- D( \tilde{\boldsymbol{p}} )))
\end{equation*}
of $0 \in H^1 (\cC_M, \Theta_{\cC_M/M}(-D(\tilde{\bt})))$ 
by a flat family 
$( \cC_M \times_M M [\epsilon],
\mathrm{Supp}(F_{v'_{\eta}}^*(D(\tilde{\bt})+D( \tilde{\boldsymbol{p}} )) ))$
of curves with points as in the definition of the lift \eqref{mulift} of $\mu$.
Let 
$\tilde{F}_{v'_{\eta}} \colon \cC^{\hat{\mu}(0)}_{v'_{\eta}} \rightarrow \cC_M$
be the morphism corresponding a representative of $\hat{\mu}(0)$.
We consider the pull-back of $(\tilde{E}, \tilde{\nabla}_0^{\sigma_M})$ 
by $\tilde{F}_{v'_{\eta}} \colon \cC^{\hat{\mu}(0)}_{v'_{\eta}} \rightarrow \cC_M$
and we apply Lemma \ref{2020.5.15.11.56} to
$\tilde{F}_{v'_{\eta}} \colon \cC^{\hat{\mu}(0)}_{v'_{\eta}} \rightarrow \cC_M$
and the lift $\hat{\mu}(0) $.
Then we can take $\chi_{\alpha}'\in\mathcal{A}_{\tilde{E}} (U_{\alpha})$ 
(for any $\alpha$) such that $\tilde{\nabla}_0^{\sigma_M} \circ \chi_{\alpha}'
- \chi_{\alpha}' \circ \tilde{\nabla}_0^{\sigma_M} 
=(v_{\alpha}^0)'$ for any $\alpha$.
As in the proof of Lemma \ref{2020.5.24.14.29}, we can check the following equality:
\begin{equation}\label{2020.5.20.23.36}
\tilde{\nabla}_0^{\sigma_M} \circ \eta_0(\chi_{\alpha}')
- \eta_0(\chi_{\alpha}') \circ \tilde{\nabla}_0^{\sigma_M} 
=(v_{\alpha}^0)'.
\end{equation}
Remark that $\eta_0(\chi_{\alpha}')$ may have 
poles on the components of $D(\tilde{\boldsymbol{p}})$,
since $\tilde{\nabla}^{\sigma_M}_0$ has 
poles on the components of $D(\tilde{\boldsymbol{p}})$.
Now we consider the class (\ref{nabla0 var symp}).
By \eqref{2020.5.15.23.02},  
the class \eqref{nabla0 var symp} coincides with
\begin{equation}\label{2020.5.20.23.42}
\left[ \left(- \left\{ 
\Tr \left( x_{\alpha\beta}^{\nabla_0}
\circ \eta( u_{\beta\gamma}') \right)   \right\}_{\alpha\beta\gamma} , 
\left\{  \Tr \left( x_{\alpha\beta}^{\nabla_0} 
\circ (v_{\beta}^0)'  \right)  \right\}_{\alpha\beta} 
\right) \right].
\end{equation}
We may check this coincidence as follows.
We take the image of cochain 
$(\{\Tr( \eta_0 (\chi_{\alpha}) \circ \eta(u'_{\alpha\beta}) )\}_{\alpha\beta} ,
 \{-\Tr( \eta_0(\chi_{\alpha}) \circ (v_{\alpha}^0)') ) \}_{\alpha})$:
\begin{equation*}{\small
\xymatrix{
  \{ \Tr( (\eta_0( \chi_{\beta}) - \eta_0( \chi_{\alpha}) )
 \circ \eta( u_{\beta\gamma}')\}_{\alpha\beta\gamma} 
 &   \\
\{\Tr( \eta_0 (\chi_{\alpha}) \circ \eta(u'_{\alpha\beta}) )\}_{\alpha\beta}
 \ar@{|->}[u]^-{\delta_1} \ar@{|->}[r]^-{d}
&
 \{ d\Tr( \eta_0 (\chi_{\alpha}) \circ \eta(u'_{\alpha\beta}) ) \}_{\alpha\beta} 
- \{ 
\Tr( \eta_0(\chi_{\beta}) (v_{\beta}^0)') )
-\Tr( \eta_0(\chi_{\alpha}) (v_{\alpha}^0)') )
 \}_{\alpha\beta} 
\\
& \{ - \Tr( \eta_0(\chi_{\alpha}) \circ (v_{\alpha}^0)') ) \}_{\alpha}
\ar@{|->}[u]^-{\delta_0}\ar@{|->}[r]^-d
& 0
 }}
\end{equation*}
Here 
\begin{equation*}
\begin{aligned}
\delta_0 \colon \prod_{\alpha} \Omega^1_{\cC_M/M} (U_{\alpha})
&\longrightarrow \prod_{\alpha\neq \beta} \Omega^1_{\cC_M/M} (U_{\alpha}\cap U_{\beta})
\text{ and} \\ 
\delta_1 \colon \prod_{\alpha\neq \beta} \mathcal{O}_{\cC_M} (U_{\alpha} \cap U_{\beta})
&\longrightarrow \prod_{\alpha\neq \beta\neq \gamma} \mathcal{O}_{\cC_M} 
(U_{\alpha}\cap U_{\beta}\cap U_{\gamma})
\end{aligned}
\end{equation*}
are the coboundary operators, and $d$ is the (relative) exterior derivative on
$\cC_M\rightarrow M$.
Then we have a coboundary.
We add this coboundary to the cocycle of \eqref{nabla0 var symp}:
\begin{equation*}
-\Tr( \eta( -u_{\alpha\beta}^{\hat{\mu}_k\nabla_0}) \circ \eta( u_{\beta\gamma}')) 
+\Tr( (\eta_0( \chi_{\beta}) - \eta_0( \chi_{\alpha}) )
 \circ \eta( u_{\beta\gamma}')=-\Tr \left( x_{\alpha\beta}^{\nabla_0}
\circ \eta( u_{\beta\gamma}') \right)
\end{equation*}
and
\begin{equation}\label{2020.11.3.17.42}
\begin{aligned}
&\omega_{\alpha\beta}^0
+d\Tr( \eta_0 (\chi_{\alpha}) \circ \eta(u'_{\alpha\beta}) ) 
-\Tr( \eta_0(\chi_{\beta})\circ (v_{\beta}^0)') )
+\Tr( \eta_0(\chi_{\alpha}) \circ (v_{\alpha}^0)') )\\
&= \Tr ( - \eta ( u_{\alpha\beta}^{\hat{\mu}_k\nabla_0} )\circ (v_{\beta}^0)' )
 - \Tr ( \left(\tilde{\nabla}_0^{\sigma_M} \circ  \eta_0(\chi_{\alpha})
- \eta_0( \chi_{\alpha}  )\circ \tilde{\nabla}_0^{\sigma_M} \right) \circ \eta( u'_{\alpha\beta}))\\
&\quad+d\Tr( \eta_0 (\chi_{\alpha}) \circ \eta(u'_{\alpha\beta}) ) 
-\Tr( \eta_0(\chi_{\beta})\circ (v_{\beta}^0)') ) \\
&\qquad 
+\Tr( \eta_0(\chi_{\alpha}) \circ 
\left( (v_{\beta}^{0})' - \tilde{\nabla}_0^{\sigma_M}\circ (\eta( u'_{\alpha\beta})) 
+ (\eta( u'_{\alpha\beta})) \circ \tilde{\nabla}_0^{\sigma_M} \right) ) \\
&=\Tr \left( x_{\alpha\beta}^{\nabla_0} 
\circ (v_{\beta}^0)'  \right).
\end{aligned}
\end{equation}
Here the first equality of \eqref{2020.11.3.17.42} follows from
 \eqref{cocycle cond. nabla0 1} and \eqref{2020.11.3.17.35}.
 The second equality of \eqref{2020.11.3.17.42}
 is checked as in \eqref{2020.11.3.17.41}.
 Then we have this coincidence 
 between \eqref{nabla0 var symp} and \eqref{2020.5.20.23.42}.
By the equality \eqref{2020.5.20.23.36} and the cocycle condition
$\tilde{\nabla}_0^{\sigma_M} \circ 
x_{\alpha\beta}^{\nabla_0}
-x_{\alpha\beta}^{\nabla_0}
\circ \tilde{\nabla}_0^{\sigma_M} =0$
of \eqref{2020.5.15.23.02},
we have 
\begin{equation}\label{2020.5.22.11.09}
d \left( \Tr\left(
x_{\alpha\beta}^{\nabla_0}
\circ \eta_0(\chi_{\beta}') \right) \right)  =
 \Tr \left( x_{\alpha\beta}^{\nabla_0}
\circ (v_{\beta}^0)'  \right) .
\end{equation}
Indeed, 
\begin{equation*}
\begin{aligned}
 \Tr \left(x_{\alpha\beta}^{\nabla_0}
\circ (v_{\beta}^0)'  \right)
&=
 \Tr \left( x_{\alpha\beta}^{\nabla_0}
\circ \left(\tilde{\nabla}_0^{\sigma_M} \circ \eta_0(\chi_{\beta}')
- \eta_0(\chi_{\beta}') \circ \tilde{\nabla}_0^{\sigma_M}  \right) \right) \\
&=
 \Tr \left(x_{\alpha\beta}^{\nabla_0}
\circ d \left( \eta_0(\chi_{\beta}'   ) \right) \right)  +
 \Tr \left( d\left( x_{\alpha\beta}^{\nabla_0}\right) 
\circ  \eta_0(\chi_{\beta}')  \right)\\
&\quad - \Tr \left( \left(
\tilde{\nabla}_0^{\sigma_M} \circ 
x_{\alpha\beta}^{\nabla_0}
-x_{\alpha\beta}^{\nabla_0}
\circ \tilde{\nabla}_0^{\sigma_M} \right) \circ  \eta_0(\chi_{\beta}')  \right)\\
&=d \left( \Tr\left(
x_{\alpha\beta}^{\nabla_0} 
\circ \eta_0(\chi_{\beta}') \right) \right) .
\end{aligned}
\end{equation*}
Note that the trace of $x_{\alpha\beta}^{\nabla_0} 
\circ \eta_0(\chi_{\beta}')$ 
has no pole on $U_{\alpha\beta}$.
That is, the trace of $x_{\alpha\beta}^{\nabla_0} 
\circ \eta_0(\chi_{\beta}')$ is an element 
of $\mathcal{O}_{\cC_M}(U_{\alpha\beta})$.
We may check that
$$
x_{\beta\gamma}^{\nabla_0} 
\circ \eta_0(\chi_{\gamma}')
-
x_{\alpha\gamma}^{\nabla_0} 
\circ \eta_0(\chi_{\gamma}')
+
x_{\alpha\beta}^{\nabla_0} 
\circ \eta_0(\chi_{\beta}')
=x_{\alpha\beta}^{\nabla_0} 
\circ \left( \eta_0(\chi_{\beta}') - \eta_0(\chi_{\gamma}' ) \right) 
$$
by the cocycle condition of 
$\{ x_{\alpha\beta}^{\nabla_0}  \}$.
By this equality and the equality \eqref{2020.5.22.11.09}, 
we may check the class \eqref{2020.5.20.23.42}
coincides with 
\begin{equation}\label{2020.11.7.16.32}
\left[ \left(- \left\{ 
\Tr \left(x_{\alpha\beta}^{\nabla_0} 
\circ \left(\eta( u_{\beta\gamma}') -\eta_0(\chi_{\beta}') + \eta_0(\chi_{\gamma}')
 \right) \right)   \right\}_{\alpha\beta\gamma} , 
\left\{  0  \right\}_{\alpha\beta} 
\right) \right] \in \bH^2(\cC_M, \Omega_{\cC_M/M}^{\bullet}).
\end{equation}
Remark that $d(\Tr \left(x_{\alpha\beta}^{\nabla_0} 
\circ \left(\eta( u_{\beta\gamma}') -\eta_0(\chi_{\beta}') + \eta_0(\chi_{\gamma}')
 \right) \right))=0$, by the cocycle condition of $\bH^2(\cC_M, \Omega_{\cC_M/M}^{\bullet})$.
 That is $\Tr \left(x_{\alpha\beta}^{\nabla_0} 
\circ \left(\eta( u_{\beta\gamma}') -\eta_0(\chi_{\beta}') + \eta_0(\chi_{\gamma}')
 \right) \right)$ is constant. 
Then we may consider this class in
 $\bH^2(\cC_M, \Omega_{\cC_M/M}^{\bullet})$
 as a class in a second cohomology of a locally constant sheaf.
 Since this locally constant sheaf is with the Zariski topology,
 we obtain that the class \eqref{2020.11.7.16.32} vanishes.
Finally, we have that the class \eqref{nabla0 var symp} vanishes. 
We obtain the equality \eqref{2020.5.19.17.57}
by the equality \eqref{equation 2 in proof of Theorem}.
\end{proof}

\begin{Thm}\label{main theorem 0}
\textit{
Let $\omega_M$ be the restriction of} (\ref{2019.7.28.14.00}) 
\textit{
on the affine open subset $M \subset M_{\cC/T}^{\balpha}(\tilde{\bt} ,r,d)_{\boldsymbol{\nu}}$.
We define $dH_k \in \cH om_{\cO_M}( \bH^1(\cG_M^{\bullet}),\cO_M)$ as 
\begin{equation*}
dH_k \colon \bH^1(\cG^{\bullet}_M) \lra \cO_M; \ [(\{ u_{\alpha\beta} \},\{ v_{\alpha} \} ) ] 
\longmapsto d_{M/T}H_k [ (\{ \eta(u_{\alpha\beta}) \},\{ v_{\alpha} \}  ) ].
\end{equation*}
Then the 1-form $\omega_M (v_{\hat{\mu}_k}, \cdot) 
\in \cH om_{\cO_M}( \bH^1(\cG_M^{\bullet}),\cO_M)$ coincides with the 1-form $dH_k$.
}
\end{Thm}

\begin{proof}
We claim that the 1-form
$dH_i \in \bH^1 (\cG_M^{\bullet})^{\vee}$ is described as
\begin{equation}\label{2020.5.19.16.40}
\begin{aligned}
d H_k \colon \bH^1 (\cG^{\bullet}_M)  &\lra 
H^1(\Omega^1_{\cC_M/M}) \cong 
 H^0(\cO_M ) \\
[(\{ u_{\alpha\beta}' \}, \{ v_{\alpha}' \})] &\longmapsto
\left[ \left\{ \mathrm{Tr} \left( \left(  \bar{\phi}_{\beta}^{-1} \circ
\left< (\hat{d}_k)_{\alpha\beta},
\left(   \tilde{A}_{\beta}\frac{d f_{\beta}}{f_{\beta}} 
-\tilde{A}_{\beta}^0\frac{d z_{\beta}}{z_{\beta}} \right) 
\right>\circ \bar{\phi}_{\beta} \right)
(v_{\beta}' -  (v_{\beta}^0)' )  \right)  \right\} \right] .
\end{aligned}
\end{equation}
We show this claim.
We consider the pull-back of 
$\tilde{\nabla} - \tilde{\nabla}_0^{\sigma_m} $
under the morphism 
$F_{v'_{\eta}} \colon \cC_M \times_M M[\epsilon] \rightarrow \cC_M$.
We have the following equalities
\begin{equation*}
\begin{aligned}
 \Tr \left(F_{v'_{\eta}}^*\tilde{\nabla} 
- F_{v'_{\eta}}^*\tilde{\nabla}_0^{\sigma_m} \right)^2 
&= \Tr \left( \varphi_{\beta}(v'_{\eta}) \circ \left(F_{v'_{\eta}}^*\tilde{\nabla} 
- F_{v'_{\eta}}^*\tilde{\nabla}_0^{\sigma_m} \right)^2
\circ \varphi_{\beta}(v'_{\eta})^{-1} \right) \\
&= \Tr \left( \left( \tilde{\nabla}  + \epsilon v_{\beta}'
- \tilde{\nabla}_0^{\sigma_m} - \epsilon (v_{\beta}^0)' \right)^2 \right)\\
&= \Tr \left( \tilde{\nabla}  - \tilde{\nabla}_0^{\sigma_m} \right)^2  
+ 2 \epsilon \Tr \left( \left( \tilde{\nabla}  - \tilde{\nabla}_0^{\sigma_m} \right)(v_{\beta}' -(v_{\beta}^0)')
 \right)\\
 &= \Tr \left( \tilde{\nabla}  - \tilde{\nabla}_0^{\sigma_m} \right)^2  
+ 2 \epsilon \Tr \left( \left( 
 \bar{\phi}_{\beta}^{-1} \circ
\left(   \tilde{A}_{\beta}\frac{d f_{\beta}}{f_{\beta}} 
-\tilde{A}_{\beta}^0\frac{d z_{\beta}}{z_{\beta}} \right) 
 \circ \bar{\phi}_{\beta} \right)
 (v_{\beta}' -(v_{\beta}^0)')
 \right).
\end{aligned}
\end{equation*}
Let $D_{v'_{\eta}} \colon \mathcal{O}_M \rightarrow \mathcal{O}_M$
be the derivative corresponding to 
the vector field $v'_{\eta}$.
We compute $D_{v'_{\eta}} \left( \frac{1}{2} [\{ \Tr ((H_k)_{\alpha\beta}) \}] \right) $.
By the equalities \eqref{2020.5.22.11.42} and \eqref{2020.6.3.10.57},
we have
\begin{equation}\label{2020.11.7.17.11}
\begin{aligned}
&D_{v'_{\eta}} \left( \frac{1}{2} [\{ \Tr ((H_k)_{\alpha\beta}) \}] \right) \\
&=\left[ \left\{  \mathrm{Tr} \left( \left(  \bar{\phi}_{\beta}^{-1} \circ
\left< (d_k)_{\alpha\beta},
\left(   \tilde{A}_{\beta}\frac{d f_{\beta}}{f_{\beta}} 
-\tilde{A}_{\beta}^0\frac{d z_{\beta}}{z_{\beta}} \right) 
\right>\circ \bar{\phi}_{\beta} \right)
(v_{\beta}' -  (v_{\beta}^0)' )  \right) \right\} \right]\\
&\quad +\left[ \left\{  \mathrm{Tr} \left( \left(  \bar{\phi}_{\beta}^{-1} \circ
\left< \left(  \tilde{\delta}^{(0)}_{k,\beta} \frac{\partial}{\partial z_{\beta}}
- \tilde{\delta}^{(0)}_{k,\alpha}  \frac{\partial}{\partial z_{\alpha}} \right),
\left(   \tilde{A}_{\beta}\frac{d f_{\beta}}{f_{\beta}} 
-\tilde{A}_{\beta}^0\frac{d z_{\beta}}{z_{\beta}} \right) 
\right>\circ \bar{\phi}_{\beta} \right)
(v_{\beta}' -  (v_{\beta}^0)' )  \right) \right\} \right]\\
&\qquad \in H^1(\Omega^1_{\cC_M/M}) .
\end{aligned}
\end{equation}
Since $D_{v'_{\eta}} \left( \frac{1}{2} [\{ \Tr ((H_k)_{\alpha\beta}) \}] \right)
=d_{M/T}H_k ([(\{ \eta(u_{\alpha\beta}') \}, \{ v_{\alpha}' \} ) ])$
and 
$$d_{M/T}H_k ([(\{ \eta(u_{\alpha\beta}') \}, \{ v_{\alpha}' \} ) ]) \\
= dH_k ([(\{ u_{\alpha\beta}' \}, \{ v_{\alpha}' \} ) ]),$$
we have the description \eqref{2020.5.19.16.40}.

By the definition \eqref{2019.7.31.11.05} of $u_{\alpha\beta}^{\hat{\mu}_k\nabla_0}$
for $\hat{\mu}=\hat{\mu}_k$ and the definition \eqref{2020.11.3.19.35}
of $\eta$, we have 
$$
\eta ( u_{\alpha\beta}^{\hat{\mu}_k\nabla_0} )
=
\bar{\phi}_{\beta}^{-1} \circ
\left< (\hat{d}_k)_{\alpha\beta},
\left(  \tilde{A}_{\beta}^0\frac{d z_{\beta}}{z_{\beta}} 
-\tilde{A}_{\beta}\frac{d f_{\beta}}{f_{\beta}} \right) 
\right>\circ \bar{\phi}_{\beta}.
$$
Here,
$\tilde{A}^0_{\alpha}  z_{\alpha}^{-1} dz_{\alpha}$ 
and $\tilde{A}_{\alpha}  f_{\alpha}^{-1} df_{\alpha}$ 
be connection matrices of $\tilde{\nabla}_0^{\sigma_M}$ 
and $\tilde{\nabla}$ on $U_{\alpha}$ via the trivialization $\bar{\phi}_{\alpha}$, respectively.
By the description \eqref{2020.5.19.16.40} of $dH_k$, 
we have the following equalities:
\begin{equation*}
\begin{aligned}
dH_k ([(\{ u_{\alpha\beta}' \}, \{ v_{\alpha}' \} ) ])
&=\left[ \left\{ \mathrm{Tr} \left( \left(  \bar{\phi}_{\beta}^{-1} \circ
\left< (\hat{d}_k)_{\alpha\beta},
\left(   \tilde{A}_{\beta}\frac{d f_{\beta}}{f_{\beta}} 
-\tilde{A}_{\beta}^0\frac{d z_{\beta}}{z_{\beta}} \right) 
\right>\circ \bar{\phi}_{\beta} \right)
(v_{\beta}' -  (v_{\beta}^0)' )  \right)  \right\} \right] \\
&=\left[ 
\left\{ \Tr ( \left( -\eta ( u_{\alpha\beta}^{\hat{\mu}_k\nabla_0} ) \right) 
\circ ( v_{\beta}' - (v_{\beta}^0)')) \right\}_{\alpha\beta}  \right] \in H^1(\Omega^1_{\cC_{M}/M}) .
\end{aligned}
\end{equation*}
We take the image of this element of $H^1(\Omega^1_{\cC_{M}/M})$
under the natural morphism $H^1(\Omega^1_{\cC_{M}/M})
\rightarrow \bH^2 (\Omega^{\bullet}_{\cC_M/M})$.
This image is
\begin{equation*}
 \left[ \left( \{ 0 \}_{\alpha\beta\gamma},
\left\{ \Tr ( \left( -\eta ( u_{\alpha\beta}^{\hat{\mu}_k\nabla_0} ) \right) 
\circ ( v_{\beta}' - (v_{\beta}^0)')) \right\}_{\alpha\beta} \right) \right] 
\in \bH^2 (\Omega^{\bullet}_{\cC_M/M}).
\end{equation*}
By Lemme \ref{2020.5.24.14.16}, this image coincides with
$\omega(v_{\hat{\mu}_k}, [ (\{ \eta( u'_{\alpha\beta} )\} , \{v_{\alpha}' \}) ] )$. 
For the natural morphism 
$H^1(\Omega^1_{\cC_{M}/M})
\rightarrow \bH^2 (\Omega^{\bullet}_{\cC_M/M})$,
the diagram
\begin{equation*}
\xymatrix{
H^1(\Omega^1_{\cC_{M}/M}) \ar[r] \ar[rd]^-{\cong}
& \bH^2 (\Omega^{\bullet}_{\cC_M/M}) \ar[d]^-{\cong} \\
&H^0(\cO_M ) 
}
\end{equation*}
is commutative.
Then the image of $dH_k ([(\{ u_{\alpha\beta}' \}, \{ v_{\alpha}' \} ) ])$
under the isomorphism 
$H^1(\Omega^1_{\cC_{M}/M}) \xrightarrow{\cong} H^0(\cO_M )$
coincides with the image of 
$\omega(v_{\hat{\mu}_k}, [ (\{ \eta( u'_{\alpha\beta} )\} , \{v_{\alpha}' \}) ] )$ 
under the isomorphism 
$\bH^2 (\Omega^{\bullet}_{\cC_M/M}) \xrightarrow{\cong} H^0(\cO_M )$.
This coincidence means the coincidence in the statement of this theorem.
\end{proof}

In the proof of Lemma \ref{2020.5.24.14.16}, 
we have the equality \eqref{cocycle cond. nabla0 2}.
Now we consider meaning of this equality 
by using the morphism $h_{\nabla_0}$.
We take $\mu_k$ and the lift $\hat{\mu}_k$ of $\mu_k$.
We take the vector field
 $$ [(\{ \hat{u}_{\alpha\beta}^{\text{IMD}}\} , 
 \{\hat{v}_{\alpha}^{\text{IMD}} \} )]$$
 of isomonodromic deformations in $M$ associated to $\mu_k$.
We take the image of $ [(\{ \hat{u}_{\alpha\beta}^{\text{IMD}}\} , 
 \{\hat{v}_{\alpha}^{\text{IMD}} \} )]$
under the morphism 
$\Theta_M \rightarrow h^*_{\nabla_0} \Theta_{M'}$.
Then the infinitesimal deformation of $(\tilde{E}, \tilde{\nabla}_0^{\sigma_M},
\boldsymbol{l}^{\sigma_m})$
 associated to this image corresponds to the pair
 \begin{equation}\label{2020.5.11.21.19}
 (\{ \hat{u}_{\alpha\beta}^{\text{IMD}}\}_{\alpha\beta} , 
 \{v_{\alpha}^{\tilde{\mu}_k\nabla_0} \}_{\alpha} ).
 \end{equation}
We consider the vector field of the isomonodromic deformation for 
$\varpi' \colon M' \rightarrow \mathcal{M}_{g,n+N}$.
That is, we define a splitting $(\varpi')^* \Theta_{\mathcal{M}_{g,n+N}} \rightarrow \Theta_{M'}$.
By $h_{\nabla_0}$, 
we have a morphism 
$h_{\nabla_0}^* (\varpi')^* \Theta_{\mathcal{M}_{g,n+N}} \rightarrow h_{\nabla_0}^* \Theta_{M'}$.
The lift $\hat{\mu}_k$ of $\mu_k$ gives 
a section of $h_{\nabla_0}^* (\varpi')^* \Theta_{\mathcal{M}_{g,n+N}}$.
We take the image of $\hat{\mu}_k$ by this
morphism $h_{\nabla_0}^* (\varpi')^* \Theta_{\mathcal{M}_{g,n+N}} 
\rightarrow h_{\nabla_0}^* \Theta_{M'}$.
This image gives
an infinitesimal deformation of $\tilde{\nabla}_0^{\sigma_M}$ 
parametrized by $M[\epsilon]$.
By this infinitesimal deformation of $(\tilde{E}, \tilde{\nabla}_0^{\sigma_M},
\boldsymbol{l}^{\sigma_m})$,
we have a pair 
\begin{equation}\label{2020.5.11.21.20}
( \{ u_{\alpha\beta}^{\hat{\mu}_k\nabla_0} \}_{\alpha\beta}, 
\{  0 \}_{\alpha} ).
\end{equation}
We consider the difference of \eqref{2020.5.11.21.19} and \eqref{2020.5.11.21.20},
which is described by 
\begin{equation}\label{2020.5.12.12.17}
(\{ \eta_0(\hat{u}_{\alpha\beta}^{\text{IMD}})\}_{\alpha\beta} , 
\{v_{\alpha}^{\tilde{\mu}_k\nabla_0} \}_{\alpha} ) .
\end{equation}
Here $\eta_0$ is defined in \eqref{2020.5.20.23.24}.
The pair \eqref{2020.5.12.12.17} means 
an infinitesimal deformation of $(\tilde{E}, \tilde{\nabla}_0^{\sigma_M},
\boldsymbol{l}^{\sigma_m})$.
The equality \eqref{cocycle cond. nabla0 2} means the
cocycle condition of 
$(\{ \eta_0(\hat{u}_{\alpha\beta}^{\text{IMD}})\} , \{v_{\alpha}^{\tilde{\mu}_k\nabla_0} \} )$.
Here note that 
$$-\eta ( u_{\alpha\beta}^{\hat{\mu}_k\nabla_0} )
=\eta_0(\hat{u}_{\alpha\beta}^{\text{IMD}}).
$$
That is, $(\{ \eta_0(\hat{u}_{\alpha\beta}^{\text{IMD}})\} , \{v_{\alpha}^{\tilde{\mu}_k\nabla_0} \} )$
gives an element of $h^*_{\nabla} \Theta_{M'}$,
which is described by using some hypercohomology.
The isomonodromic deformations gives the 
isomonodromic splitting 
$\Theta_{M'} \cong \Theta_{M'/\mathcal{M}_{g,n+N}} 
\oplus (\varpi')^* \Theta_{\mathcal{M}_{g,n+N}}$.
Roughly speaking,
the class of \eqref{2020.5.12.12.17} means
the first component of the decomposition of 
the class of \eqref{2020.5.11.21.19} induced by the 
isomonodromic splitting.

Next we consider meaning of \eqref{equation 2 in proof of Theorem}
and \eqref{nabla0 var symp}.
We may define a relative 2-form
\begin{equation*}
\omega_{M'/\mathcal{M}_{g,n+N}}^0\colon 
\Theta_{M'/\mathcal{M}_{g,n+N}} \otimes \Theta_{M'/\mathcal{M}_{g,n+N}} \longrightarrow 
\mathcal{O}_{M'}
\end{equation*}
on $\Theta_{M'/\mathcal{M}_{g,n+N}}$ as in Section \ref{S2 symp form}.
By the 
isomonodromic splitting $\Theta_{M'} 
\cong \Theta_{M'/\mathcal{M}_{g,n+N}} \oplus (\varpi')^* \Theta_{\mathcal{M}_{g,n+N}}$
and taking the first component of vector fields in $\Theta_{M'}$,
we may define a 2-form 
\begin{equation*}
\omega_{M'}^0\colon 
\Theta_{M'} \otimes \Theta_{M'} \longrightarrow 
\mathcal{O}_{M'}.
\end{equation*}
We consider the pull-back $h_{\nabla_0}^* \omega^0_{M'}$
 of the $2$-form $\omega^0_{M'}$
under the morphism $h_{\nabla_0} \colon M \rightarrow M'$.
We put  
$$
\Phi_0 := h_{\nabla_0}^* \omega^0_{M'}( [(\{ \hat{u}_{\alpha\beta}^{\text{IMD}}\} , 
 \{\hat{v}_{\alpha}^{\text{IMD}} \} )], \cdot),
$$ 
which is a 1-form on $M$.
 We denote by $\Phi_1$ the 1-form 
$\omega_M (v_{\hat{\mu}_k} , \cdot)$ on $M$.
Moreover, 
we denote by $\tilde{\Phi}_0$ the relative 1-form 
induced by $\Phi_0$ on $M\rightarrow T \times \{ \boldsymbol{\nu}\}$
and we denote by $\tilde{\Phi}_1$ the relative 1-form 
induced by $\Phi_1$ on $M\rightarrow T \times \{ \boldsymbol{\nu}\}$.
We consider the relative 1-form $d_{M/T}H_k$
on $\varpi_M\colon M \rightarrow T$.
By the isomonodromic splitting
we take the 2-forms 
$\tilde{\Phi}^{\mathrm{IMD}}_1$, $dH_k$, and $\tilde{\Phi}_0^{\mathrm{IMD}}$ on $M$
induced by $\tilde{\Phi}_1$, $d_{M/T}H_k$, and $\tilde{\Phi}_0$,
respectively. 
We may check that 
the class \eqref{nabla0 var symp} coincides with the value
$\tilde{\Phi}_0^{\mathrm{IMD}}([\{u_{\alpha\beta}'\} , \{v_{\alpha}' \}])$
by the argument as in the proof of Lemma \ref{2020.5.24.14.16}.
The equality \eqref{equation 2 in proof of Theorem} means the equality
\begin{equation*}
\tilde{\Phi}^{\mathrm{IMD}}_1=dH_k+\tilde{\Phi}_0^{\mathrm{IMD}}.
\end{equation*}

\subsection{Hamiltonian description of isomonodromic deformations}

Set $N=r^2(g-1)+nr(r-1)/2+1$.
Note that $\dim M =2N$.
Let $ \partial / \partial q_i \in \bold{H}^1 (\cF^{\bullet}_M)$ 
and $\partial / \partial p_i \in \bold{H}^1  (\cF^{\bullet}_M) $ 
($i=1,\ldots , N$) be vector fields on $M$
such that
the morphism 
\begin{equation*}\label{trivialization relative tangent}
\begin{aligned}
\cO_M^{\oplus 2N} &\lra  \bold{H}^1  (\cF^{\bullet}_M)\\ 
(f_1,\ldots, f_{2N}) &\longmapsto f_1\partial / \partial q_1 + \ldots+ f_{N}\partial / \partial q_{N}  
 +f_{N+1}\partial / \partial p_1 + \ldots+ f_{2N}\partial / \partial p_{N}
\end{aligned}
\end{equation*}
gives a trivialization of $\bold{H}^1 (\cF^{\bullet}_M)$
and the vector fields satisfy the conditions
$\omega_M(\partial / \partial q_i,\partial / \partial q_j)
= \omega_M(\partial / \partial p_i,\partial / \partial p_j)=0$ and 
$\omega_M(\partial / \partial q_i,\partial / \partial p_j) 
= \delta_{i,j}$, where $\delta_{i,j}$ is the Kronecker's symbol.
Here we also denote by $ \partial / \partial q_i$ and $ \partial / \partial p_i$ the images of
 $ \partial / \partial q_i \in \bold{R}^1 \pi_* (\cF^{\bullet})(M)$ 
 and $\partial / \partial p_i \in \bold{R}^1 \pi_* (\cF^{\bullet}_M) $ by the tangent morphism
\begin{equation*}
\bold{H}^1  (\cF^{\bullet}_M) \lra \bold{H}^1  (\cG^{\bullet}_M).
\end{equation*}

\begin{Rem}
We do not have a proof of existence of such algebraic vector fields 
$ \partial / \partial q_i \in \bold{H}^1  (\cG^{\bullet}_M)$ and 
$\partial / \partial p_i \in \bold{H}^1  (\cG^{\bullet}_M) $.
For some special cases, 
there exists studies of giving algebraic Darboux coordinates 
for some symplectic structures on moduli spaces of connections.
For example, one of techniques to give algebraic
 Darboux coordinates is the 
\textit{theory of apparent singularities} 
(\cite{Okamoto4}, \cite{Iwa1}, \cite{DM}, \cite{LS}, \cite{FL}).
In this paper, we do not discuss that there exist 
algebraic Darboux coordinates for the symplectic structure.
We assume that existence, although this argument is optimistic.
\end{Rem}

\begin{Cor}\label{main theorem}
\textit{
If we may take vector fields $ \partial / \partial q_i 
\in \bold{H}^1  (\cG^{\bullet}_M)$ and $\partial / \partial p_i 
\in \bold{H}^1  (\cG^{\bullet}_M) $ as above,
then the vector field determined by the isomonodromic deformation on $M$ is described as
\begin{equation}\label{Hamil flow of IMD}
v_{\hat{\mu}_k}
 - \sum_{i=1}^N \left( dH_k \left( \frac{\partial}{\partial p_i} \right)
 \frac{\partial}{\partial q_i} - dH_k \left( \frac{\partial}{\partial q_i} \right)\frac{\partial}{\partial p_i} \right)
\end{equation}
for $k=1,\ldots, 3g-3+n$.
}
\end{Cor}

\begin{proof}
Let $X$ be the vector field (\ref{Hamil flow of IMD}).
We show that $X \in \mathrm{Ker}(\omega_M)$.
Note that
\begin{equation*}
\begin{aligned}
 \omega_M\left(\partial/\partial q_i, \cdot \right) 
&\colon [(\{ u_{\alpha\beta}\}, \{ v_{\alpha}\})] 
\longmapsto  \omega\left(\partial/\partial q_i, [(\{ \eta(u_{\alpha\beta})\}, \{ v_{\alpha} \})] \right),\\
 \omega_M\left(\partial/\partial p_i, \cdot \right) 
&\colon [(\{ u_{\alpha\beta}\}, \{ v_{\alpha}\})] 
\longmapsto  \omega\left(\partial/\partial p_i, [(\{ \eta(u_{\alpha\beta})\}, \{  v_{\alpha} \})] \right), \\
 dH_k &\colon [(\{ u_{\alpha\beta}\}, \{ v_{\alpha}\})]
  \longmapsto d_{M/(T\times N^{(n)}_r(d))}H_k [(\{ \eta(u_{\alpha\beta})\}, \{  v_{\alpha} \})].
\end{aligned}
\end{equation*}
We have
\begin{equation*}
dH_k -\sum_i \left( dH_k \left( \frac{\partial}{\partial p_i} \right) 
\omega_M\left(\partial/\partial q_i, \cdot \right)
- dH_k \left( \frac{\partial}{\partial q_i} \right) 
\omega_M \left(\partial/\partial p_i, \cdot \right) \right)  =0.
\end{equation*}
By Theorem \ref{main theorem 0}, we have
\begin{equation*}
\omega_M\left(v_{\hat{\mu}_k}, \cdot \right) 
-\sum_i \left( dH_k \left( \frac{\partial}{\partial p_i} \right)
 \omega_M\left(\partial/\partial q_i, \cdot \right)
- dH_k \left( \frac{\partial}{\partial q_i} \right) 
\omega_M \left(\partial/\partial p_i, \cdot \right) \right) =0.
\end{equation*}
Then we obtain that $X \in \mathrm{Ker}(\omega_M)$.
Proposition \ref{Kernel is IMD} implies that $X$ is the vector field determined 
by the isomonodromic deformation.
By uniqueness of the isomonodromic deformation for a Kodaira--Spencer class, 
we obtain this corollary.
\end{proof}

\noindent
{\bf Acknowledgments.}
The author warmly thanks an anonymous referee for pointing out 
a gap in the proof of Proposition \ref{Prop d t s of m by c c} and
several other incorrect points.
Specially, Section 3 have been corrected essentially 
by the anonymous referee's comments.
He appreciates all referees' comment.
They help him to improve the quality of this manuscript.
He thanks Masa-Hiko Saito for guiding me to the problem treated in this paper 
and for warm encouragement.
He also thanks Michi-Aki Inaba for giving him valuable comments 
and for warm encouragement.
He is supported by JSPS KAKENHI Grant Numbers
17H06127, 18J00245 and 19K14506.

\end{document}